\documentclass[10pt]{amsart} % From Prof. Rathan
\title{\texttt{mcode.sty} Demo}
\usepackage[a4paper, width=18.2cm, height=26.2cm, centering]{geometry}

\usepackage{amsmath}
\usepackage{amssymb}
\usepackage{amsthm}
\usepackage{amsfonts}
\usepackage{latexsym}
\usepackage{fixmath}
\usepackage{bm}

\usepackage{esint}       % <-- provides \dashint

\usepackage{graphicx}
\usepackage{epsfig}
\usepackage{epstopdf}
\usepackage{graphics}
\usepackage{subcaption}
\usepackage{float}

\usepackage{siunitx}
\usepackage{booktabs}
\usepackage{multirow}
\usepackage{color}
\usepackage{xcolor}
\usepackage{pagecolor}
\usepackage{soul}

\usepackage[all]{xy}

\usepackage[vlined,ruled]{algorithm2e}

\usepackage{caption}
\usepackage{subcaption}
\usepackage{setspace}
\usepackage{ragged2e}

\usepackage{tikz}
\usetikzlibrary{arrows.meta}

\usepackage{hyperref}

\hypersetup{
  colorlinks=true,
  urlcolor=red,
  linkcolor=blue,
  citecolor=blue
}

\makeatletter
\newcommand{\dashint}{\mathop{\int\!\!\!\!\!-}}
\makeatother

\numberwithin{equation}{section}
\numberwithin{equation}{section}
\newcommand{\baa}{\begin{eqnarray*}}
\newcommand{\eaa}{\end{eqnarray*}}
\newcommand{\ba}{\begin{equation}}
\newcommand{\ea}{\end{equation}}

\pagecolor{white}

\usepackage{lineno, blindtext}
\usepackage{float}

\usepackage{amsmath}

\begin{document}

\author{Sudipta Sahu$^{\S}$, Rathan Samala$^{\S}$*}
%\title[Three-Stage Semi-Implicit Central Scheme of Relaxation type]{Uniformly Accurate Three-Stage Semi-Implicit Central Schemes for Hyperbolic Systems with Stiff Relaxation}

\title[Three-Stage Semi-Implicit Central Scheme of Relaxation type]{
Semi-Implicit Finite Volume Central Scheme for Hyperbolic Systems of Balance Laws with Adjustable Free Parameter approximation to relaxed source term }

\thanks{ $^{\S}$Department of Humanities and Sciences, Indian Institute of Petroleum and Energy-Visakhapatnam,
India-531035\\
Corresponding author Email: rathans.math@iipe.ac.in (Rathan Samala)}
\date{\today}
\vspace{-1.0cm}
\begin{abstract}
This study develops a generalized framework for constructing a three-stage predictor-corrector finite volume scheme for solving hyperbolic balance laws with stiff relaxation source terms. The proposed scheme extends the central Nessyahu-Tadmor (NT) framework for the nonhomogeneous part of PDE and the integrals of the source terms are computed using an undetermined coefficient approach to introduce a free parameter \(\alpha\), which offers flexibility in formulation. A unified theoretical analysis is then carried out to characterize the properties of the resulting scheme and determine the admissible range of \(\alpha\). Specifically, uniform consistency and asymptotic properties are investigated along with $A$ -stability and $L$-stability, with the corresponding theoretical constraints combined to establish a common admissible range for the free parameter. The numerical performance of the scheme is subsequently evaluated through a sequence of test problems. Smooth-solution tests are used to assess its accuracy, while one-dimensional benchmark problems arising from different mathematical models are considered to examine its performance in comparison with existing numerical schemes. The assessment is further extended to ten two-dimensional Riemann problems to investigate the ability of the scheme to capture shocks and complex wave interactions along with computational cost. \\

\noindent\textbf{Keywords:} Hyperbolic balance laws; Central finite-volume scheme; Stiff relaxation; Asymptotic-preserving property; Stability analysis; Riemann problems; Uniform consistency.
\end{abstract}

\maketitle

\section{Introduction}\label{section:1}
Hyperbolic systems of balance laws provide a mathematical framework for describing transport processes governed by wave propagation and local source-driven mechanisms. They arise in numerous applications, including compressible fluid dynamics, gas kinetics, shallow-water flows, and relaxation phenomena. In one spatial dimension, a general system of balance laws with a relaxation source can be written as
\begin{equation}
\frac{\partial Q}{\partial t}
+\frac{\partial \mathrm{F}(Q)}{\partial x}
=\mathrm{S}(Q,\varepsilon),
\label{eq:1d}
\end{equation}
where $Q\in\mathbb{R}^{m}$ denotes the vector of conserved variables, $\mathrm{F}(Q)$ is the flux function, and $\mathrm{S}(Q,\varepsilon)$ represents the relaxation source term. The corresponding two-dimensional system takes the form
\begin{equation}
\frac{\partial Q}{\partial t}
+\frac{\partial \mathrm{F}_1(Q)}{\partial x}
+\frac{\partial \mathrm{F}_2(Q)}{\partial y}
=\mathrm{S}(Q,\varepsilon),
\label{eq:2d}
\end{equation}
where \(\mathrm{F}_1(Q)\) and \(\mathrm{F}_2(Q)\) denote the fluxes in the \(x\)- and \(y\)-directions, respectively. Numerically solving these systems is challenging because discontinuous waves and source-driven processes may evolve on very different timescales. In particular, when the source term is stiff, the relaxation timescale may be much smaller than the characteristic timescale of wave propagation. A suitable numerical method must therefore preserve the conservative structure, resolve wave interactions accurately, and reproduce the appropriate limiting behaviour as the relaxation parameter approaches zero.

Finite-volume methods are widely used for hyperbolic conservation and balance laws because their conservative formulation is well suited to problems involving discontinuous solutions and complex wave interactions. Classical Godunov-type schemes employ exact or approximate Riemann solvers to capture the characteristic structure of the governing equations \cite{Godunov1959,Roe1981,HartenLaxVanLeer1983}. High-resolution methods based on MUSCL reconstruction, TVD principles, and flux limiters improve spatial accuracy while controlling spurious oscillations near discontinuities \cite{vanLeer1979,Harten1983,Sweby1984}. ENO and WENO methods provide further improvements in the accuracy and resolution of complex wave structures \cite{HartenENO1987,LiuOsherChan1994,JiangShu1996}. In contrast, central schemes provide an alternative approach that does not require explicit Riemann solvers at cell interfaces. In particular, the second-order Nessyahu-Tadmor (NT) scheme employs piecewise-linear reconstruction and staggered control volumes to obtain non-oscillatory approximations of hyperbolic conservation laws \cite{NessyahuTadmor1990}. Subsequent developments introduced non-staggered central schemes and extended central finite-volume methods to multidimensional problems, improving their applicability to nonlinear hyperbolic conservation laws \cite{JiangTadmor1998,JiangLevyLinOsherTadmor1998,KurganovTadmor2000}. Central schemes have also been developed for two-dimensional Riemann problems in gas dynamics without requiring conventional Riemann problem solvers \cite{KurganovTadmor2002}. Their broader development and applicability to nonlinear hyperbolic partial differential equations are reviewed by Kurganov \cite{Kurganov2016}. Central schemes have also been investigated for hyperbolic systems with stiff relaxation terms, highlighting their potential for coupling wave propagation with source-driven dynamics \cite{Pareschi2001,DumbserEnauxToro2008}. These developments provide a useful foundation for constructing central finite-volume methods for balance laws in which the source term introduces additional temporal stability and consistency requirements.

Stiff source terms require careful treatment because an explicit source-term update can impose severe time-step restrictions or cause instability, even when the transport discretization satisfies the usual Courant-Friedrichs-Lewy (CFL) condition \cite{JinLevermore1996,CaflischJinRusso1997}. Implicit and semi-implicit methods can help reduce these restrictions, but their effectiveness depends on how the time discretization, transport terms, and source-term updates work together. Moreover, stability for finite relaxation parameters does not, by itself, guarantee an accurate approximation of the limiting system. This motivates the development of asymptotic-preserving (AP) methods, which are designed to recover a consistent discretization of the limiting dynamics without requiring the small relaxation timescale to be resolved numerically. Foundational contributions to the numerical treatment of hyperbolic systems with stiff relaxation include the methods of Jin \cite{Jin1995}, Jin and Levermore \cite{JinLevermore1996}, and Caflisch, Jin, and Russo \cite{CaflischJinRusso1997}. Recent developments and theoretical perspectives on AP methods for hyperbolic systems with stiff relaxation are reviewed by Boscarino and Russo \cite{BoscarinoRusso2024}. These studies provide the theoretical background for assessing the consistency, stability, and asymptotic behaviour of numerical schemes across different relaxation regimes.

Within the central-scheme framework, Liotta, Romano, and Russo developed second-order central schemes for balance laws with stiff relaxation, including a formulation designed to maintain accuracy across different relaxation regimes \cite{Liotta, russo}. In our previous work, we introduced CS-EBT2, a second-order central finite-volume scheme based on the staggered-grid Nessyahu-Tadmor framework, combining midpoint-based spatial treatment with an implicit trapezoidal update and a backward semi-implicit Taylor expansion for the source-term contribution \cite{SahuMaccaSamala2026}. Building on these developments and our previous work, the present study treats the intermediate-stage location and quadrature coefficients as design parameters and systematically investigates their effects on the consistency, temporal accuracy, stability, and asymptotic-preserving properties of central schemes for stiff balance laws.

To address this objective, we develop a generalized, parameter-dependent second-order central finite-volume scheme using a predictor-corrector formulation and an integral quadrature rule with undetermined coefficients. The intermediate-stage location, \(t^{n+\alpha}=t^n+\alpha\Delta t\), is treated as a free parameter, while the admissible stage locations and corresponding coefficients are determined from consistency and stability requirements. The resulting scheme is analyzed for temporal accuracy, stability, and asymptotic-preserving properties across non-stiff, intermediate, and stiff relaxation regimes. The absolute-stability region of the proposed scheme is compared with that of the classical second-order explicit Runge-Kutta (RK2) method \cite{Butcher2008}, and UCS2 \cite{Liotta} is recovered as a particular realization under appropriate coefficient relations.

The formulation is first developed for one-dimensional balance laws and subsequently extended to two dimensions through a staggered-grid finite-volume construction that accounts for control-volume geometry, directional fluxes, and source-term contributions. Numerical experiments involving one- and two-dimensional Riemann problems and benchmark balance-law models are evaluated to test consistency, wave resolution, and performance across different relaxation regimes. The proposed scheme is compared with UCS2 and a suitable second-order implicit-explicit Runge-Kutta (IMEX-RK2) method \cite{PareschiRusso} to assess its ability to capture shocks, and complex wave dynamics. The effects of intermediate-stage placement and source-term treatment on the accuracy and stability of the scheme are also examined.

The remainder of the paper is organized as follows. Section~\ref{sec:numerical} introduces the proposed generalized central finite-volume scheme in one and two dimensions. Section~\ref{sec:consistencyStability} establishes the theoretical properties of the scheme, including its consistency, temporal accuracy, stability, and asymptotic-preserving behavior. Section~\ref{sec:Numerical} presents a comprehensive set of numerical experiments to evaluate the accuracy, convergence, and performance of the proposed method for various one- and two-dimensional test problems. Finally, Section~\ref{sec:conclusion} summarizes the main findings of the study. Appendix~\ref{appendix} presents the formulations of the UCS2 central scheme and the IMEX Runge–Kutta scheme in one and two dimensions for comparison with the proposed method.

\section{Numerical Approximation and Validation}\label{sec:numerical}
This section describes the numerical schemes considered for the one- and two-dimensional hyperbolic balance laws introduced in \ref{section:1}. We consider a second-order generalized central finite-volume scheme to investigate the accuracy, stability, and numerical behavior of the proposed method across different relaxation regimes. The numerical formulation is presented separately for the one- and two-dimensional cases.
\subsection{One-Dimensional Numerical Scheme}
We first consider the one-dimensional balance law
\eqref{eq:1d}. For the scalar case, $m=1$, the conserved variable $Q$ reduces to a scalar variable $q$. Introducing the source-term
notation
\begin{equation}
\mathrm{s}(q,\varepsilon)
=\frac{1}{\varepsilon}\mathrm{g}(q),
\label{eq:source_scaling}
\end{equation}
the governing equation can be written as
\begin{equation}
q_t+\partial_x\mathrm{f}(q) =\frac{1}{\varepsilon}\mathrm{g}(q).
\label{eq:hyperbolic}
\end{equation}
For the one-dimensional formulation, let $\bar{q}_j(t)$ denote the cell average of $q$ over the interval $[x_{j-\frac12},x_{j+\frac12}]$,
\begin{equation}
\bar{q}_j(t) = \frac{1}{\Delta x} \int_{x_{j-\frac12}}^{x_{j+\frac12}} q(x,t)\,dx.
\end{equation}
Within each cell, a piecewise linear reconstruction is employed,
\begin{equation}
\mathcal{L}_j(x,t) = \bar{q}_j(t) + (x-x_j)\frac{\bar{q}'_j}{\Delta x},
\qquad x\in[x_{j-\frac12},x_{j+\frac12}],
\end{equation}
where the slopes $\bar{q}'_j$ are computed using the MinMod limiter
\begin{equation}
\bar{q}'_j = \operatorname{MM}\!\left(\bar{q}_{j+1}-\bar{q}_j,\; \bar{q}_j-\bar{q}_{j-1}\right),
\end{equation}
with limiter
\begin{equation}\label{eq:minmod}
\operatorname{MM}(a,b) =
\begin{cases}
\operatorname{sgn}(a)\,\min(|a|,|b|), & \text{if } \operatorname{sgn}(a)=\operatorname{sgn}(b),\\
0, & \text{otherwise},
\end{cases}
\end{equation}
which provides a non-oscillatory reconstruction while retaining second-order accuracy in smooth regions.
Integrating the balance law over the space-time control volume 
$[x_j, x_{j+1}] \times [t^n, t^{n+1}]$, and using the piecewise linear polynomial reconstruction for $q(x,t^n)$ in $[x_j, x_{j+\frac12}]$ and $[x_{j+\frac12}, x_{j+1}],$ the cell average over the spatial cell $[x_j, x_{j+1}]$ can be expressed as
\begin{equation}\label{corrector1}
    \frac{1}{\Delta x} \int_{x_j}^{x_{j+1}} q(x,t^n)\, dx =
    \frac{1}{\Delta x} \bigg(\int_{x_j}^{x_{j+\frac12}} q(x,t^n)\, dx +\int_{x_{j+\frac12}}^{x_{j+1}} q(x,t^n)\, dx\bigg) = \frac{1}{2} \big(\bar{q}_j^n + \bar{q}_{j+1}^n \big) 
      + \frac{1}{8} \big(\bar{q}'_j - \bar{q}'_{j+1} \big).
\end{equation}
and the contribution of the flux term over the space-time cell is given by
\begin{equation}\label{corrector2}
\frac{1}{\Delta x}\int_{t^n}^{t^{n+1}}
\Bigl[\mathrm{f}(q(x_{j+1},t))-\mathrm{f}(q(x_j,t))\Bigr]\,dt.
\end{equation}
This temporal integral is approximated by the midpoint rule,
\begin{equation}\label{corrector3}
\int_{t^n}^{t^{n+1}}
\Bigl[\mathrm{f}(q(x_{j+1},t))-\mathrm{f}(q(x_j,t))\Bigr]\,dt
\;\approx\;
\Delta t\Bigl[
\mathrm{f}(\bar{q}_{j+1}^{\,n+\frac12})
-
\mathrm{f}(\bar{q}_j^{\,n+\frac12})
\Bigr],
\end{equation}
where the staggered half-time interface states
$\bar{q}_j^{\,n+\frac12}$ are obtained from the NT-type predictor step.

The space-time integral of the source term,
\begin{equation}
\mathcal{I}_{\mathrm{g}} =
\int_{t^n}^{t^{n+1}}
\int_{x_{j}}^{x_{j+1}}
\mathrm{g}(q(x,t))\,dx\,dt,
\end{equation}
is evaluated approximately by applying the midpoint rule for the spatial integration and an undetermined-coefficient quadrature rule for the temporal integration.
\begin{equation}
\mathcal{I}_{\mathrm{g}} \approx
\Delta t
\left(
w_1 \mathrm{g}(\bar{q}_{j+1/2}^{\,n+\alpha}) + w_2 \mathrm{g}(\bar{q}_{j+\frac12}^{\,n+1})\right).
\end{equation}
Here $x_{j+1/2} = \dfrac{x_j + x_{j+1}}{2}$. Then
\begin{equation}
\mathcal{I}_{\mathrm{g}} \approx
\Delta x\,\Delta t
\left(
\frac{w_1}{2} \mathrm{g}(\bar{q}_j^{\,n+\alpha}) +\frac{w_1}{2}  \mathrm{g}(\bar{q}_{j+1}^{\,n+\alpha}) + w_2 \mathrm{g}(\bar{q}_{j+\frac12}^{\,n+1})\right).
\end{equation}
The first two evaluations rely on predictor states, while the last term is treated implicitly at the new time level, ensuring stability in the stiff relaxation regime.

\noindent Combining the above discretizations yields the predictor-corrector formulation. The predictor stage reads
\begin{align}
\bar{q}_j^{\,n+\frac12}
&=
\bar{q}_j^n
+
\frac{\Delta t}{2}
\left(\frac{1}{\varepsilon}
\mathrm{g}(\bar{q}_j^{\,n+\frac12})
-
\frac{\mathrm{f}'_j}{\Delta x}
\right),\\
\bar{q}_j^{\,n+\alpha}
&=
\bar{q}_j^n+ \alpha \Delta t
\left(\frac{1}{\varepsilon}
\mathrm{g}(\bar{q}_j^{\,n+\alpha})
-
\frac{\mathrm{f}'_j}{\Delta x}
\right).
\end{align}
The corrector stage of Predictor-Corrector scheme is given by
\begin{align}
\bar{q}_{j+\frac12}^{\,n+1}
&=
\frac{1}{2}
(\bar{q}_j^n+\bar{q}_{j+1}^n)
+ \frac{1}{8}
(\bar{q}'_j-\bar{q}'_{j+1})
- \lambda
\Bigl[
\mathrm{f}(\bar{q}_{j+1}^{\,n+\frac12})
- \mathrm{f}(\bar{q}_j^{\,n+\frac12})
\Bigr]
\notag\\
&\quad
+\frac{\Delta t}{\varepsilon} \left(\frac{w_1}{2}\mathrm{g}(\bar{q}_j^{\,n+\alpha}) + \frac{w_1}{2}\mathrm{g}(\bar{q}_{j+1}^{\,n+\alpha}) + w_2 {\mathrm{g}}(\bar{q}_{j+\frac12}^{\,n+1}) \right),
\end{align}
where $\lambda=\Delta t/\Delta x$ and the MinMod limiter is used in the computation of $\bar{q}'_j$ and $\mathrm{f}'_j$.

Consider the two-point quadrature formula
\begin{equation}\label{1}
    \int_n^{n+1} q(t)\,dt
    \approx
    w_1 q(n+\alpha) + w_2 q(n+1),
\end{equation}
where $\alpha$ is a prescribed parameter and $w_1,w_2$ are
undetermined weights. The weights are determined by requiring the quadrature rule to be exact for $\mathbb{P}_1$, the space of polynomials of degree at most one.

It is sufficient to impose exactness for the basis
$\{1,t\}$ of $\mathbb{P}_1$. For $q(t)=1$, equation~(\ref{1}) yields
\begin{equation*}
    \int_n^{n+1}1\,dt
    =
    w_1+w_2,
\end{equation*}
and hence
\begin{equation}\label{2}
    w_1+w_2=1.
\end{equation}
For $q(t)=t$, the exactness condition becomes
\begin{equation*}
    \int_n^{n+1}t\,dt
    =
    w_1(n+\alpha)+w_2(n+1).
\end{equation*}
Since
\begin{equation*}
    \int_n^{n+1}t\,dt
    =
    \left.\frac{t^2}{2}\right|_n^{n+1}
    =
    n+\frac{1}{2},
\end{equation*}
we have
\begin{equation*}
    n+\frac{1}{2}
    =
    w_1(n+\alpha)+w_2(n+1).
\end{equation*}
Using~(\ref{2}), this can be reduced to
\begin{equation}\label{3}
    \alpha w_1+w_2=\frac{1}{2}.
\end{equation}
Consequently, the weights satisfy the linear system
\begin{equation*}
    \begin{pmatrix}
        1 & 1\\
        \alpha & 1
    \end{pmatrix}
    \begin{pmatrix}
        w_1\\
        w_2
    \end{pmatrix}
    =
    \begin{pmatrix}
        1\\
        \dfrac{1}{2}
    \end{pmatrix}.
\end{equation*}
Solving this system gives
\begin{equation}\label{coeff_quad}
    w_1=\frac{1}{2(1-\alpha)},
    \quad
    w_2=\frac{1-2\alpha}{2(1-\alpha)},
    \qquad \alpha\neq1.
\end{equation}
Therefore, the corresponding quadrature formula is
\begin{equation*}\label{4}
    \int_n^{n+1}q(t)\,dt
    \approx
    \frac{1}{2(1-\alpha)}q(n+\alpha)
    +
    \frac{1-2\alpha}{2(1-\alpha)}q(n+1).
\end{equation*}
Then the corrector stage is given by
\begin{align}\label{corrector:general}
\bar{q}_{j+\frac12}^{\,n+1}
&=
\frac{1}{2}
(\bar{q}_j^n+\bar{q}_{j+1}^n)
+
\frac{1}{8}
(\bar{q}'_j-\bar{q}'_{j+1})
-
\lambda
\Bigl[
\mathrm{f}(\bar{q}_{j+1}^{\,n+\frac12})
-
\mathrm{f}(\bar{q}_j^{\,n+\frac12})
\Bigr]
\notag\\
&\quad
+\frac{\Delta t}{\varepsilon} \left(\frac{1}{4(1-\alpha)}\mathrm{g}(\bar{q}_j^{\,n+\alpha}) + \frac{1}{4(1-\alpha)}\mathrm{g}(\bar{q}_{j+1}^{\,n+\alpha}) + \frac{1-2\alpha}{2(1-\alpha)} \mathrm{g}(\bar{q}_{j+\frac12}^{\,n+1}) \right),
\end{align}
where $\lambda=\Delta t/\Delta x$. The MinMod limiter is employed in the evaluation of $\bar{q}'_j$ and $\mathrm{f}'_j$, with the corresponding predictor stages incorporated into the computational procedure.

% After presenting the one-dimensional formulations, we proceed to their extension to two-dimensional non-homogeneous systems. The associated
% finite-volume discretization is described
% in the following subsection.
\subsection{Two-Dimensional Numerical Scheme}
We consider the two-dimensional hyperbolic balance law with relaxation
\begin{equation}\label{2d:EQuation}
q_t + \partial_x \mathrm{f}_1(q) + \partial_y \mathrm{f}_2(q) = \frac{1}{\varepsilon} \mathrm{g}(q),
\end{equation}
supplemented with the initial condition $q(x,y,0)=q_0(x,y)$. This equation is discretized using a central finite-volume framework.\\
The computational domain is partitioned into a uniform Cartesian grid composed of rectangular control volumes
\[
C_{rs} = [x_{r-\frac{1}{2}},\,x_{r+\frac{1}{2}}] \times [y_{s-\frac{1}{2}},\,y_{s+\frac{1}{2}}],
\]
where the cell centers are defined by 
\[
x_r = \frac{1}{2}(x_{r-\frac{1}{2}} + x_{r+\frac{1}{2}}), 
\qquad
y_s = \frac{1}{2}(y_{s-\frac{1}{2}} + y_{s+\frac{1}{2}}).
\]
The corresponding mesh sizes are $\Delta x = x_{r+\frac{1}{2}} - x_{r-\frac{1}{2}}$ and $\Delta y = y_{s+\frac{1}{2}} - y_{s-\frac{1}{2}}$.\\
Within the finite-volume setting, the numerical approximation is expressed in terms of cell averages
\begin{equation*}
\bar{q}_{rs}^n = \frac{1}{\Delta x \Delta y}
\int_{y_{s-\frac{1}{2}}}^{y_{s+\frac{1}{2}}}
\int_{x_{r-\frac{1}{2}}}^{x_{r+\frac{1}{2}}}
q(x,y,t^n)\,dx\,dy.
\end{equation*}

\medskip

\noindent\textbf{Reconstruction.}
Starting from the piecewise constant data $\bar{q}_{rs}^n$, a second-order accurate representation is obtained through a piecewise linear reconstruction in each cell:
\begin{equation}\label{homo1}
q(x,y,t^n) =
\left[
\bar{q}_{rs}^n 
+ q^{\dagger}_{rs}\frac{x - x_r}{\Delta x} 
+ q^{\ast}_{rs}\frac{y - y_s}{\Delta y}
\right]\chi_{rs}(x,y),
\end{equation}
where $\chi_{rs}(x,y)$ is the characteristic function of $C_{rs}$. The coefficients $q^{\dagger}_{rs}$ and $q^{\ast}_{rs}$ denote the discrete slopes in the $x$- and $y$-directions, respectively, obtained from neighboring cell averages.
For second-order accuracy, these slopes satisfy
\begin{equation*}
q^{\dagger}_{rs} = \Delta x\, q_x(x_p,y_q,t^n) + \mathcal{O}(\Delta x^2),
\qquad
q^{\ast}_{rs} = \Delta y\, q_y(x_p,y_q,t^n) + \mathcal{O}(\Delta y^2).
\end{equation*}

\begin{figure}[ht]
\centering
\begin{tikzpicture}[
    >=Latex,
    every node/.style={inner sep=2pt}
]

\draw[->, thick] (-2,0) -- (7,0) node[right] {$x$};
\draw[->, thick] (0,-2) -- (0,7) node[above] {$y$};
\draw (4,0.12) -- (4,-0.12);
\node[below=6pt] at (4,0) {$x_{j+\frac{1}{2}}$};
\draw (0.12,4) -- (-0.12,4);
\node[left=6pt] at (0,4) {$y_{k+\frac{1}{2}}$};

% Guide lines
\draw[densely dotted] (4,0) -- (4,4);
\draw[densely dotted] (0,4) -- (4,4);
\draw[dashed] (2,2) rectangle (6,6);
\draw[dashed] (4,2) -- (4,6);
\draw[dashed] (2,4) -- (6,4);

\node at (3,5.3) {$C_{j,k+1}$};
\node at (5,5.3) {$C_{j+1,k+1}$};
\node at (3,2.7) {$C_{j,k}$};
\node at (5,2.7) {$C_{j+1,k}$};

\draw[thick] (3.0,3.0) rectangle (5.0,5.0);

% Internal split
\draw[dashed] (4,3.0) -- (4,5.0);
\draw[dashed] (3.0,4) -- (5.0,4);

\node[font=\footnotesize] at (3.5,4.5) {NW};
\node[font=\footnotesize] at (4.5,4.5) {NE};
\node[font=\footnotesize] at (3.5,3.5) {SW};
\node[font=\footnotesize] at (4.5,3.5) {SE};

\node[below] at (4,2) {$\omega_{j+\frac{1}{2}}$};
\node[left] at (2,4) {$\Omega_{k+\frac{1}{2}}$};

\coordinate (center) at (4,4);

\fill (4,4) circle (1.5pt);

\draw[->, thick] (6.8,4.8) -- (center);
\node[right] at (6.8,4.8) {$C_{j+\frac{1}{2},k+\frac{1}{2}}$};
\end{tikzpicture}
\caption{Geometric decomposition of the staggered control volume.}
\label{fig:staggered-grid}
\end{figure}
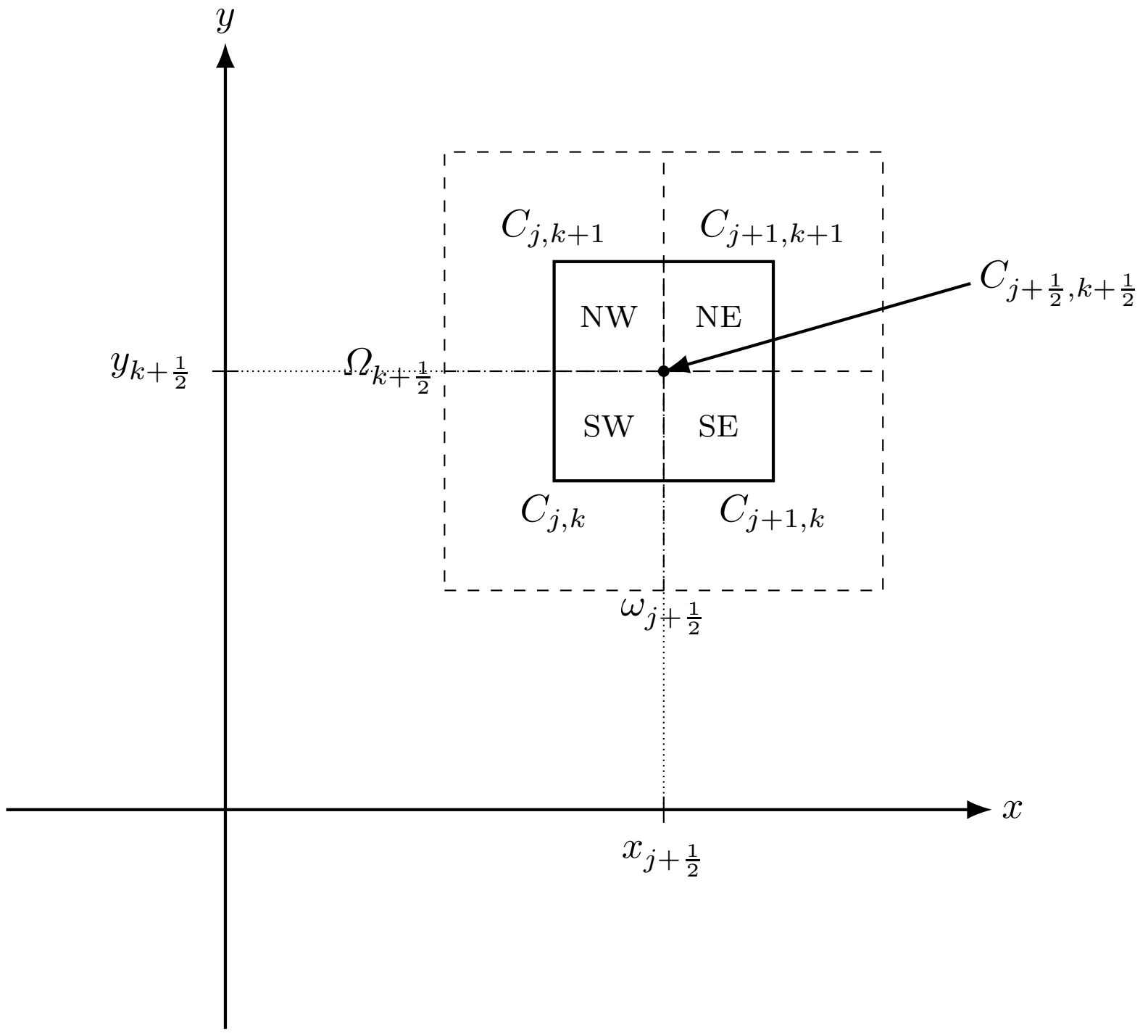

\medskip
\noindent\textbf{Staggered evolution.}
Following the central Nessyahu-Tadmor (NT) approach, the scheme evolves the reconstructed solution by averaging over staggered control volumes. These staggered cells are defined as
\[
C_{j+\frac{1}{2},\,k+\frac{1}{2}} = [x_j,\,x_{j+1}] \times [y_k,\,y_{k+1}] = \omega_{j+\frac{1}{2}} \times \Omega_{k+\frac{1}{2}},
\]
with associated staggered averages
\begin{equation*}
\bar{q}_{j+\frac{1}{2},\,k+\frac{1}{2}}(t)
=
\frac{1}{\Delta x \Delta y}
\int_{y_k}^{y_{k+1}}
\int_{x_j}^{x_{j+1}}
q(x,y,t)\,dx\,dy.
\end{equation*}
Define the mesh ratios
$\lambda = \frac{\Delta t}{\Delta x}, 
\quad 
\mu = \frac{\Delta t}{\Delta y}.$
Integrating \eqref{2d:EQuation} over the spacetime control volume 
$C_{j+\frac{1}{2},\,k+\frac{1}{2}} \times [t^n, t^{n+1}]$ yields the finite-volume update formula
\begin{equation}\label{integral:eq}
\begin{aligned}
\bar{q}_{j+\frac{1}{2},\,k+\frac{1}{2}}^{\,n+1}
&= \dashint_{C_{j+\frac{1}{2},\,k+\frac{1}{2}}} q(x,y,t^n)\,dx\,dy- \lambda 
\int_{t^n}^{t^{n+1}} \dashint_{y \in \Omega_{k+\frac{1}{2}}}
\left[
\mathrm{f}_1(q(x_{j+1},y,\tau)) - \mathrm{f}_1(q(x_j,y,\tau))
\right] dy\, d\tau \\
&\; - \mu 
\int_{t^n}^{t^{n+1}} 
\dashint_{x \in \omega_{j+\frac{1}{2}}}
\left[
\mathrm{f}_2(q(x,y_{k+1},\tau)) - \mathrm{f}_2(q(x,y_k,\tau))
\right] dx\, d\tau + \frac{1}{\varepsilon}
\int_{t^n}^{t^{n+1}}\dashint_{C_{j+\frac{1}{2},\,k+\frac{1}{2}}}
\mathrm{g}(q(x,y,\tau))\,dx\,dy\, d\tau.
\end{aligned}
\end{equation}
\noindent\textbf{Evaluation of the staggered average.}
We begin by evaluating the normalized integral
\[
\dashint_{C_{j+\frac{1}{2},\,k+\frac{1}{2}}}
q(x,y,t^n)\,dx\,dy,
\]
which represents the average of the reconstructed solution over the
staggered control volume \(C_{j+\frac{1}{2},\,k+\frac{1}{2}}\). 
Here, the notation
\[
\dashint_{\mathcal{D}}(\cdot)\,d\mathbf{x}
:=\frac{1}{|\mathcal{D}|}\int_{\mathcal{D}}(\cdot)\,d\mathbf{x}
\]
denotes the normalized integral over a domain \(\mathcal{D}\), where
\(|\mathcal{D}|\) represents its measure. As illustrated in
Figure~\ref{fig:staggered-grid} \cite{JiangTadmor1998}, the staggered
control volume receives contributions from the four adjacent cells,
namely \(C_{j,k}\), \(C_{j+1,k}\), \(C_{j+1,k+1}\), and \(C_{j,k+1}\). Using the piecewise linear reconstruction \eqref{homo1}, the integral over the staggered cell is decomposed into contributions from four subcells, corresponding to the southwest (SW), southeast (SE), northeast (NE), and northwest (NW) regions. Each contribution is computed by integrating the reconstructed polynomial over the corresponding subregion.

\medskip
\noindent Southwest subcell:
\begin{equation}
\dashint_{C_{SW}} q(x,y,t^n)\,dx\,dy
= \frac{1}{4}\bar{q}_{j,k}^n
+ \frac{1}{16}\left( q^{\dagger}_{j,k} + q^{\ast}_{j,k} \right),
\end{equation}
\noindent Southeast subcell:
\begin{equation}
\dashint_{C_{SE}} q(x,y,t^n)\,dx\,dy
= \frac{1}{4}\bar{q}_{j+1,k}^n
+ \frac{1}{16}\left( -q^{\dagger}_{j+1,k} + q^{\ast}_{j+1,k} \right),
\end{equation}
\noindent Northeast subcell:
\begin{equation}
\dashint_{C_{NE}} q(x,y,t^n)\,dx\,dy
= \frac{1}{4}\bar{q}_{j+1,k+1}^n
- \frac{1}{16}\left( q^{\dagger}_{j+1,k+1} + q^{\ast}_{j+1,k+1} \right),
\end{equation}
\noindent Northwest subcell:
\begin{equation}
\dashint_{C_{NW}} q(x,y,t^n)\,dx\,dy
= \frac{1}{4}\bar{q}_{j,k+1}^n
+ \frac{1}{16}\left( q^{\dagger}_{j,k+1} - q^{\ast}_{j,k+1} \right).
\end{equation}
By summing the above four contributions, we obtain the average over the staggered control volume:
\begin{equation}
\begin{aligned}
\dashint_{C_{j+\frac{1}{2},\,k+\frac{1}{2}}} & q(x,y,t^n)\,dx\,dy
= \frac{1}{4}\left(
\bar{q}_{j,k}^n + \bar{q}_{j+1,k}^n
+ \bar{q}_{j,k+1}^n + \bar{q}_{j+1,k+1}^n
\right) \\
& \quad+ \frac{1}{16} \Big[(q^{\dagger}_{j,k} - q^{\dagger}_{j+1,k})
+ (q^{\dagger}_{j,k+1} - q^{\dagger}_{j+1,k+1})+ (q^{\ast}_{j,k} - q^{\ast}_{j,k+1})
+ (q^{\ast}_{j+1,k} - q^{\ast}_{j+1,k+1})
\Big].
\end{aligned}
\end{equation}
\noindent\textbf{Flux approximation.}
Next, we approximate the flux integrals appearing in \eqref{integral:eq}. To maintain second-order accuracy in both space and time, we employ the midpoint rule in time together with a second-order quadrature in space.

Accordingly, the flux across the right interface is approximated as
\begin{equation}
\begin{aligned}
\int_{t^n}^{t^{n+1}} 
\dashint_{y \in \Omega_{k+\frac{1}{2}}} 
\mathrm{f}_1(q(x_{j+1},y,\tau))\,dy\,d\tau
\approx
\frac{\Delta t}{2}
\left[
\mathrm{f}_1(q_{j+1,k}^{n+\frac{1}{2}})
+
\mathrm{f}_1(q_{j+1,k+1}^{n+\frac{1}{2}})
\right].
\end{aligned}
\end{equation}
Similarly, the remaining flux contributions are approximated by
\begin{equation}
\int_{t^n}^{t^{n+1}} 
\dashint_{y \in \Omega_{k+\frac{1}{2}}} 
\mathrm{f}_1(q(x_j,y,\tau))\,dy\,d\tau
\approx
\frac{\Delta t}{2}
\left[
\mathrm{f}_1(q_{j,k}^{n+\frac{1}{2}})
+
\mathrm{f}_1(q_{j,k+1}^{n+\frac{1}{2}})
\right],
\end{equation}
\begin{equation}
\int_{t^n}^{t^{n+1}} 
\dashint_{x \in \omega_{j+\frac{1}{2}}} 
\mathrm{f}_2(q(x,y_{k+1},\tau))\,dx\,d\tau
\approx
\frac{\Delta t}{2}
\left[
\mathrm{f}_2(q_{j,k+1}^{n+\frac{1}{2}})
+
\mathrm{f}_2(q_{j+1,k+1}^{n+\frac{1}{2}})
\right],
\end{equation}

\begin{equation}
\int_{t^n}^{t^{n+1}} 
\dashint_{x \in \omega_{j+\frac{1}{2}}} 
\mathrm{f}_2(q(x,y_k,\tau))\,dx\,d\tau
\approx
\frac{\Delta t}{2}
\left[
\mathrm{f}_2(q_{j,k}^{n+\frac{1}{2}})
+
\mathrm{f}_2(q_{j+1,k}^{n+\frac{1}{2}})
\right].
\end{equation}
Finally, the source term contribution, denoted by \(\mathcal{I}_g\), must be approximated. Its temporal treatment depends on the chosen time discretization and the stiffness parameter \(\varepsilon\), as detailed below.

\medskip
\noindent\textbf{Source term approximation.}
The spatial discretization, reconstruction, and flux evaluation remain unchanged, while the source term is approximated using different formulations.
Using the midpoint rule in space and a undetermined coefficient quadrature rule in time,
\[
\mathcal{I}_{\mathrm{g}} \approx \dfrac{\Delta x \Delta y \Delta t}{\varepsilon} 
\left[
w_1\,\mathrm{g}(\bar{q}_{j+\frac{1}{2},k+\frac{1}{2}}^{n+\alpha})
+w_2\,\mathrm{g}(\bar{q}_{j+\frac{1}{2},k+\frac{1}{2}}^{n+1})
\right].
\]
The intermediate value is approximated as
\begin{eqnarray*} 
\begin{aligned} \mathrm{g}(\bar{q}_{i+1/2,j+1/2}^{n+\alpha})&\approx \dfrac{1}{2}\left(\mathrm{g}(\bar{q}_{i,j+1/2}^{n+\alpha}) + \mathrm{g}(\bar{q}_{i+1,j+1/2}^{n+\alpha})\right),\\ 
&\approx \dfrac{1}{4}\left(\mathrm{g}(\bar{q}_{i,j+1}^{n+\alpha})+\mathrm{g}(\bar{q}_{i+1,j+1}^{n+\alpha})+ \mathrm{g}(\bar{q}_{i,j}^{n+\alpha})+ \mathrm{g}(\bar{q}_{i+1,j}^{n+\alpha})\right). 
\end{aligned} 
\end{eqnarray*}
The additional predictor is
\begin{equation*}
\bar{q}_{j,k}^{n+\alpha}
=\bar{q}_{j,k}^n
-{\alpha \lambda}\, \mathrm{f}_1(q)^\dagger_{j,k}
-{\alpha \mu}\, \mathrm{f}_2(q)^\ast_{j,k}
+ \frac{\alpha \Delta t}{\varepsilon} \mathrm{g}(\bar{q}_{j,k}^{\,n+\alpha}).
\end{equation*}
The corrector step follows the same convective structure as above, with the source term replaced by
\begin{equation*}
\begin{aligned} 
\bar{q}_{j+\frac{1}{2},k+\frac{1}{2}}^{\,n+1} 
&= \frac{1}{4}\Big(\bar{q}_{j,k}^n + \bar{q}_{j+1,k}^n + \bar{q}_{j,k+1}^n + \bar{q}_{j+1,k+1}^n \Big) \\ 
&+ \frac{1}{16}\Big[(q^{\dagger}_{j,k} - q^{\dagger}_{j+1,k}) + (q^{\dagger}_{j,k+1} - q^{\dagger}_{j+1,k+1}) + (q^{\ast}_{j,k} - q^{\ast}_{j,k+1}) + (q^{\ast}_{j+1,k} - q^{\ast}_{j+1,k+1}) \Big] \\ 
&- \frac{\lambda}{2}\Big[\mathrm{f}_1(\bar{q}_{j+1,k}^{n+\frac{1}{2}}) - \mathrm{f}_1(\bar{q}_{j,k}^{n+\frac{1}{2}}) + \mathrm{f}_1(\bar{q}_{j+1,k+1}^{n+\frac{1}{2}}) - 
\mathrm{f}_1(\bar{q}_{j,k+1}^{n+\frac{1}{2}}) \Big] \\ 
&- \frac{\mu}{2}\Big[\mathrm{f}_2(\bar{q}_{j,k+1}^{n+\frac{1}{2}}) - \mathrm{f}_2(\bar{q}_{j,k}^{n+\frac{1}{2}}) + \mathrm{f}_2(\bar{q}_{j+1,k+1}^{n+\frac{1}{2}}) - \mathrm{f}_2(\bar{q}_{j+1,k}^{n+\frac{1}{2}}) \Big] \\
&+ \dfrac{\Delta t}{\varepsilon} \left[\dfrac{w_1}{4} \left(\mathrm{g}(\bar{q}_{j,k+1}^{n+\alpha}) + \mathrm{g}(\bar{q}_{j+1,k+1}^{n+\alpha})+ \mathrm{g}(\bar{q}_{j,k}^{n+\alpha})+ \mathrm{g}(\bar{q}_{j+1,k}^{n+\alpha})\right) + w_2 \, \mathrm{g}\left(\bar{q}_{j+\frac{1}{2},k+\frac{1}{2}}^{n+1}\right)\right]. 
\end{aligned} 
\end{equation*}
Here $w_1=\dfrac{1}{2(1-\alpha)}$ and $w_2=\dfrac{1-2\alpha}{2(1-\alpha)}$. The proposed scheme for the two-dimensional problem is formulated through the following predictor and corrector steps, analogous to the one-dimensional case:

\medskip
\noindent\textbf{Predictor step:}
\label{generalised:scheme2d} 
Solve for $\bar{q}_{j,k}^{\,n+1/2}$ and $\bar{q}_{j,k}^{\,n+\alpha}$ from 
\begin{equation}
\bar{q}_{j,k}^{\,n+1/2}= \bar{q}_{j,k}^n + \dfrac{\Delta t}{2}\left[\dfrac{1}{\varepsilon} \mathrm{g}(\bar{q}_{j,k}^{\,n+1/2}) - \dfrac{(f_1)^{\dagger}_{j,k}}{\Delta x}-\dfrac{(f_2)^{\ast}_{j,k}}{\Delta y}\right]. 
\end{equation} 
\begin{equation} 
\bar{q}_{j,k}^{\,n+\alpha}= \bar{q}_{j,k}^n + {\alpha \Delta t}\left[\dfrac{1}{\varepsilon} \mathrm{g}(\bar{q}_{j,k}^{\,n+\alpha}) - \dfrac{(f_1)^{\dagger}_{j,k}}{\Delta x}-\dfrac{(f_2)^{\ast}_{j,k}}{\Delta y}\right]. 
\end{equation} 
\noindent\textbf{Corrector step:} Then update to $\bar{q}_{j+1/2,k+1/2}^{n+1}$ by 
\begin{equation} 
\begin{aligned} 
\bar{q}_{j+1/2,k+1/2}^{\,n+1}
&= \left\langle\dfrac{1}{4}\left(\bar{q}_{j,\cdot}^n + \bar{q}_{j+1,\cdot}^n\right) + \dfrac{1}{8}\left(q^{\dagger}_{i,\cdot} - q^{\dagger}_{j+1,\cdot}\right) - \lambda\left[\mathrm{f}_1(\bar{q}_{j+1,\cdot}^{\,n+1/2}) - \mathrm{f}_1(\bar{q}_{j,\cdot}^{\,n+1/2})\right] \right.\\ 
& \qquad \left.+\dfrac{\Delta t}{\varepsilon}\left(\dfrac{1}{8(1-\alpha)}\mathrm{g}(\bar{q}_{j,\cdot}^{n+\alpha})+\dfrac{1}{8(1-\alpha)}\mathrm{g}(\bar{q}_{j+1,\cdot}^{n+\alpha})\right)\right\rangle_{k+1/2}\\ 
&\quad +\left\langle\dfrac{1}{4}\left(\bar{q}_{\cdot,k}^n + \bar{q}_{\cdot,k+1}^n\right) + \dfrac{1}{8}\left(q^{\ast}_{\cdot,k} - q^{\ast}_{\cdot,k+1}\right) - \mu\left[\mathrm{f}_2(\bar{q}_{\cdot,k+1}^{\,n+1/2}) - \mathrm{f}_2(\bar{q}_{\cdot,k}^{\,n+1/2})\right] \right.\\ 
& \qquad \left.+\dfrac{\Delta t}{\varepsilon}\left(\dfrac{1}{8(1-\alpha)}\mathrm{g}(\bar{q}_{\cdot,k}^{n+\alpha})+\dfrac{1}{8(1-\alpha)}\mathrm{g}(\bar{q}_{\cdot,k+1}^{n+\alpha})\right)\right\rangle_{j+1/2} +\frac{\Delta t (1-2\alpha)}{2(1-\alpha)\varepsilon} \mathrm{g}(\bar{q}_{j+1/2,k+1/2}^{n+1}). 
\end{aligned} 
\end{equation}
The proposed scheme follows a three-step predictor-corrector procedure. Starting from the cell averages $\bar{q}_{j,k}^n$, intermediate values are first computed using the predictor step, followed by the corrector step to obtain $\bar{q}_{j,k}^{n+1}$. The resulting scheme is designed to achieve second-order accuracy and non-oscillatory behavior without requiring Riemann solvers. The corresponding discretization is employed throughout the numerical experiments presented in the following section.
\section{Theoretical Analysis}
\label{sec:consistencyStability}
In this section, we investigate the consistency, asymptotic-preserving, and stability properties of the proposed numerical scheme, with particular emphasis on the roles of the intermediate-stage parameter \(\alpha\) and the relaxation parameter \(\varepsilon\).
\subsection{Consistency Analysis in Different Asymptotic Regimes}
\label{subsec:consistency}
The consistency analysis follows the approach of \cite{Liotta}. For brevity, we present the main expansions and consistency arguments, omitting intermediate calculations that proceed analogously to those in \cite{Liotta}. We consider the rarefied, fluid-dynamic, and thin-relaxation regimes under the assumption of sufficiently smooth solutions and the time-step relation $\Delta t=\lambda h$, where $h=\Delta x$ and $\lambda>0$ is fixed. The relaxation system is
\begin{equation}\label{model}
\begin{aligned}
\begin{cases}
u_t+v_x=0,\\[2mm]
v_t+u_x=-\frac{1}{\varepsilon}(v-au),
\end{cases}
\end{aligned}
\end{equation}
where $a\in\mathbb R$, $|a|<1$, and $\varepsilon>0$, with initial conditions
\[
u(x,0)=u_0(x),\quad v(x,0)=v_0(x).
\]
The corresponding state vector and flux are
\[
q=\begin{pmatrix}u\\v\end{pmatrix},
\qquad
\mathrm f(q)=\begin{pmatrix}v\\u\end{pmatrix}.
\]
The generalized scheme employs the weights
\[
w_1=\frac{1}{2(1-\alpha)},\quad
w_2=\frac{1-2\alpha}{2(1-\alpha)},
\]
where $0<\alpha<\tfrac12$. These weights are nonnegative and well-defined throughout the specified range.\\\\
\noindent{\textbf{Cell averages and reconstruction assumptions.}} Since the scheme evolves cell averages, consistency is established by comparing the numerical updates with the corresponding exact cell-average expansions. For a smooth function $u$, define
\[
\bar u(x,t)=\frac1h\int_{-h/2}^{h/2}
u(x+\nu,t)\,d\nu,
\quad
\bar u^0(x)=\frac1h\int_{-h/2}^{h/2}
u_0(x+\nu)\,d\nu,
\]
and analogously define $\bar v(x,t)$ and $\bar v^0(x)$. When evaluating at the interface $x_{j+\frac12}$, we write
\[
\bar u^0=\bar u^0(x_{j+\frac12}),
\quad
\bar v^0=\bar v^0(x_{j+\frac12}).
\]
Barred derivatives denote cell averages of the corresponding spatial derivatives, evaluated at the initial time unless otherwise indicated. Expanding the neighboring cell averages about the interface gives
\begin{equation}\label{cell:neighbor}
\begin{aligned}
\begin{cases}
\bar u_j^0
=\bar u^0-\frac h2\bar u_x^0
+\frac{h^2}{8}\bar u_{xx}^0+\mathcal O(h^3),\\[2mm]
\bar u_{j+1}^0
=\bar u^0+\frac h2\bar u_x^0
+\frac{h^2}{8}\bar u_{xx}^0+\mathcal O(h^3),
\end{cases}
\end{aligned}
\end{equation}
with analogous expansions for $\bar v_j^0$ and $\bar v_{j+1}^0$.

We assume first-order accurate, smoothly varying reconstruction slopes:
\begin{equation}\label{ass:2}
\begin{aligned}
\begin{cases}
\frac{u'_j}{h}
=u_x(x_j)+c(x_j)h+\mathcal O(h^2),\\[2mm]
\frac{v'_j}{h}
=v_x(x_j)+d(x_j)h+\mathcal O(h^2),
\end{cases}
\end{aligned}
\end{equation}
where $c(x)$ and $d(x)$ are smooth. The reconstruction consequently satisfies
\begin{equation}\label{ass:3}
\begin{aligned}
\begin{cases}
\bar u^0(x_{j+\frac12})
=\frac12(\bar u_j^0+\bar u_{j+1}^0)
-\frac18(u'_{j+1}-u'_j)+\mathcal O(h^3),\\[2mm]
\bar u_x^0(x_{j+\frac12})
=\frac{\bar u_{j+1}^0-\bar u_j^0}{h}
+\mathcal O(h^2),\\[2mm]
\bar v_{xx}^0(x_{j+\frac12})
=\frac{\bar v'_{j+1}-\bar v'_j}{h^2}
+\mathcal O(h).
\end{cases}
\end{aligned}
\end{equation}
Analogous relations are assumed for the corresponding quantities involving $v$ and higher derivatives of $u$ whenever required. These estimates are used in the following regime-wise consistency analysis.\\\\
\noindent{\textbf{Rarefied Regime: $\varepsilon=1$.}} Writing $u^0=u(x,0)$ and $v^0=v(x,0)$, Taylor expansion in time yields
\[
\begin{aligned}
\begin{cases}
u(x,\Delta t)
=u^0+\Delta t\,u_t^0
+\frac{\Delta t^2}{2}u_{tt}^0
+\mathcal O(\Delta t^3),\\[2mm]
v(x,\Delta t)
=v^0+\Delta t\,v_t^0
+\frac{\Delta t^2}{2}v_{tt}^0
+\mathcal O(\Delta t^3).
\end{cases}
\end{aligned}
\]
From the governing equations,
\[
u_t^0=-v_x^0,\qquad
v_t^0=-u_x^0-v^0+au^0.
\]
Differentiating with respect to time gives
\[
\begin{aligned}
u_{tt}^0=u_{xx}^0+v_x^0-au_x^0, \quad v_{tt}^0=v_{xx}^0+u_x^0-av_x^0+v^0-au^0.
\end{aligned}
\]
Hence, the pointwise expansions are
\[
\begin{aligned}
\begin{cases}
u(x,\Delta t)
=u^0-\Delta t\,v_x^0
+\frac{\Delta t^2}{2}
\left(u_{xx}^0+v_x^0-au_x^0\right)
+\mathcal O(\Delta t^3),\\[2mm]
v(x,\Delta t)
=v^0-\Delta t
\left(u_x^0+v^0-au^0\right)+\frac{\Delta t^2}{2}
\left(v_{xx}^0+u_x^0-av_x^0+v^0-au^0\right)
+\mathcal O(\Delta t^3).
\end{cases}
\end{aligned}
\]
Integrating over the cell and using $\Delta t=\lambda h$, we obtain
\begin{equation}\label{rar:exact}
\begin{aligned}
\begin{cases}
\bar u(x_{j+\frac12},\Delta t)
=\bar u^0-\Delta t\,\bar v_x^0
+\frac{\Delta t^2}{2}
\left(\bar u_{xx}^0+\bar v_x^0-a\bar u_x^0\right)
+\mathcal O(h^3),\\[2mm]
\bar v(x_{j+\frac12},\Delta t)
=\bar v^0-\Delta t
\left(\bar u_x^0+\bar v^0-a\bar u^0\right)+\frac{\Delta t^2}{2}
\left(\bar v_{xx}^0+\bar u_x^0-a\bar v_x^0
+\bar v^0-a\bar u^0\right)
+\mathcal O(h^3).
\end{cases}
\end{aligned}
\end{equation}
The numerical update for $u$ is
\begin{equation}\label{u:nn}
\begin{aligned}
\bar u_{j+\frac12}^{\,1}
={}&\frac12(\bar u_j^0+\bar u_{j+1}^0)
+\frac18(\bar u'_j-\bar u'_{j+1})-\lambda(\bar v_{j+1}^{1/2}-\bar v_j^{1/2}).
\end{aligned}
\end{equation}
Substituting the neighboring-cell expansions, reconstruction estimates, and intermediate-stage expansion gives
\begin{equation}\label{rareu:numerical}
\begin{aligned}
\bar u_{j+\frac12}^{\,1}
={}&\bar u^0-\Delta t\,\bar v_x^0+\frac{\Delta t^2}{2}
\left(\bar u_{xx}^0+\bar v_x^0-a\bar u_x^0\right)
+\mathcal O(h^3).
\end{aligned}
\end{equation}
Thus, the numerical and exact expansions agree through second order, provided the intermediate-stage estimate holds.

The numerical update for $v$ is
\[
\begin{aligned}
\bar v_{j+\frac12}^{\,1}
={}&\frac12(\bar v_j^0+\bar v_{j+1}^0)
+\frac18(\bar v'_j-\bar v'_{j+1})-\lambda(\bar u_{j+1}^{1/2}-\bar u_j^{1/2})\\
&+\frac{\Delta t(1-2\alpha)}{2(1-\alpha)}
\left(a\bar u_{j+\frac12}^{1}
-\bar v_{j+\frac12}^{1}\right)+\frac{\Delta t}{4(1-\alpha)}
\left[
a(\bar u_{j+1}^{\alpha}+\bar u_j^{\alpha})
-(\bar v_{j+1}^{\alpha}+\bar v_j^{\alpha})
\right].
\end{aligned}
\]
Using the corresponding stage expansions and reconstruction estimates yields the proposed expansion
\begin{equation}\label{rarev:numerical}
\begin{aligned}
\bar v_{j+\frac12}^{\,1}
={}&\bar v^0-\Delta t
\left(\bar u_x^0+\bar v^0-a\bar u^0\right)+\frac{\Delta t^2}{2}
\left(\bar v_{xx}^0+\bar u_x^0-a\bar v_x^0
+\bar v^0-a\bar u^0\right)
+\mathcal O(h^3).
\end{aligned}
\end{equation}
Comparison with \eqref{rar:exact} establishes consistency through second order if this numerical-stage expansion is verified.\\\\
\noindent{\textbf{Fluid-Dynamic Regime: $\varepsilon\to0.$}} In the fluid-dynamic limit, the relaxation term enforces $v=au$. The limiting equation is therefore
\begin{equation}\label{fluid:eq}
u_t+au_x=0.
\end{equation}
Assuming equilibrium-compatible initial data $v_0=au_0$, Taylor expansion gives
\[
u(x,\Delta t)
=u_0-a\Delta t\,u_x^0
+\frac{a^2\Delta t^2}{2}u_{xx}^0
+\mathcal O(\Delta t^3).
\]
The corresponding exact cell-average expansion is
\begin{equation}\label{fluid:exact}
\begin{aligned}
\bar u(x_{j+\frac12},\Delta t)
={}&\bar u^0-a\Delta t\,\bar u_x^0
+\frac{a^2\Delta t^2}{2}\bar u_{xx}^0
+\mathcal O(h^3).
\end{aligned}
\end{equation}
In this limit, the numerical update for $u$ becomes
\[
\begin{aligned}
\bar u_{j+\frac12}^{\,1}
={}&\frac12(\bar u_j^0+\bar u_{j+1}^0)
+\frac18(\bar u'_j-\bar u'_{j+1})-a\lambda(\bar u_{j+1}^0-\bar u_j^0)
+\frac a2\lambda^2(\bar v'_{j+1}-\bar v'_j).
\end{aligned}
\]
Applying the reconstruction estimates gives
\begin{equation}\label{fluid:numerical}
\begin{aligned}
\bar u_{j+\frac12}^{\,1}
={}&\bar u^0-a\Delta t\,\bar u_x^0
+\frac a2\Delta t^2\bar v_{xx}^0
+\mathcal O(h^3).
\end{aligned}
\end{equation}
Since $v_0=au_0$, we have $\bar v_{xx}^0=a\bar u_{xx}^0$. Consequently,
\[
\bar u_{j+\frac12}^{\,1}
=\bar u^0-a\Delta t\,\bar u_x^0
+\frac{a^2\Delta t^2}{2}\bar u_{xx}^0
+\mathcal O(h^3),
\]
which agrees with \eqref{fluid:exact}. Hence, the limiting update is consistent with the advection equation which retains the order, subject to the stated reconstruction assumptions.\\\\
\noindent{\textbf{Thin-Relaxation Regime: $\varepsilon\ll1$.}}
For small $\varepsilon$, the Chapman-Enskog expansion yields the effective equation
\begin{equation}\label{thin}
u_t+au_x=\varepsilon(1-a^2)u_{xx},
\end{equation}
together with
\begin{equation}\label{thin:v}
v=au-\varepsilon(1-a^2)u_x+\mathcal O(\varepsilon^2).
\end{equation}
We assume well-prepared initial data satisfying
\[
v_0=au_0+\mathcal O(\varepsilon),
\qquad
v_x^0=au_x^0+\mathcal O(\varepsilon),
\]
along with the preceding smoothness and reconstruction assumptions. Setting $D=1-a^2$, we have
\[
u_t=-au_x+\varepsilon D u_{xx}.
\]
Differentiating in time and substituting the effective equation gives
\[
u_{tt}
=a^2u_{xx}-2\varepsilon aD u_{xxx}
+\mathcal O(\varepsilon^2).
\]
Thus,
\begin{equation}\label{thin:pointwise}
\begin{aligned}
u(x,\Delta t)
={}&u^0-a\Delta t\,u_x^0
+\frac{a^2\Delta t^2}{2}u_{xx}^0+\varepsilon D\Delta t
\left(u_{xx}^0-a\Delta t\,u_{xxx}^0\right)+\mathcal O(\Delta t^3,\varepsilon\Delta t^2,
\varepsilon^2).
\end{aligned}
\end{equation}
Integrating over the cell and substituting \(\Delta t=\lambda h\), we obtain the exact expansion
\begin{equation}\label{thin:exact}
\begin{aligned}
\bar u(x_{j+\frac12},\lambda h)
={}&\bar u^0-a\lambda h\,\bar u_x^0
+\frac{a^2\lambda^2h^2}{2}\bar u_{xx}^0+\varepsilon D\lambda h
\left(\bar u_{xx}^0-a\lambda h\,\bar u_{xxx}^0\right)+\mathcal O(h^3,\varepsilon h^2,\varepsilon^2).
\end{aligned}
\end{equation}
Consistency in this regime requires the numerical update to reproduce both the leading-order advection terms and the first-order diffusive correction. To verify this, we expand the intermediate numerical stages and substitute these expansions into the complete numerical update. After collecting the resulting terms up to the required orders, we obtain
\begin{equation}\label{thin:numerical}
\begin{aligned}
\bar u_{j+\frac12}^{\,1}
={}&\bar u^0-a\lambda h\,\bar u_x^0
+\frac{a^2\lambda^2h^2}{2}\bar u_{xx}^0+\varepsilon D\lambda h
\left(\bar u_{xx}^0-a\lambda h\,\bar u_{xxx}^0\right)+\mathcal O(h^3,\varepsilon h^2,\varepsilon^2).
\end{aligned}
\end{equation}
Comparing the exact and numerical expansions, we observe that their corresponding remainder terms are both of order \(\mathcal{O}(h^3,\varepsilon h^2,\varepsilon^2)\). Thus, under the well-prepared initial-data assumptions, the numerical and exact cell-average solutions agree up to the stated asymptotic order in the thin-relaxation regime. This establishes the consistency of the numerical scheme with the Chapman-Enskog approximation.

% In summary, the consistency analysis establishes the agreement between the exact and numerical expansions in the rarefied, fluid-dynamic, and thin-relaxation regimes, subject to the stated assumptions and verification of the corresponding numerical-stage expansions.

We now turn to the asymptotic-preserving (AP) property, which examines whether the numerical scheme correctly recovers the limiting fluid-dynamic model as $\varepsilon\to0$ without imposing a time-step restriction that depends on the relaxation parameter.
\subsection{Asymptotic-Preserving Property}
As \(\varepsilon\to0\), the relaxation term in the system enforces the equilibrium relation $v=au$. Consequently, the limiting equation is
\[
u_t+au_x=0.
\]
To show that the scheme is asymptotic-preserving, we first consider the predictor step. From the second component of the predictor equation, we obtain
\[
\begin{aligned}
\bar v_j^{n+\frac12}-a\bar u_j^{n+\frac12}
={}&\frac{\varepsilon}{\varepsilon+\frac{\Delta t}{2}}
\left[
\bar v_j^n-a\bar u_j^{n+\frac12}
-\frac{\Delta t}{2}\frac{u'_j}{\Delta x}
\right].
\end{aligned}
\]
Assuming that the numerical solution and its reconstructed slopes remain bounded independently of \(\varepsilon\), the expression in brackets remains bounded for fixed \(\Delta t\) and \(\Delta x\). Hence,
\[
\bar v_j^{n+\frac12}-a\bar u_j^{n+\frac12}
=\mathcal{O}(\varepsilon),
\]
and therefore,
\[
\bar v_j^{n+\frac12}-a\bar u_j^{n+\frac12}
\longrightarrow 0
\qquad\text{as }\varepsilon\to0.
\]
Thus, the predictor step satisfies the equilibrium relation in the limit.

Next, define the relaxation residual at the corrector step by
\[
r_{j+\frac12}^{n+1}
=\bar v_{j+\frac12}^{n+1}
-a\bar u_{j+\frac12}^{n+1},
\qquad
r_j^{n+\alpha}
=\bar v_j^{n+\alpha}-a\bar u_j^{n+\alpha}.
\]
The corrector step of the proposed scheme given in \eqref{corrector:general},
\begin{align}
\bar{u}_{j+\frac12}^{\,n+1}
&=
\frac{1}{2}
(\bar{u}_j^n+\bar{u}_{j+1}^n)
+ \frac{1}{8}
(\bar{u}'_j-\bar{u}'_{j+1})
- \lambda
\Bigl[(\bar{v}_{j+1}^{\,n+\frac12})
- (\bar{v}_j^{\,n+\frac12})
\Bigr]\\
\bar{v}_{j+\frac12}^{\,n+1}
&=
\frac{1}{2}
(\bar{v}_j^n+\bar{v}_{j+1}^n)
+ \frac{1}{8}
(\bar{v}'_j-\bar{v}'_{j+1})
- \lambda
\Bigl[(\bar{u}_{j+1}^{\,n+\frac12})
-(\bar{u}_j^{\,n+\frac12})
\Bigr]
\notag\\
&\quad
-\frac{\Delta t}{\varepsilon} \left(\frac{w_1}{2}(\bar{v}_j^{\,n+\alpha} -a \bar{u}_j^{\,n+\alpha}) + \frac{w_1}{2}(\bar{v}_{j+1}^{\,n+\alpha}-a\bar{u}_{j+1}^{\,n+\alpha}) + w_2 (\bar{v}_{j+\frac12}^{\,n+1}-a\bar{u}_{j+\frac12}^{\,n+1}) \right),
\end{align}
with $w_1=\dfrac{1}{2(1-\alpha)}$ and $w_2=\dfrac{1-2\alpha}{2(1-\alpha)}$.
Subtracting \(a\) times the \(u\)-corrector equation from the \(v\)-corrector equation 
% and introducing the relaxation variable
% \[
% r_j^n=\bar{v}_j^n-a\bar{u}_j^n,
% \qquad
% r_{j+\frac12}^{n+1}
% =\bar{v}_{j+\frac12}^{n+1}
% -a\bar{u}_{j+\frac12}^{n+1},
% \]
% we obtain
\begin{align}
r_{j+\frac12}^{n+1}
={}&\frac12(r_j^n+r_{j+1}^n)
+\frac18(r'_j-r'_{j+1})-\lambda\Big[
(\bar{u}_{j+1}^{n+\frac12}-\bar{u}_j^{n+\frac12})
-a(\bar{v}_{j+1}^{n+\frac12}-\bar{v}_j^{n+\frac12})
\Big] \nonumber\\
&-\frac{\Delta t}{\varepsilon}
\left[
\frac{w_1}{2}(r_j^{n+\alpha}+r_{j+1}^{n+\alpha})
+w_2r_{j+\frac12}^{n+1}
\right].
\label{eq:relaxation_corrector}
\end{align}
Defining
\begin{align}
\mathcal{C}_{j+\frac12}^{n+1}
={}&\frac12(r_j^n+r_{j+1}^n)
+\frac18(r'_j-r'_{j+1})-\lambda\Big[
(\bar{u}_{j+1}^{n+\frac12}-\bar{u}_j^{n+\frac12})
-a(\bar{v}_{j+1}^{n+\frac12}-\bar{v}_j^{n+\frac12})
\Big],
\end{align}
and collecting the terms involving \(r_{j+\frac12}^{n+1}\), equation \eqref{eq:relaxation_corrector} becomes
\[
\left(1+\frac{w_2\Delta t}{\varepsilon}\right)
r_{j+\frac12}^{n+1}
=
\mathcal{C}_{j+\frac12}^{n+1}
-\frac{w_1\Delta t}{2\varepsilon}
\left(r_j^{n+\alpha}+r_{j+1}^{n+\alpha}\right).
\]
Substituting
\[
w_1=\frac{1}{2(1-\alpha)},
\qquad
w_2=\frac{1-2\alpha}{2(1-\alpha)},
\]
which yields
\begin{align}
&\left(1+\frac{1-2\alpha}{2(1-\alpha)}
\frac{\Delta t}{\varepsilon}\right)
r_{j+\frac12}^{n+1}=
\mathcal{C}_{j+\frac12}^{n+1}
-\frac{\Delta t}{4(1-\alpha)\varepsilon}
\left(r_j^{n+\alpha}+r_{j+1}^{n+\alpha}\right).
\label{eq:relaxation_final}
\end{align}
% Subtracting \(a\) times the \(u\)-corrector equation from the \(v\)-corrector equation gives
% \[
% \begin{aligned}
% \left(1+\frac{1-2\alpha}{2(1-\alpha)}
% \frac{\Delta t}{\varepsilon}\right)
% r_{j+\frac12}^{n+1}
% ={}&\mathcal{C}_{j+\frac12}^{n+1}-\frac{\Delta t}{4(1-\alpha)\varepsilon}
% \left(r_j^{n+\alpha}+r_{j+1}^{n+\alpha}\right),
% \end{aligned}
% \]
% where \(\mathcal{C}_{j+\frac12}^{n+1}\) collects the non-stiff terms in the corrector equation. 
Under the same boundedness assumption, these terms satisfy
\[
\mathcal{C}_{j+\frac12}^{n+1}=\mathcal{O}(1).
\]
Moreover, applying the predictor relation at the intermediate stage \(n+\alpha\) gives
\[
r_j^{n+\alpha}=\mathcal{O}(\varepsilon),
\qquad
r_{j+1}^{n+\alpha}=\mathcal{O}(\varepsilon),
\]
for fixed \(\alpha\in(0,\frac12)\), provided the corresponding predictor terms remain bounded. Thus, the right-hand side of the residual equation is bounded independently of \(\varepsilon\). Since
\[
\frac{1-2\alpha}{2(1-\alpha)}>0
\qquad\text{for }0<\alpha<\frac12,
\]
we obtain, for fixed \(\Delta t>0\),
\[
r_{j+\frac12}^{n+1}=\mathcal{O}(\varepsilon)
\longrightarrow0
\qquad\text{as }\varepsilon\to0.
\]
Consequently,
\[
\bar v_{j+\frac12}^{n+1}
\longrightarrow
a\bar u_{j+\frac12}^{n+1}.
\]
Finally, substituting the equilibrium relation into the \(u\)-corrector equation gives
\[
\begin{aligned}
\bar u_{j+\frac12}^{n+1}
={}&\frac12(\bar u_j^n+\bar u_{j+1}^n)
+\frac18(\bar u'_j-\bar u'_{j+1})-a\lambda
\left(\bar u_{j+1}^{n+\frac12}
-\bar u_j^{n+\frac12}\right),
\end{aligned}
\]
which is the corresponding corrector equation for the limiting advection equation.

Hence, for fixed time and space steps, the scheme reduces to a discretization of the limiting equation as \(\varepsilon\to0\). Therefore, the scheme possesses the asymptotic-preserving property.
\subsection{Stability Analysis}
We next investigate the stability properties of the proposed scheme. To this end, we consider evolution problems comprising stiff and non-stiff contributions. Specifically, let
\begin{equation}
\mathrm{y}'=\mathrm{f}(\mathrm{y})+\mathrm{g}(\mathrm{y}),
\label{eq:split_system}
\end{equation}
subject to the initial condition
\begin{equation}
\mathrm{y}(0)=\mathrm{y}_0,
\label{eq:split_initial}
\end{equation}
where $\mathrm{f}(\mathrm{y})$ and $\mathrm{g}(\mathrm{y})$ denote the non-stiff and stiff contributions, respectively. This decomposition is particularly relevant to hyperbolic balance laws with relaxation, in which the source term may introduce time scales substantially smaller than those associated with the transport dynamics.

For such problems, an effective time discretizations should accommodate the stiffness of the source contribution without imposing unnecessarily restrictive time-step conditions. Moreover, evaluating the non-stiff contribution may involve computationally expensive operations, making its explicit treatment desirable whenever appropriate. These considerations motivate a partitioned time discretizations in which the non-stiff term $\mathrm{f}(\mathrm{y})$ is treated explicitly, while the stiff term $\mathrm{g}(\mathrm{y})$ is treated implicitly. This explicit-implicit treatment forms the basis for the stability analysis developed below.
\par The numerical scheme can be formulated as the following system
of equations:
\begin{align}
\mathrm{y}_{\frac12} &=\mathrm{{y}}_0+\frac{\mathrm{k}}{2}\mathrm{f}(\mathrm{y}_0) +\frac{\mathrm{k}}{2}\mathrm{g}\left(\mathrm{y}_{\frac12}\right),\notag\\[4pt]
\mathrm{y}_{\alpha} &=\mathrm{y}_0+\alpha \mathrm{k}\mathrm{f}(\mathrm{y}_0) +\alpha \mathrm{k}\mathrm{g}\left(\mathrm{y}_{\alpha}\right),\notag\\[4pt]
\mathrm{y}_1 &=\mathrm{y}_0+\mathrm{k}\mathrm{f}\left(\mathrm{y}_{\frac12}\right) +\frac{\mathrm{k}}{2(1-\alpha)}\mathrm{g}\left(\mathrm{y}_{\alpha}\right) +\frac{\mathrm{k}(1-2\alpha)}{2(1-\alpha)}\mathrm{g}(\mathrm{y}_1).
\end{align}
Here, $\mathrm{k}$ represents the time-step size, $\mathrm{y}_0$ denotes the initial numerical solution, and $\mathrm{y}_1$ is the approximation after one time step. Standard local truncation error analysis gives
\[
\mathrm{y}(\mathrm{k})-\mathrm{y}_1=\mathcal{O}(\mathrm{k}^3).
\]
Consequently, the method is second-order accurate.
To investigate the stability properties of the scheme, we
consider the linear test problem
\[
\begin{aligned}
&\mathrm{y}'=\lambda_1\mathrm{y}+\lambda_2\mathrm{y},\\
&\mathrm{y}(0)=1.
\end{aligned}
\]
Introducing the dimensionless parameters
\[
\mathrm{z}_1=\lambda_1 \mathrm{k}, \;\; \mathrm{z}_2= \lambda_2\mathrm{k},
\]
the intermediate stages are given by
\[
\begin{aligned}
\mathrm{y}_{\frac12}& = \mathrm{y}_0\frac{2+\mathrm{z}_1}{2-\mathrm{z}_2},\\[4pt]
\mathrm{y}_{\alpha}&=\mathrm{y}_0\frac{1+\alpha \mathrm{z}_1}{1-\alpha \mathrm{z}_2}.
\end{aligned}
\]
Substituting these expressions into the final stage yields
\[
\mathrm{y}_1=\left(\frac{1+\dfrac{\mathrm{z}_1(2+\mathrm{z}_1)}{2-\mathrm{z}_2}+\dfrac{\mathrm{z}_2(1+\alpha \mathrm{z}_1)}{2(1-\alpha)(1-\alpha \mathrm{z}_2)}}{1-\dfrac{\mathrm{z}_2(1-2\alpha)}{2(1-\alpha)}}\right)\mathrm{y}_0.
\]
Therefore, the amplification factor is
\[
R(\mathrm{z}_1,\mathrm{z}_2)= \frac{1+\dfrac{\mathrm{z}_1(2+\mathrm{z}_1)}{2-\mathrm{z}_2}+\dfrac{\mathrm{z}_2(1+\alpha \mathrm{z}_1)}{2(1-\alpha)(1-\alpha \mathrm{z}_2)}}{1-\dfrac{\mathrm{z}_2(1-2\alpha)}{2(1-\alpha)}}.
\]
To examine the stability associated with the $\mathrm{g}$-component, we set $\mathrm{f}\equiv0$, which corresponds to $\mathrm{z}_1=0$. The resulting stability function is
\[
R(0,\mathrm{z}_2)=\frac{1+\dfrac{\mathrm{z}_2}{2(1-\alpha)(1-\alpha \mathrm{z}_2)}}{1-\dfrac{\mathrm{z}_2(1-2\alpha)}{2(1-\alpha)}}.
\]
Combining the terms in the numerator and denominator gives
\[
R(0,\mathrm{z}_2)=\frac{2(1-\alpha)+\mathrm{z}_2\left(1-2\alpha(1-\alpha)\right)}{\left(2(1-\alpha)-\mathrm{z}_2(1-2\alpha)\right)(1-\alpha \mathrm{z}_2)}.
\]
For simplicity, let $\mathrm{z}=\mathrm{z}_2$. The stability function can be expressed in rational form as
\begin{equation}
R(\mathrm{z})=\frac{\mathrm{N}(\mathrm{z})}{\mathrm{D}(\mathrm{z})},
\end{equation}
where
\begin{align}
\mathrm{N}(\mathrm{z})&=2(1-\alpha)+\mathrm{z}\left(1-2\alpha(1-\alpha)\right),\\
\mathrm{D}(\mathrm{z})&=\left(2(1-\alpha)-\mathrm{z}(1-2\alpha)\right)(1-\alpha \mathrm{z}).
\end{align}
To investigate the stability properties along the imaginary axis, we set
\begin{equation}
\mathrm{z}=i\mathrm{y},\quad \mathrm{y}\in\mathbb{R}.
\end{equation}
Consequently,
\begin{equation}
R(i\mathrm{y})=\frac{\mathrm{N}(i\mathrm{y})}{\mathrm{D}(i\mathrm{y})},
\end{equation}
and, provided that $D(iy)\neq 0$, we obtain
\begin{equation}
|R(i\mathrm{y})|^2=\frac{|\mathrm{N}(i\mathrm{y})|^2}{|\mathrm{D}(i\mathrm{y})|^2}.
\end{equation}
Substituting $\mathrm{z}=i\mathrm{y}$ into the numerator gives
\begin{equation*}
\mathrm{N}(i\mathrm{y})=2(1-\alpha) +i\mathrm{y}\left(1-2\alpha(1-\alpha)\right).
\end{equation*}
Taking the squared modulus, we obtain
\begin{equation*}
|\mathrm{N}(i\mathrm{y})|^2=4(1-\alpha)^2 +\mathrm{y}^2\left(1-2\alpha+2\alpha^2\right)^2.
\end{equation*}
Similarly, the denominator becomes
\begin{equation*}
\mathrm{D}(i\mathrm{y})= \left(2(1-\alpha)-i\mathrm{y}(1-2\alpha)\right)(1-i\alpha \mathrm{y}).
\end{equation*}
Its squared modulus is therefore
\begin{equation*}
|\mathrm{D}(i\mathrm{y})|^2 = \left[4(1-\alpha)^2+\mathrm{y}^2(1-2\alpha)^2\right] \left(1+\alpha^2\mathrm{y}^2\right).
\end{equation*}
Comparing the squared moduli of the denominator and numerator yields
\begin{align*}
|\mathrm{D}(i\mathrm{y})|^2-|\mathrm{N}(i\mathrm{y})|^2 &=\left[4(1-\alpha)^2+\mathrm{y}^2(1-2\alpha)^2\right] \left(1+\alpha^2 \mathrm{y}^2\right)- \left[4(1-\alpha)^2 +\mathrm{y}^2\left(1-2\alpha(1-\alpha)\right)^2\right] \nonumber\\
& = \alpha ^2 (1-2 \alpha)^2 \mathrm{y}^4.
\end{align*}
For the parameter range $0<\alpha<\frac12$, the imaginary-axis
stability condition can be established by verifying that
\begin{equation}
|\mathrm{D}(i\mathrm{y})|^2-|\mathrm{N}(i\mathrm{y})|^2\geq 0,
\; \text{for all} \;\; y \in\mathbb{R}.
\end{equation}
Whenever this inequality holds and $\mathrm{D}(i\mathrm{y})\neq 0$, it follows that
\begin{equation}
|R(i\mathrm{y})|\leq 1.
\end{equation}
At $\mathrm{y}=0$, we have
\begin{equation}
R(0)=1.
\end{equation}
Next, the poles of the stability function are determined by solving $\mathrm{D}(\mathrm{z})=0$, that is,
\begin{equation*}
\left(2(1-\alpha)-\mathrm{z}(1-2\alpha)\right)(1-\alpha \mathrm{z})=0.
\end{equation*}
For $0<\alpha<\frac12$, the poles are
\begin{equation}
\mathrm{z}_1=\frac{1}{\alpha},
\quad
\mathrm{z}_2=\frac{2(1-\alpha)}{1-2\alpha}.
\end{equation}
Both poles are positive real numbers, since
\begin{equation}
\frac{1}{\alpha}>0,
\quad
\frac{2(1-\alpha)}{1-2\alpha}>0.
\end{equation}
Therefore, the stability function has no poles in the closed left half-plane and is analytic in a neighborhood of that region. A numerical method is said to be $A$-stable if its stability function satisfies
\begin{equation}
|R(\mathrm{z})|\leq 1,
\quad \operatorname{Re}(\mathrm{z})\leq 0,
\end{equation}
and has no poles in the closed left half-plane. Although the preceding analysis establishes the stability inequality along the imaginary axis, an additional argument is required to extend this result to the entire closed left half-plane. We establish this extension using the maximum modulus principle. For \(M>0\), consider the left half-disk
\begin{equation}
\Theta_M=\left\{\mathrm{z}\in\mathbb{C}:
\operatorname{Re}(\mathrm{z})<0,\ |\mathrm{z}|<M\right\}.
\end{equation}
Here, $\overline{\Theta_M}$ denotes the closure of $\Theta_M$, including its boundary. For \(0<\alpha<\frac12\), the stability function \(R(\mathrm{z})\) has no poles in the closed left half-plane and is therefore analytic in a neighbourhood of \(\overline{\Theta_M}\). By the maximum modulus principle, its maximum modulus over \(\overline{\Theta_M}\) is attained on the boundary, which consists of a segment of the imaginary axis and a semicircular arc in the left half-plane. On the imaginary-axis segment, the previously established inequality gives
\begin{equation}
\label{modulus}
|R(i\mathrm{y})|\leq 1,\qquad \mathrm{y}\in\mathbb{R}.
\end{equation}
Furthermore, the numerator of \(R(\mathrm{z})\) has degree at most one, whereas its denominator has degree two, with nonzero leading coefficient \(\alpha(1-2\alpha)\). Consequently,
\begin{equation}
\lim_{|\mathrm{z}|\to\infty}R(\mathrm{z})=0.
\end{equation}
Thus, the maximum modulus of \(R(\mathrm{z})\) on the semicircular arc tends to zero as \(M\to\infty\). Applying the maximum modulus principle and taking this limit yields
\begin{equation*}
|R(\mathrm{z})|\leq 1,\qquad \operatorname{Re}(\mathrm{z})\leq 0.
\end{equation*}
Hence, the stability function is A-stable for \(0<\alpha<\frac12\), provided the imaginary-axis stability inequality holds. A numerical method is L-stable if it is A-stable and its stability function satisfies
\begin{equation*}
\lim_{|\mathrm{z}|\to\infty}R(\mathrm{z})=0.
\end{equation*}
For the present scheme, the stability function can be written as
\begin{equation}
R(\mathrm{z})=
\frac{2(1-\alpha)+\mathrm{z}\left[1-2\alpha(1-\alpha)\right]}
{\left[2(1-\alpha)-(1-2\alpha)\mathrm{z}\right](1-\alpha\mathrm{z})}.
\end{equation}
For \(0< \alpha < \frac12 \), the numerator is a polynomial of degree at most one, whereas the denominator is a polynomial of degree two with leading coefficient \(\alpha(1-2\alpha)\neq0\). Therefore,
\begin{align*}
\lim_{|\mathrm{z}|\to\infty}R(\mathrm{z})
&=\lim_{|\mathrm{z}|\to\infty}
\frac{\mathrm{z}\left[1-2\alpha(1-\alpha)\right]+\mathcal{O}(1)}
{\alpha(1-2\alpha)\mathrm{z}^2+\mathcal{O}(\mathrm{z})}
\\
&=0.
\end{align*}
It follows that the source-only stability function is L-stable for \(0<\alpha<\frac12\).

We next investigate the absolute stability region of the fully coupled transport-source scheme. Let \(\mathrm{z}_1\) and \(\mathrm{z}_2\) denote the dimensionless parameters associated with the transport and source contributions, respectively. The two-parameter absolute stability region is defined by
\begin{equation}
\Theta_A=
\left\{(\mathrm{z}_1,\mathrm{z}_2)\in\mathbb{C}^2:
|R(\mathrm{z}_1,\mathrm{z}_2)|\leq1\right\}.
\end{equation}
In general, this region need not contain the entire product \(\mathbb{C}^{-}\times\mathbb{C}^{-}\), where
\begin{equation*}
\mathbb{C}^{-}=
\left\{\mathrm{z}\in\mathbb{C}:
\operatorname{Re}(\mathrm{z})\leq0\right\}.
\end{equation*}
We therefore seek the region \(\Theta_1\) in the transport parameter \(\mathrm{z}_1\) for which the method remains stable for every \(\mathrm{z}_2\in\mathbb{C}^{-}\). Specifically, we define
\begin{equation}
\Theta_1=
\left\{\mathrm{z}_1\in\mathbb{C}:
\sup_{\mathrm{z}_2\in\mathbb{C}^{-}}
|R(\mathrm{z}_1,\mathrm{z}_2)|\leq1\right\}.
\end{equation}
For each fixed \(\mathrm{z}_1\), provided \(R(\mathrm{z}_1,\mathrm{z}_2)\) is analytic with respect to \(\mathrm{z}_2\) in the left half-plane and continuous up to its boundary, and tends to zero as \(|\mathrm{z}_2|\to\infty\) within that half-plane, the maximum modulus principle allows the supremum over the left half-plane to be determined from the imaginary-axis boundary. Under these conditions,
\begin{equation}
\Theta_1=
\left\{\mathrm{z}_1\in\mathbb{C}:
\sup_{\mathrm{y}\in\mathbb{R}}
|R(\mathrm{z}_1,i\mathrm{y})|\leq1\right\}.
\end{equation}
Accordingly, the boundary locus of this region is characterized by
\begin{equation}
\partial\Theta_1=
\left\{\mathrm{z}_1\in\mathbb{C}:
\sup_{\mathrm{y}\in\mathbb{R}}
|R(\mathrm{z}_1,i\mathrm{y})|=1\right\},
\end{equation}
with the boundary characterization understood under the regularity conditions stated above. This formulation enables the stability region in the transport parameter to be evaluated while accounting for all source parameters in the closed left half-plane.
\subsection{Selection of the Intermediate-Stage Parameter $\alpha$} \label{subsec:alpha_selection}
The consistency and stability analyses yield the common admissible range \(0< \alpha < \frac12\). Although all values within this interval satisfy the derived analytical constraints, the choice of $\alpha$ influences the geometry and extent of the stability regions. Therefore, we examine both the implicit and explicit stability properties to identify a suitable parameter value for the numerical experiments.

\begin{figure}[htb]
    \centering
    \includegraphics[scale=0.48, trim=2.0cm 0cm 9.6cm 0cm, clip]{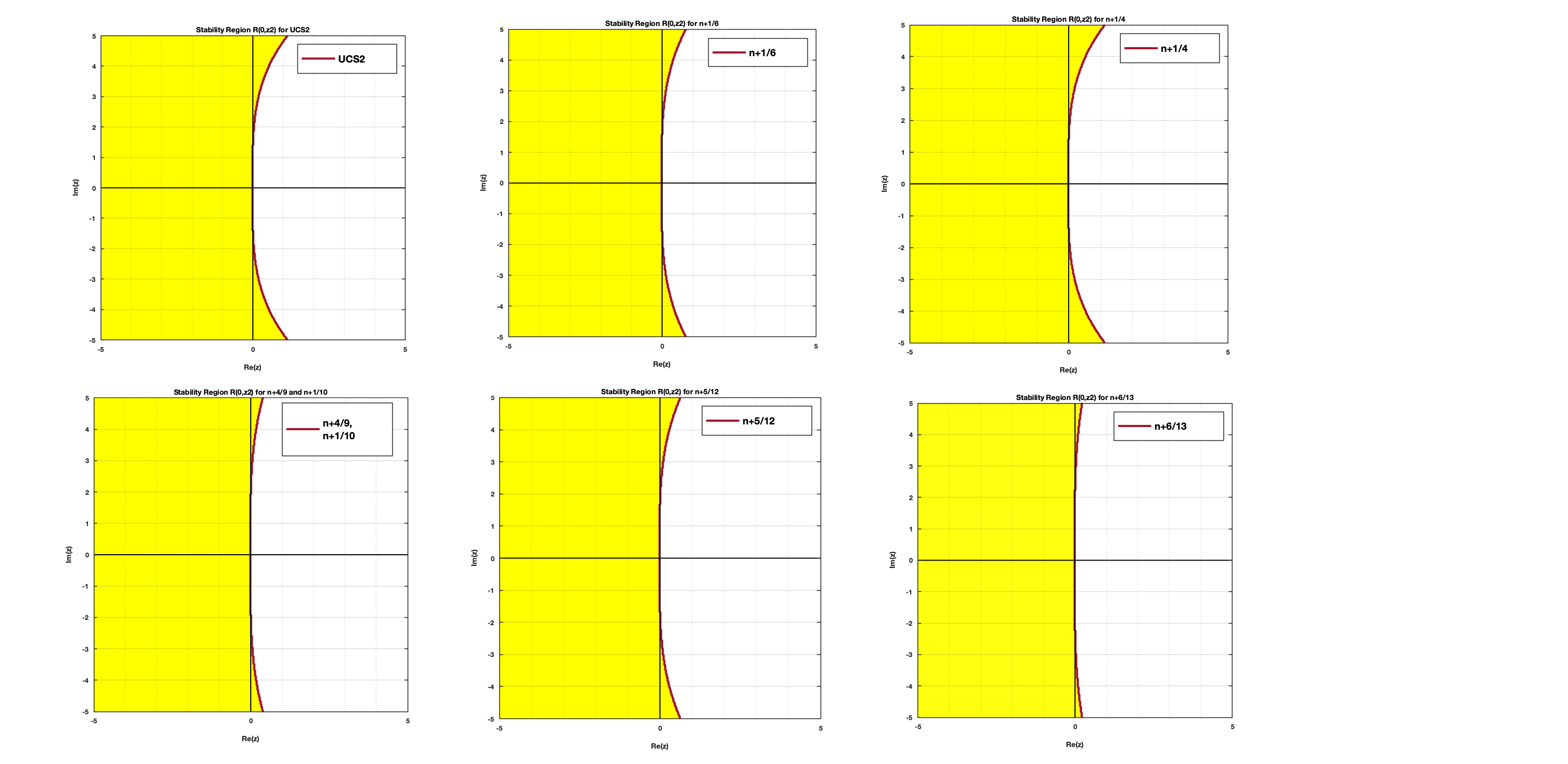}
    \caption{Stability regions \(R(0,\mathrm{z}_2)\) for UCS2 and the proposed generalized scheme corresponding to different values of \(\alpha\), namely \(\alpha \in \{1/6, 1/4, 1/10, 4/9, 5/12, 6/13\}\). The yellow-shaded regions represent the stability domain for each case.}
    \label{stabilty:implicit}
\end{figure}

\begin{figure}[htb]
    \centering
    \begin{subfigure}[b]{0.49\textwidth}
        \centering        \includegraphics[width=\linewidth,height=7.3cm,
        trim=12.0cm 0cm 6.5cm 0cm,clip,
        keepaspectratio]{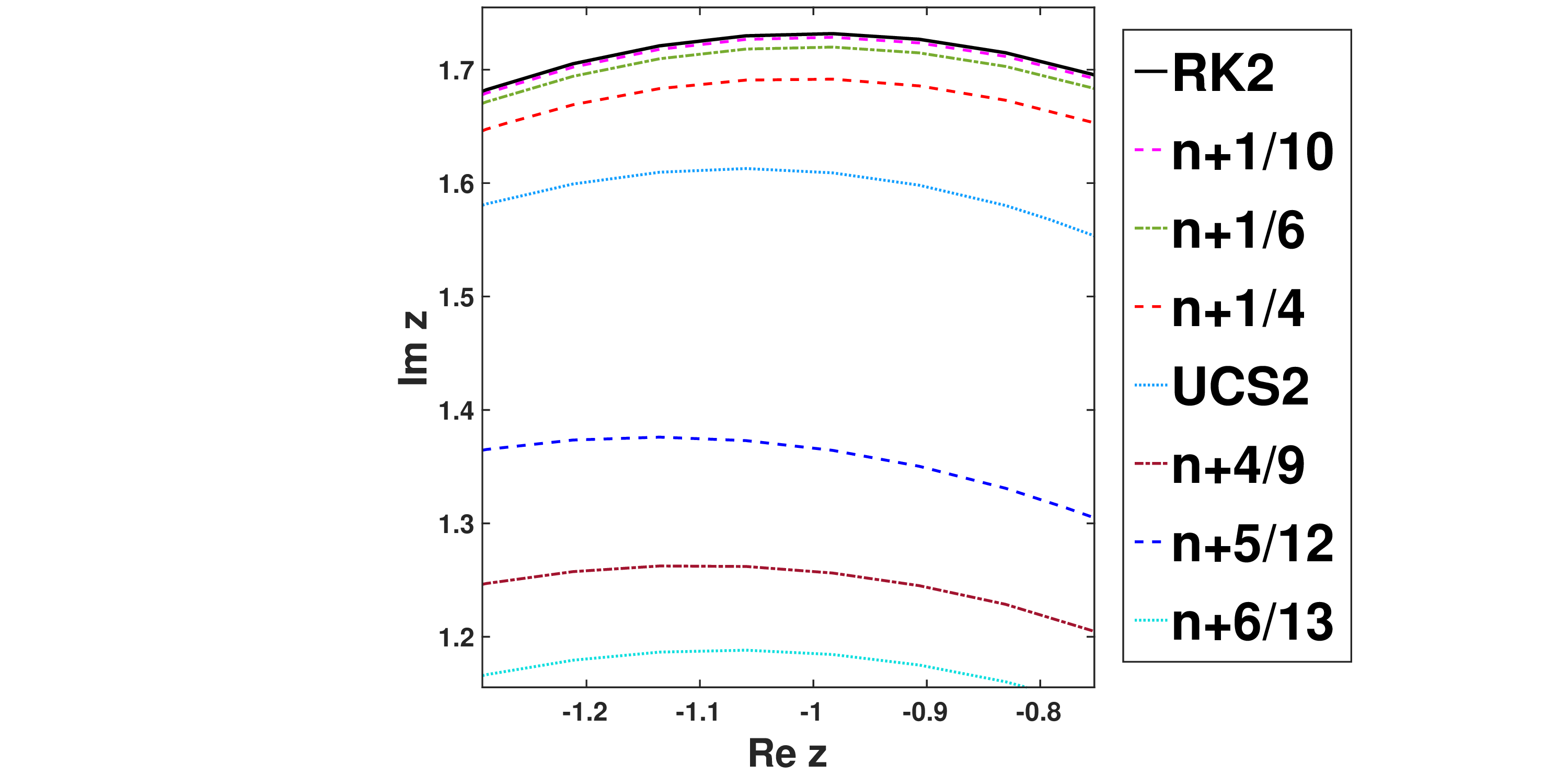}
        \caption{}
    \end{subfigure}
    \hfill
    \begin{subfigure}[b]{0.49\textwidth}
        \centering
        \includegraphics[width=\linewidth,height=7.3cm,
        trim=12.0cm 0cm 6.5cm 0cm,clip,
        keepaspectratio]{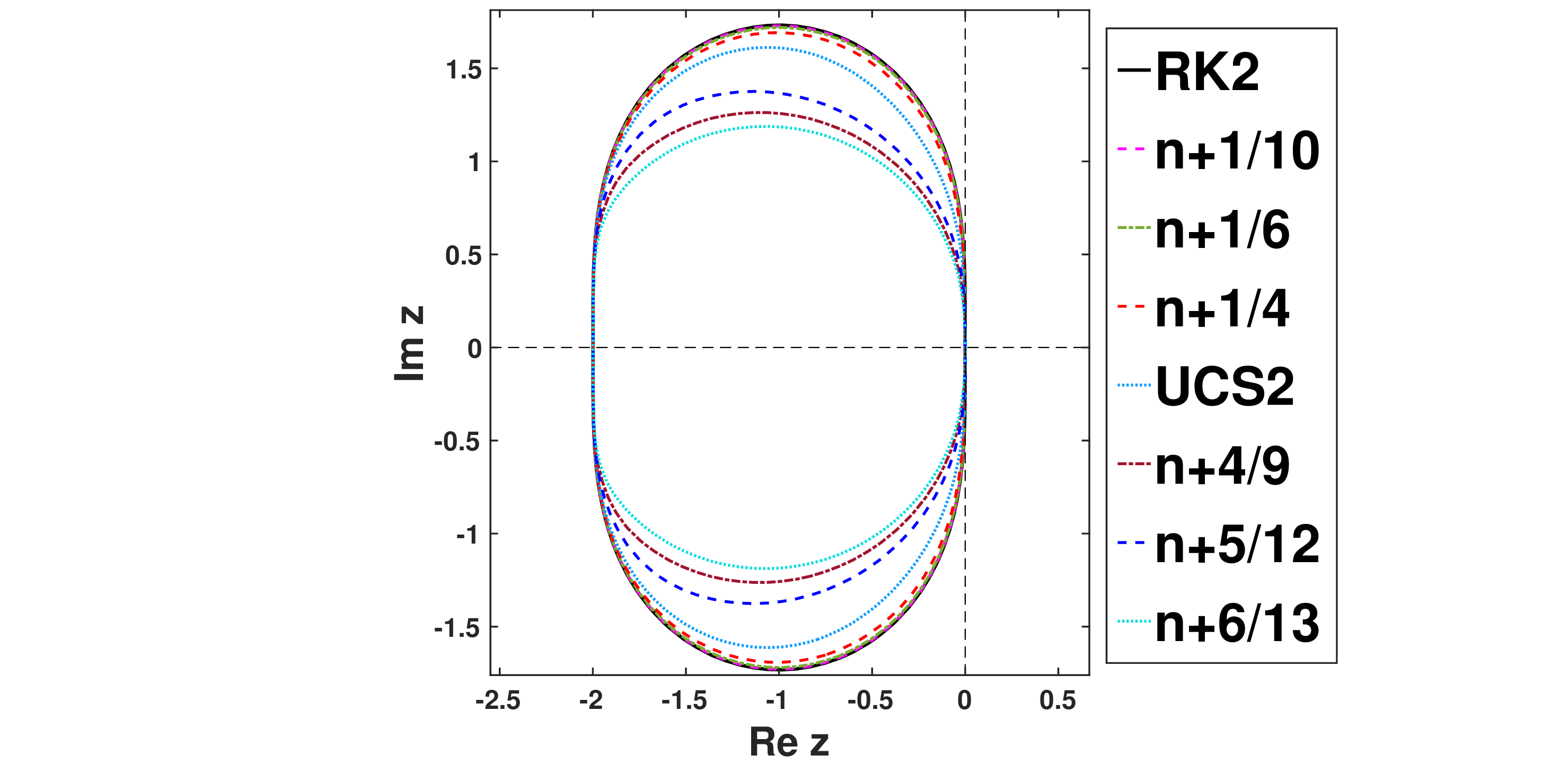}
        \caption{}
    \end{subfigure}
    \caption{Stability regions for the explicit parameter: (a) zoomed-in view of the stability region; (b) complete stability region.}
    \label{stability:explicit_zoomed}
\end{figure}

\begin{figure}[htb]
         \centering
    \includegraphics[scale=0.56, trim=2.0cm 0cm 20.0cm 0cm, clip]{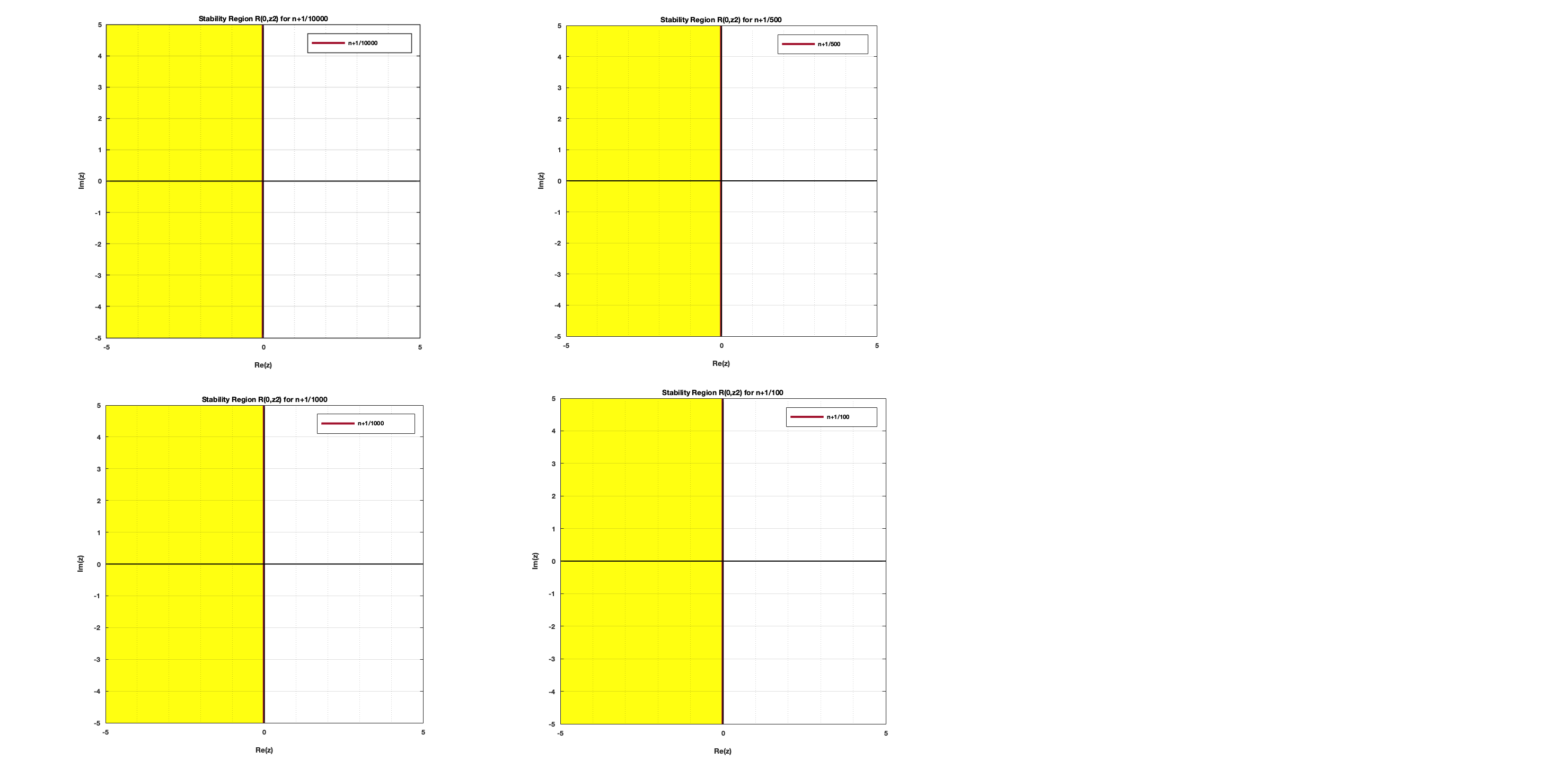}
    \caption{Stability regions \(R(0,\mathrm{z}_2)\) corresponding to \(\alpha \in \{1/100,\,1/500,\,1/1000,\,1/10000\}\). The yellow-shaded areas represent the stability regions for each value of \(\alpha\).}
     \label{stability:4implicitt}
\end{figure}

\begin{figure}[htbp]
    \centering
    \begin{subfigure}[b]{0.32\textwidth}
        \centering
        \includegraphics[width=\linewidth,height=6.0cm,
        trim=11.4cm 0cm 13.1cm 0cm,clip]{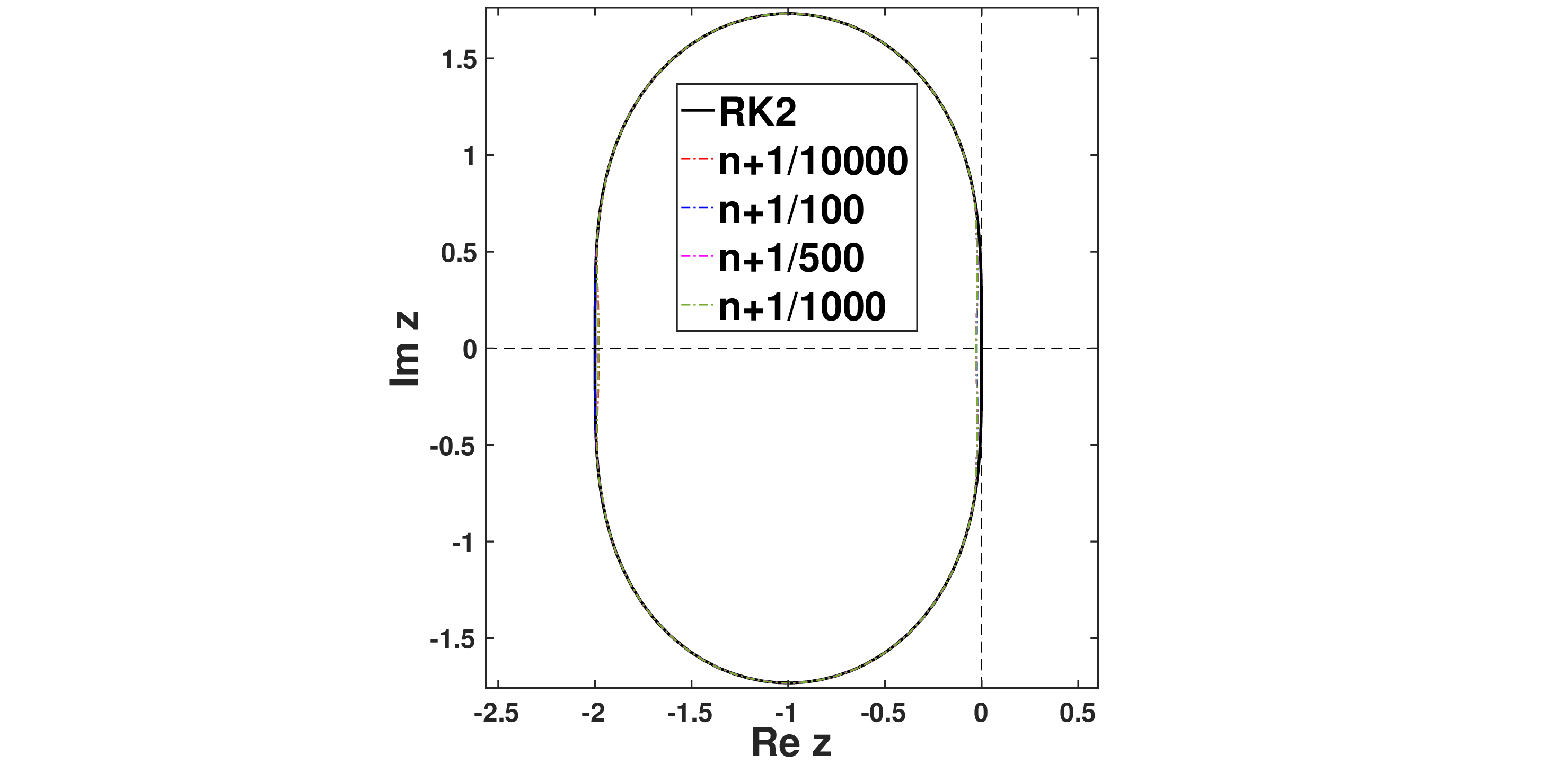}
        \caption{}
    \end{subfigure}
    \hfill
    \begin{subfigure}[b]{0.32\textwidth}
        \centering
        \includegraphics[width=\linewidth,height=6.0cm,
        trim=11.3cm 0cm 13.8cm 0cm,clip]{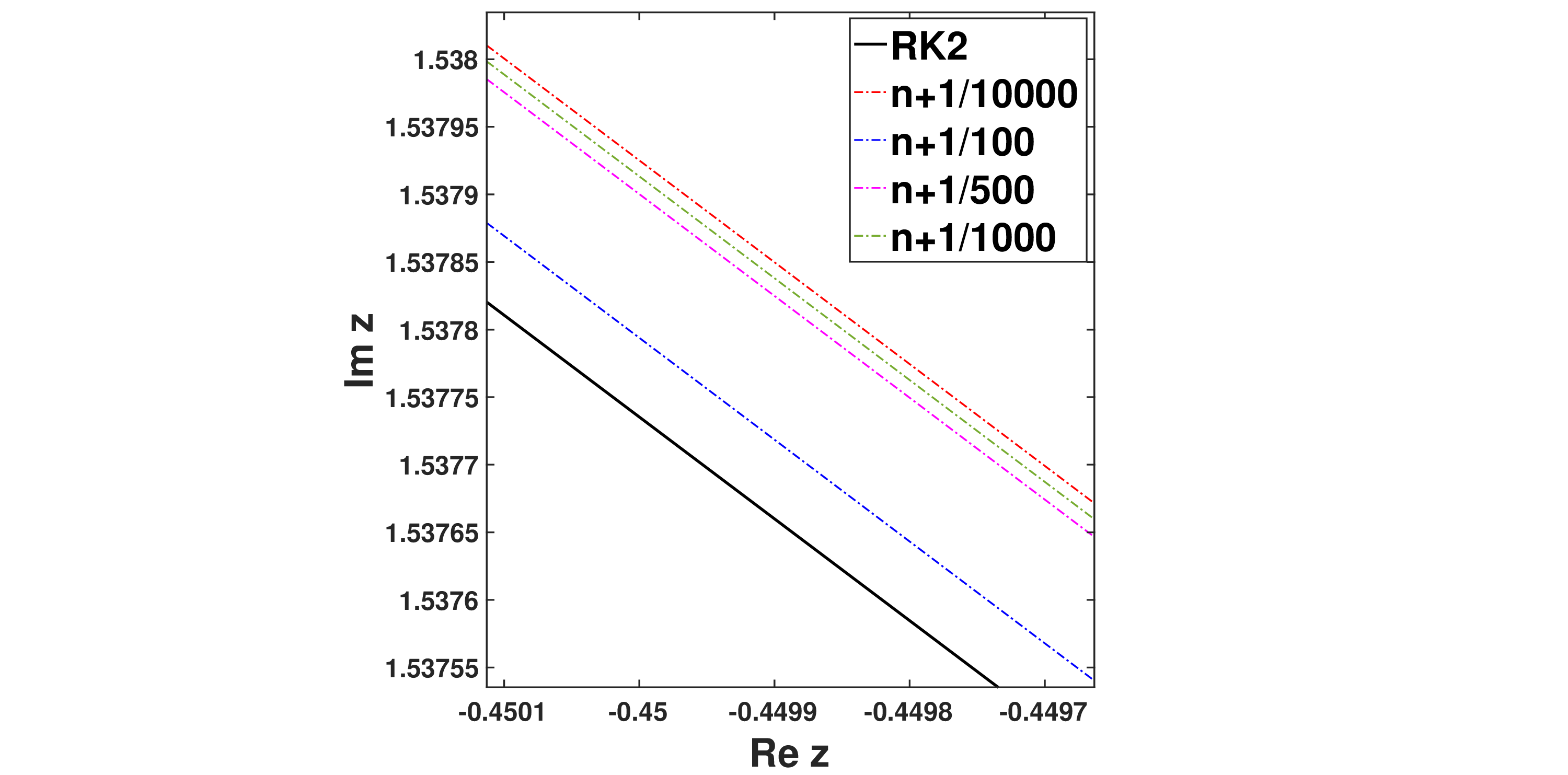}
        \caption{}
    \end{subfigure}
    \hfill
    \begin{subfigure}[b]{0.32\textwidth}
        \centering       \includegraphics[width=\linewidth,height=6.0cm,
        trim=12.0cm 0cm 13.5cm 0cm,clip]{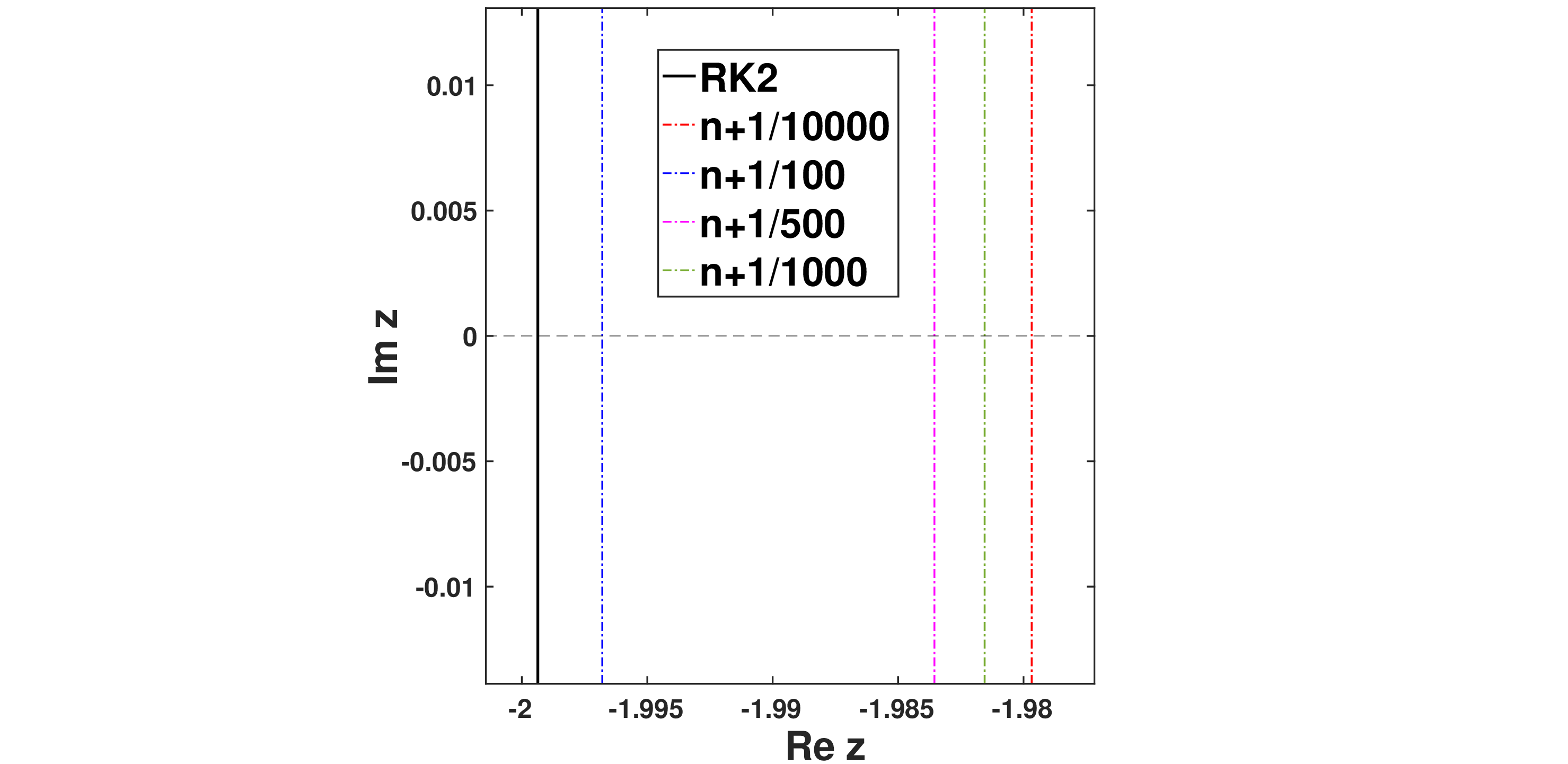}
        \caption{}
    \end{subfigure}
    \caption{Stability regions of the scheme for parameters \(1/100\), \(1/500\), \(1/1000\), and \(1/10000\): (a) stability region for the explicit parameter; (b) zoomed-in view of the stability region; (c) further zoomed-in view of the stability region.}
    \label{stability:explicitn4}
\end{figure}

Figure~\ref{stabilty:implicit} illustrates the stability regions associated with the source-only stability function \(R(0,\mathrm{z}_2)\) for \(\alpha\in\{\frac{1}{10},\frac{1}{6},\frac{1}{4},\frac{4}{9},\frac{5}{12},\frac{6}{13}\}\). The corresponding regions for smaller values of \(\alpha\), namely \(\frac{1}{100},\frac{1}{500},\frac{1}{1000}\), and \(\frac{1}{10000}\), are presented in Figure~\ref{stability:4implicitt}, where the yellow-shaded areas indicate the stability regions. The comparison shows that very small values of \(\alpha\) yield comparatively smaller implicit stability regions, whereas the selected values ranging from \(1/10\) to \(6/13\) produce more extensive regions. To complement this analysis, Figures~\ref{stability:explicit_zoomed} and~\ref{stability:explicitn4} illustrate the previously defined region \(\Theta_1\) for the selected values of \(\alpha\), together with the stability region of the classical second-order explicit Runge-Kutta method (RK2) for comparison. The results indicate that the region approaches the RK2 stability region as \(\alpha\to0\). 
Considering the implicit and explicit stability characteristics, together with the numerical performance, we select $\alpha=1/10$
as the representative parameter for the subsequent numerical investigations. This choice is motivated by the following observations:
\begin{itemize}
    \item The implicit stability region obtained for
    $\alpha=1/10$ is considerably more extensive than those corresponding to very small values of $\alpha$.
    \item Among the selected parameter values, the explicit stability region $\Theta_1$ for $\alpha=1/10$ exhibits a closer resemblance to the classical RK2 stability region.
    \item The comparative one-dimensional numerical results presented in Section~\ref{1d:test} further support this choice, with $\alpha=1/10$ exhibiting favourable numerical performance among the tested parameter values.
\end{itemize}
The selected value $\alpha=1/10$ is subsequently used in the detailed numerical investigations, including the two-dimensional problem considered in Section~\ref{sec:eulerxinjin2d}.

% Considering both stability characteristics, we select \(\alpha=1/10\) as the representative parameter for demonstrating the accuracy and convergence order in Section \ref{accuracy} of the proposed scheme. This choice offers a favourable compromise between the extent of the implicit stability region and the similarity of \(\Theta_1\) to the RK2 stability region. The same value is adopted for the two-dimensional problem in Section~\ref{sec:eulerxinjin2d}, ensuring consistency across the numerical experiments.


\section{Numerical Results}\label{sec:Numerical}
This section assesses the accuracy, convergence behavior, and stability of the proposed numerical scheme through a series of one- and two-dimensional numerical experiments. The performance of the proposed scheme is investigated for different values of the generalized parameter $\alpha$ and compared primarily with the UCS2 central scheme. In addition, the IMEX-RK2 scheme is employed on a finer spatial grid to obtain a reference solution for assessing the accuracy of the computed numerical results. Unless otherwise specified, all simulations are performed using a Courant-Friedrichs-Lewy (CFL) number of $0.9$, with the time-step size determined by
\[
\Delta t = \frac{\mathrm{CFL}\,\Delta x}{\lambda_{\max}},
\]
where $\lambda_{\max}$ denotes the maximum characteristic wave speed. The test problems include smooth and discontinuous initial conditions and cover different relaxation regimes and boundary conditions. 
%Wherever available, the computed solutions are compared with exact solutions, reference solutions, or established benchmark results to assess the accuracy and reliability of the proposed scheme.
\subsection{Accuracy and Convergence Assessment}\label{accuracy}
The accuracy and convergence properties of the proposed numerical scheme are investigated using smooth test problems in one and two spatial dimensions. Since exact solutions are available for both problems using Fourier transforms, the numerical errors can be evaluated directly, allowing the observed convergence rates to be assessed under successive mesh refinement. The two test problems are described below, followed by the computational setup and the corresponding convergence results.
\subsubsection{One-Dimensional Smooth Test Problem}
We consider the one-dimensional Jin-Xin relaxation system on the domain $[0,1]$:
\begin{equation}\label{Xin_Jin_model}
\begin{cases}
u_t+v_x=0,\\[1mm]
v_t+u_x=-\dfrac{1}{\varepsilon}(v-au).
\end{cases}
\end{equation}
Here, $u=u(x,t)$ denotes the macroscopic conserved variable, $v=v(x,t)$ is the associated auxiliary relaxation variable, $\varepsilon>0$ is the relaxation parameter, and $a$ is a prescribed constant. We consider the following smooth, well-prepared initial data:
\begin{equation}\label{smoothdata}
\begin{aligned}
\begin{cases}
u(x,0)&=\sin(2\pi x),\\
v(x,0)&=a\,u(x,0),
\end{cases}
\end{aligned}
\end{equation}
where $a=0.7$.
\subsubsection{Two-Dimensional Smooth Test Problem}
\begin{figure}[htb]
\centering
\begin{subfigure}[b]{0.49\textwidth}
\centering
\includegraphics[width=\textwidth,height=7.3cm,trim=9.0cm 0cm 13.0cm 0cm,clip]{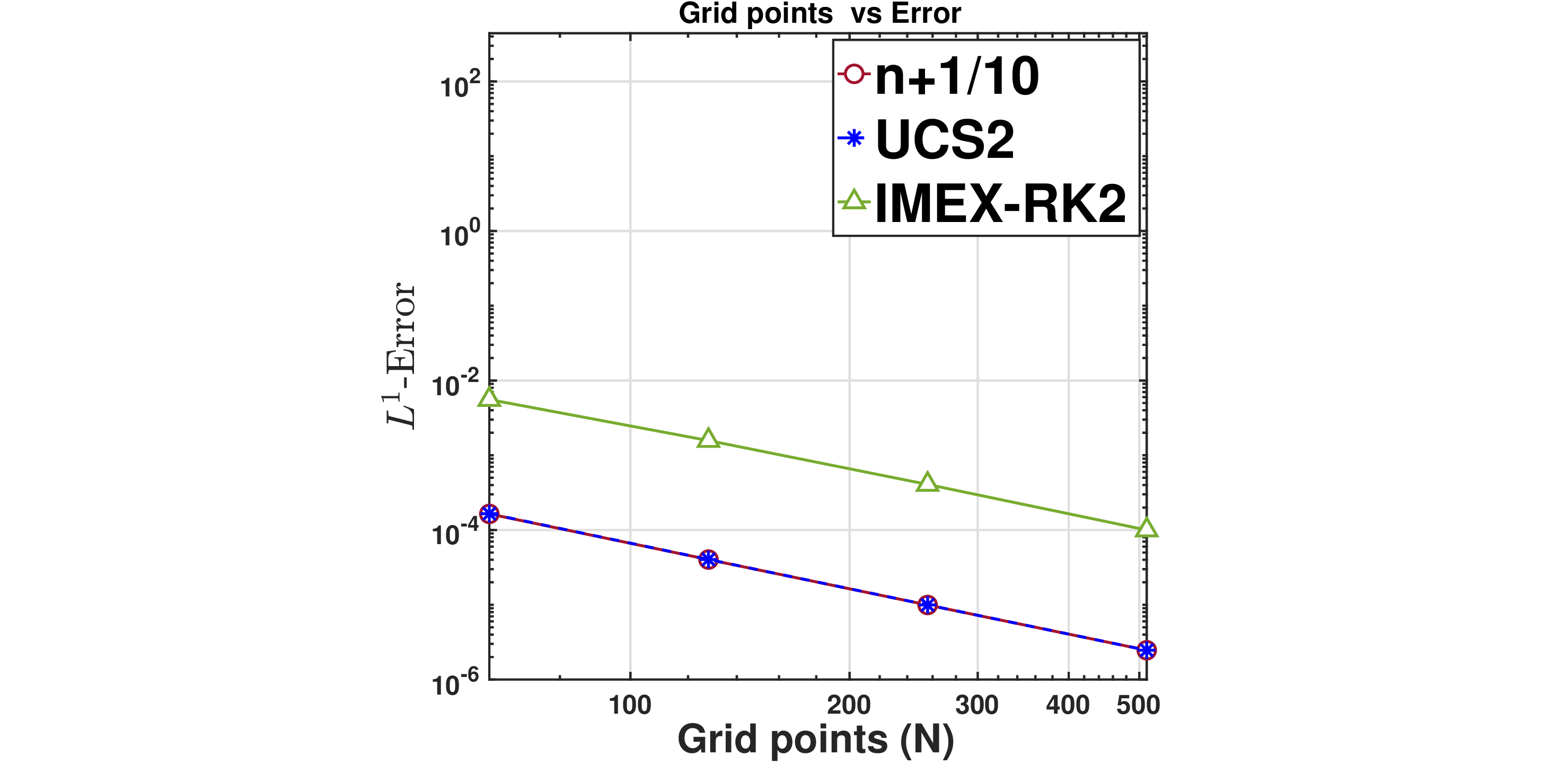}
\caption{}
\label{fig:loglog_1d}
\end{subfigure}
\hfill
\begin{subfigure}[b]{0.49\textwidth}
\centering
\includegraphics[width=\textwidth,height=7.3cm,trim=11.5cm 0cm 10.5cm 0cm,clip]{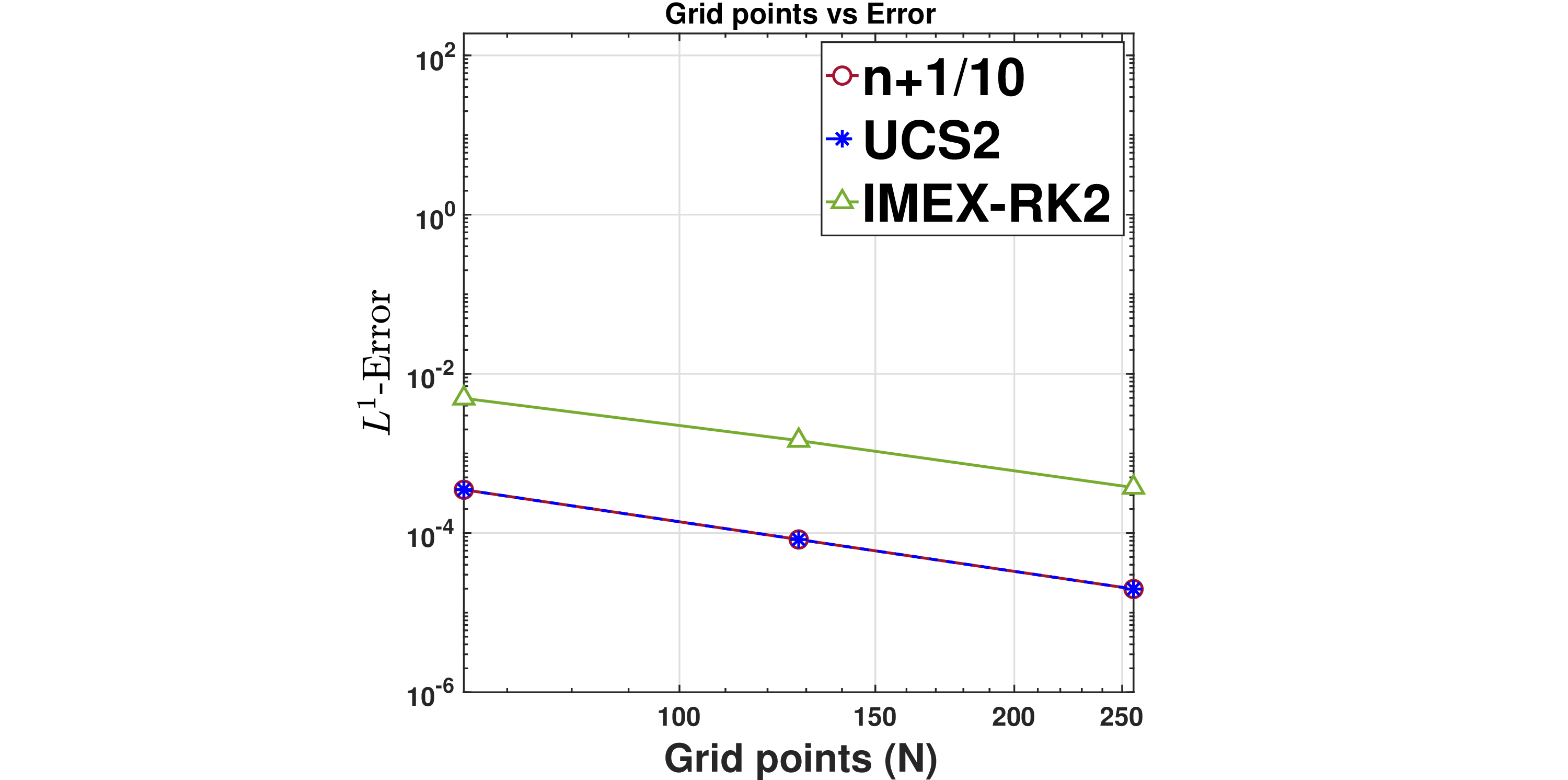}
\caption{}
\label{fig:loglog_2d}
\end{subfigure}
\caption{Log-log plots of the $L^1$-error as a function of the number of grid points for the proposed scheme ($\alpha=1/10$), UCS2, and IMEX-RK2: (a) one-dimensional problem; (b) two-dimensional problem. The results illustrate the second-order convergence of the numerical methods for both the one- and two-dimensional Jin-Xin relaxation systems.}
\label{order_convergence1_2d}
\end{figure}

To examine the accuracy of the scheme in multiple spatial dimensions, we consider the two-dimensional Jin-Xin relaxation system on $[0,1]\times[0,1]$:
\begin{equation}\label{Xin_Jin2d_model}
\begin{cases}
u_t+v_x+w_y=0,\\[2mm]
v_t+u_x=-\dfrac{1}{\varepsilon}(v-au),\\[2mm]
w_t+u_y=-\dfrac{1}{\varepsilon}(w-bu).
\end{cases}
\end{equation}
Here, $u=u(x,y,t)$ is the macroscopic variable, while $v=v(x,y,t)$ and $w=w(x,y,t)$ are auxiliary relaxation variables. The parameters $a$ and $b$ are prescribed constants. The corresponding smooth, well-prepared initial data are
\begin{equation}\label{smoothdata:XinJin2d}
\begin{cases}
u(x,y,0)=\sin(2\pi x)\sin(2\pi y),\\[1mm]
v(x,y,0)=a\,u(x,y,0),\\[1mm]
w(x,y,0)=b\,u(x,y,0),
\end{cases}
\end{equation}
with $a=0.7$ and $b=0.5$.
\paragraph{\textbf{Computational Setup.}}
For both test problems, periodic boundary conditions are imposed in all spatial directions. The relaxation parameter is set to $\varepsilon=10^{-10}$ to examine the behavior of the schemes in the stiff relaxation regime. The final time is $T=0.35$, and the CFL number is fixed at $0.9$. The time-step size is determined according to the CFL condition specified earlier. The computations are performed on sequences of uniformly refined grids, with the finest grids containing $N=320$ points in one dimension and $N_x=N_y=400$ points in two dimensions. The exact solutions are used to evaluate the numerical errors in both cases. Unless otherwise stated, the proposed scheme is considered with $\alpha=1/10$.
\paragraph{\textbf{Convergence Results.}}
For the one-dimensional test problem, Table~\ref{tab:u_v} summarizes the $L^1$- and $L^\infty$-norm errors and their corresponding observed convergence orders for the solution components $u$ and $v$. As the spatial grid is refined, the errors in both components decrease consistently, with the convergence rates remaining close to two. These results demonstrate that the proposed scheme achieves the expected second-order accuracy for the smooth solution in the stiff relaxation regime. To assess the equilibrium-preserving behavior of the generalized scheme, we compute the \(L^1\)-norm of the residual \(\|v_{\mathrm{num}}-au_{\mathrm{num}}\|\) for the one-dimensional linear Jin-Xin model. Table~\ref{tab:equilibrium_residual} shows that the residual decreases with mesh refinement, indicating the equilibrium relation \(v=au\) with the numerical solution.
\begin{table}[ht]
\caption{$L^{1}$- and $L^{\infty}$-norm errors and the corresponding orders of convergence for $u$ and $v$ with smooth initial data for $\alpha = 1/10$.}
\centering
\resizebox{0.79\textwidth}{!}{
\begin{tabular}{c c c c c c c c c}
\toprule
\multirow{2}{*}{$N$}
&
\multicolumn{4}{c}{$u$}
&
\multicolumn{4}{c}{$v$}
\\
\cmidrule(lr){2-5}
\cmidrule(lr){6-9}
&
$L^{1}$-Error & Order
&
$L^{\infty}$-Error & Order
&
$L^{1}$-Error & Order
&
$L^{\infty}$-Error & Order
\\
\midrule
40  & 6.69687e-04 & - & 3.37869e-03 & - & 1.10545e-03 & -   & 5.16353e-03 & - \\
80  & 1.64288e-04 & 2.027 & 9.18995e-04 &  1.878 & 3.58352e-04 & 1.625 & 1.64988e-03  & 1.646\\
160 & 4.02614e-05 & 2.029 & 2.37720e-04 & 1.951 & 8.55469e-05  & 2.067 &  3.90728e-04 & 2.078\\
320 & 9.93902e-06 & 2.018  & 6.04014e-05 & 1.977 & 1.90798e-05 & 2.165 & 8.85562e-05  & 2.142\\
640 & 2.46591e-06 & 2.011 & 1.52282e-05 & 1.988 & 3.58319e-06 & 2.413 & 1.77150e-05 & 2.322\\
\bottomrule
\end{tabular}
}
\label{tab:u_v}
\end{table}
\begin{table}[htbp]
\centering
\caption{$L^1$-norm of the numerical equilibrium residual for $a=0.7$ and $\alpha = 1/10$.}
\label{tab:equilibrium_residual}
\begin{tabular}{cc}
\hline
$N$ & $\|v_{\mathrm{num}}-0.7u_{\mathrm{num}}\|$ \\
\hline
320   & $1.22845e-05$ \\
640   & $1.87866e-06$ \\
1280  & $5.44603e-07$ \\
20000 & $2.23122e-09$ \\
50000 & $5.03975e-10$ \\
\hline
\end{tabular}
\end{table}
For the two-dimensional problem in the stiff relaxation regime, the convergence properties of the proposed scheme are assessed using the $L^1$- and $L^\infty$-norm errors and the corresponding observed convergence rates, as presented in Table~\ref{tab:n110}. The results show a systematic decrease in numerical errors as the mesh is refined, demonstrating the accuracy and convergence of the scheme.
% \renewcommand{\arraystretch}{1.2}
% \begin{table}[H]
% \caption{Numerical $L^{1}$-errors and observed convergence orders for the solution components $u$, $v$, and $w$ using the UCS2 scheme with the specified initial data}
% \centering
% %\resizebox{\textwidth}{!}{
% \begin{tabular}{c c c c c c c}
% \toprule
% \multirow{2}{*}{$N \times N$}
% &
% \multicolumn{2}{c}{$u$}
% &
% \multicolumn{2}{c}{$v$}
% &
% \multicolumn{2}{c}{$w$}
% \\
% \cmidrule(lr){2-3}
% \cmidrule(lr){4-5}
% \cmidrule(lr){6-7}
% &
% $L^{1}$-Error & Order 
% &
% $L^{1}$-Error & Order 
% &
% $L^{1}$-Error & Order 
% \\
% \midrule
% $32\times 32$  & 1.44183e-03 & - & 1.00928e-03 & - & 7.20916e-04 & - \\
% $64\times 64$  & 3.50034e-04  &  2.042  &  2.45024e-04  & 2.042 & 1.75017e-04 & 2.042  \\
% $128\times 128$ & 8.29779e-05 & 2.077  & 5.80846e-05 & 2.077 & 4.14890e-05 & 2.077 \\
% $256 \times 256$ & 1.98244e-05 & 2.066 & 1.38771e-05 & 2.065  & 9.91220e-06 & 2.065 \\
% \bottomrule
% \end{tabular}
% %}
% \label{tab:ucs}
% \end{table}
%\renewcommand{\arraystretch}{1.2}
\begin{table}[ht]
\caption{The $L^{1}$- and $L^\infty$- errors and corresponding convergence rates for $u$, $v$, and $w$ computed using the $n+1/10$ scheme for the given initial data}
\centering
\resizebox{0.99\textwidth}{!}{
\begin{tabular}{c c c c c c c c c c c c c}
\toprule
\multirow{2}{*}{$N \times N$}
&
\multicolumn{4}{c}{$u$}
&
\multicolumn{4}{c}{$v$}
&
\multicolumn{4}{c}{$w$}
\\
\cmidrule(lr){2-5}
\cmidrule(lr){6-9}
\cmidrule(lr){10-13}
&
$L^{1}$-Error & Order
&
$L^{\infty}$-Error & Order
&
$L^{1}$-Error & Order
&
$L^{\infty}$-Error & Order
&
$L^{1}$-Error & Order
&
$L^{\infty}$-Error & Order
\\
\midrule
$32\times 32$  & 1.44183e-03 & - & 6.05587e-03 &  - & 1.00928e-03 & - &  4.23911e-03  &  - & 7.20916e-04 & - & 3.02794e-03 & -\\
$64\times 64$  & 3.50033e-04  &  2.042 & 1.66475e-03   &  1.863 &  2.45023e-04  & 2.042 & 1.16533e-03  &   1.863 & 1.75016e-04 & 2.042 &  8.32377e-04  &  1.863\\
$128\times 128$ & 8.29779e-05 & 2.076 & 3.94323e-04 & 2.078 & 5.80845e-05 & 2.076 &  2.76026e-04 & 2.078  & 4.14889e-05 & 2.076 &  1.97162e-04  & 2.078\\
$256 \times 256$ & 1.98244e-05 & 2.065 & 9.74242e-05 & 2.017 & 1.38771e-05 & 2.065 &  6.81969e-05  & 2.017 & 9.91221e-06 & 2.065 &  4.87121e-05 & 2.017\\
\bottomrule
\end{tabular}
}
\label{tab:n110}
\end{table}
% \renewcommand{\arraystretch}{1.2}
% \begin{table}[H]
% \caption{$L^{1}$-errors and observed orders of convergence for $u$, $v$, and $w$ obtained using the UCS2 scheme for the prescribed initial data}
% \centering
% %\resizebox{\textwidth}{!}{
% \begin{tabular}{c c c c c c c}
% \toprule
% \multirow{2}{*}{$N \times N$}
% &
% \multicolumn{2}{c}{$u$}
% &
% \multicolumn{2}{c}{$v$}
% &
% \multicolumn{2}{c}{$w$}
% \\
% \cmidrule(lr){2-3}
% \cmidrule(lr){4-5}
% \cmidrule(lr){6-7}
% &
% $L^{1}$-Error & Order 
% &
% $L^{1}$-Error & Order 
% &
% $L^{1}$-Error & Order 
% \\
% \midrule
% $32\times 32$  & 1.46952e-02 & - & 1.028662e-02 & - & 7.347592e-03 & - \\
% $64\times 64$  & 4.92543e-03 & 1.577  &   3.44780e-03 & 1.577 & 2.46271e-03 & 1.577  \\
% $128\times 128$ &1.45172e-03 & 1.762  &   1.01620e-03 & 1.762 &  7.25859e-04 & 1.762 \\
% $256 \times 256$ &  3.73320e-04 & 1.959 &  2.61324e-04 & 1.959  &1.86660e-04 & 1.959
%  \\
% \bottomrule
% \end{tabular}
% %}
% \label{tab:imexrk2}
% \end{table}
Figure~\ref{order_convergence1_2d} shows the log-log plots of the $L^1$-error versus the grid resolution $N$ for the one-dimensional (left) and two-dimensional (right) problems, comparing UCS2, IMEX-RK2, and the generalized scheme with $\alpha=1/10$. Although the schemes exhibit different error magnitudes, with the generalized scheme achieving smaller errors than IMEX-RK2, their convergence slopes indicate approximately second-order accuracy. This confirms that the proposed scheme retains its expected convergence order while providing improved numerical accuracy.
Overall, the numerical results confirm the expected second-order accuracy of the proposed scheme for smooth test problems in both one and two spatial dimensions. In the two-dimensional case, its convergence behavior is comparable to that of the UCS2 scheme described in Appendix~\ref{sub:ucs2d} and the IMEX-RK2 scheme presented in Appendix~\ref{sub:imex2d}. Building on these findings, we next investigate the performance of the proposed scheme for a broader range of test problems, focusing on its ability to capture discontinuities, resolve different solution structures, and reproduce the expected behavior across various relaxation regimes. Note that comparisons with the UCS2 \ref{sub:ucs} and IMEX-RK2 \ref{sub:imex} schemes are provided wherever appropriate.
\subsection{One-Dimensional Test Cases}\label{1d:test}
In this subsection, we investigate the performance of the proposed numerical scheme for a range of one-dimensional test problems involving different relaxation models. The test cases include the Jin-Xin relaxation model, the Broadwell model, and the Euler equations with heat transfer. These experiments are designed to assess the ability of the proposed scheme to capture discontinuities and resolve different solution structures under various relaxation regimes. The numerical results obtained using the proposed generalized scheme for various values of $\alpha$ are compared with those produced by UCS2. 
%Whenever an exact solution is available, it is used as the reference; otherwise, the IMEX-RK2 solution computed on a finer grid is used as a reference for comparison.
\subsubsection{Jin-Xin Relaxation Model}
\label{subsec:jin_xin}
We consider the one-dimensional Jin-Xin relaxation model \eqref{Xin_Jin_model} with discontinuous initial data to examine the performance of the proposed scheme in different relaxation regimes. Particular attention is given to the resolution of sharp transitions and the influence of the parameter $\alpha$ on the numerical solution. The results are assessed using the reference solutions and comparison methodology described above.

\begin{figure}[htb]
    \centering
    \begin{subfigure}[t]{0.48\textwidth}
        \centering
        \includegraphics[width=\linewidth,height=5cm,
        keepaspectratio]{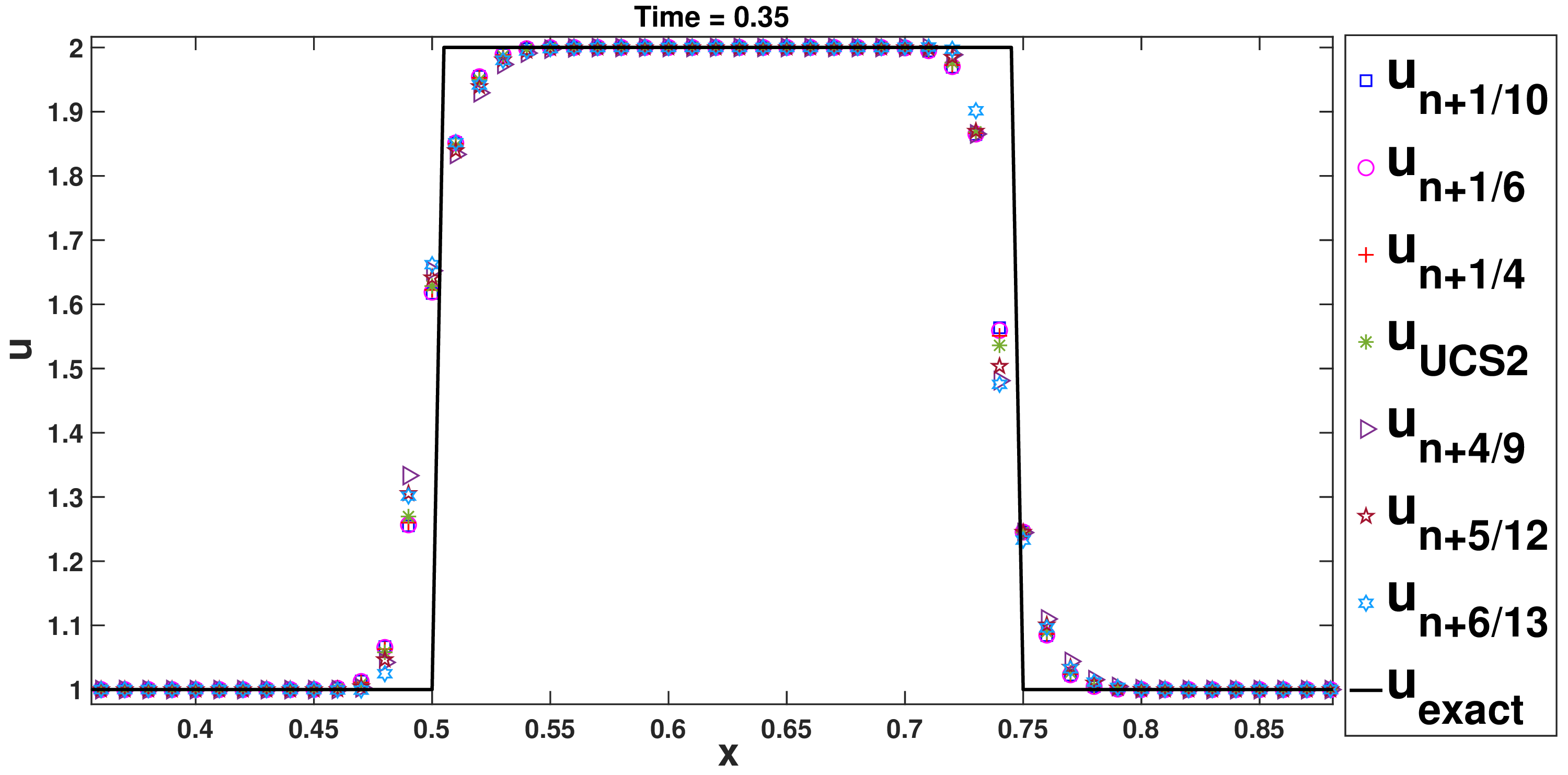}
        \caption{}
    \end{subfigure}
    \hfill
    \begin{subfigure}[t]{0.48\textwidth}
        \centering
        \includegraphics[width=\linewidth,height=4.8cm,
        keepaspectratio]{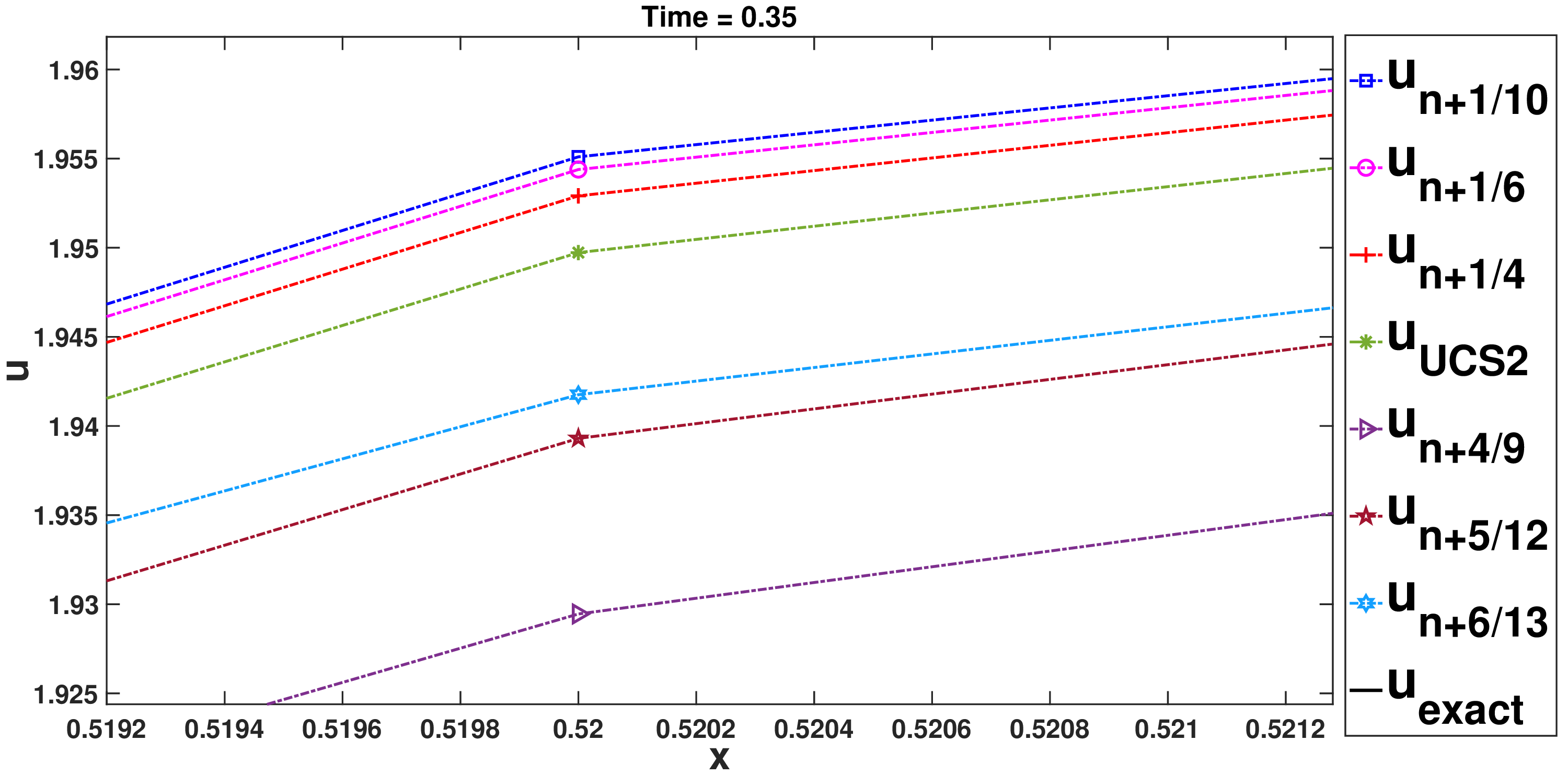}
        \caption{}
    \end{subfigure}
    \vspace{0.3cm}
    \begin{subfigure}[t]{0.48\textwidth}
        \centering
        \includegraphics[width=\linewidth,height=5cm,
        keepaspectratio]{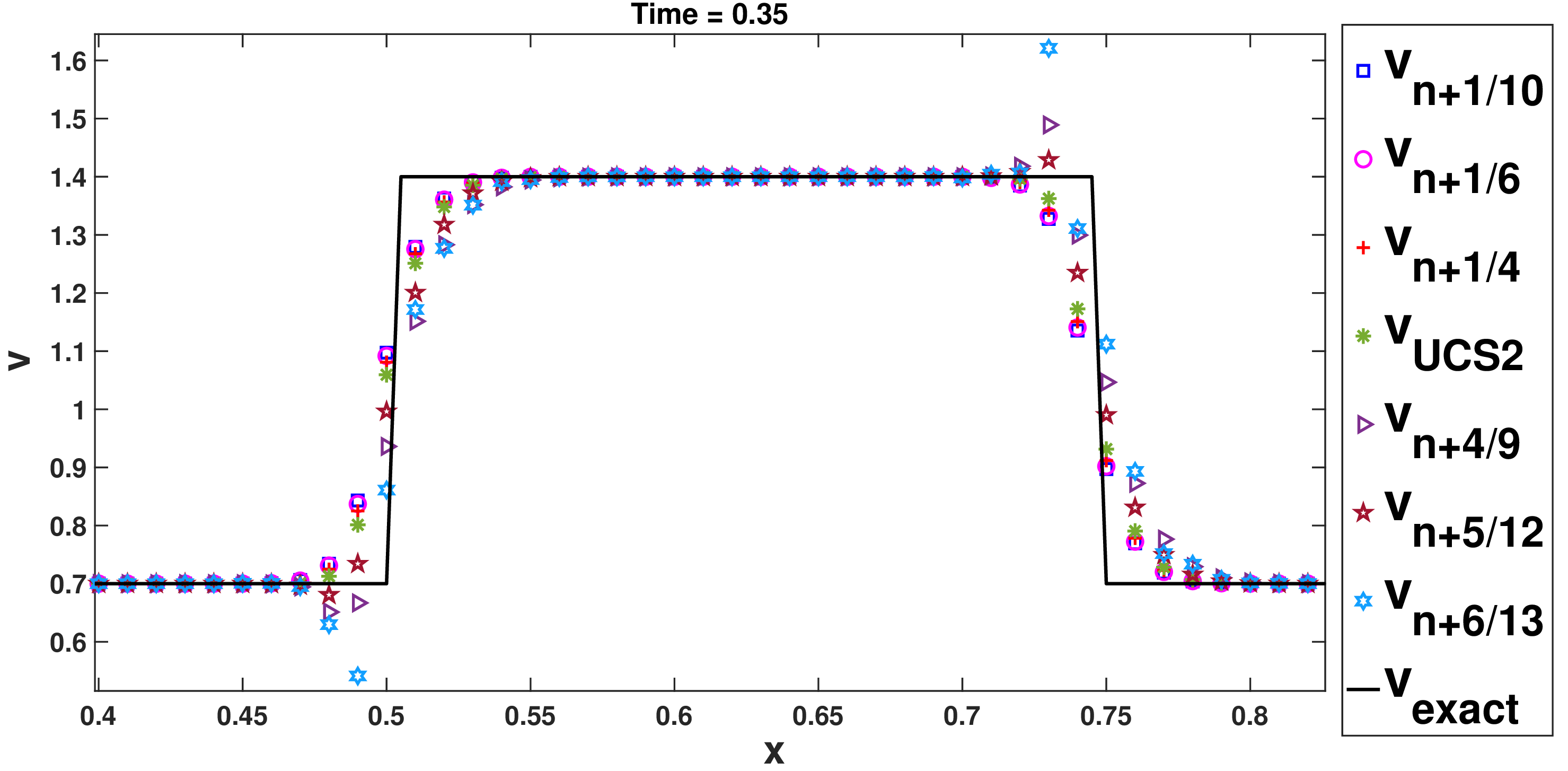}
        \caption{}
    \end{subfigure}
    \hfill
    \begin{subfigure}[t]{0.48\textwidth}
        \centering
        \includegraphics[width=\linewidth,height=5cm,
        keepaspectratio]{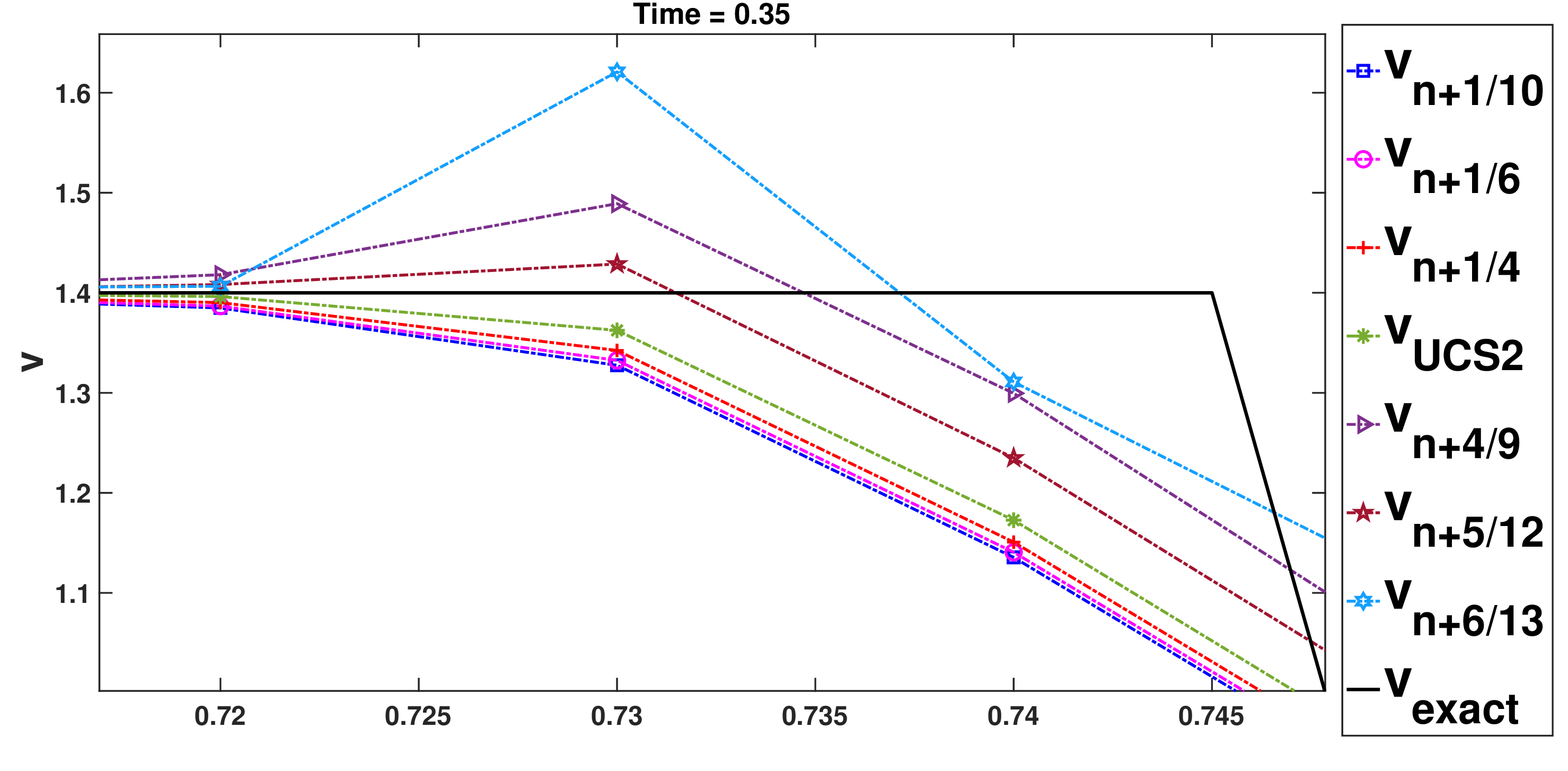}
        \caption{}
    \end{subfigure}

    \caption{Nonsmooth data: Jin-Xin model, $\varepsilon=10^{-10}$:
    (a) solution $u$; (b) zoomed solution $u$;
    (c) solution $v$; (d) zoomed solution $v$.}
    \label{nonsmooth_xinjin}
\end{figure}

\medskip
\noindent\textbf{Test 1: Discontinuous Initial Data.}
We begin with a discontinuous initial profile to examine the ability of the proposed scheme to capture sharp transitions and control spurious oscillations. The initial conditions are prescribed as
\begin{equation}\label{xinjin}
\begin{aligned}
u(x,0) &=
\begin{cases}
2, & 0.25 < x < 0.5,\\
1, & \text{otherwise},
\end{cases}\\
v(x,0) &= a\,u(x,0),
\end{aligned}
\end{equation}
where $a=0.7$. The computational domain is $[0,1]$, discretized using $N=320$ uniformly distributed grid points, with transmissive boundary conditions imposed at both ends. The relaxation parameter is set to $\varepsilon=10^{-10}$ to investigate the behavior of the numerical schemes in the stiff relaxation regime. To examine the influence of the parameter $\alpha$ on the numerical solutions, we consider
\begin{equation}\label{alpha:para}
\alpha\in
\left\{
\frac{1}{10},\frac{1}{6},\frac{1}{4},
\frac{4}{9},\frac{5}{12},\frac{6}{13}
\right\}.
\end{equation}
For each value of $\alpha$, the numerical solution produced by the proposed scheme is compared with that obtained using UCS2, with the exact solution serving as a reference. The results are presented in Figure~\ref{nonsmooth_xinjin}. The left panels display the solution profiles over the full computational domain, whereas the right panels provide enlarged views near the discontinuity to facilitate a closer examination of the local solution behavior.
The results demonstrate that the choice of $\alpha$ influences the numerical solution in the stiff relaxation regime. In particular, localized overshoots in the auxiliary variable $v$ are observed near the discontinuity for $\alpha=4/9$, $5/12$, and $6/13$, indicating the presence of non-monotone behavior for these parameter values. These observations highlight the importance of selecting $\alpha$ appropriately and examining both solution components when assessing the performance of the proposed scheme for nonsmooth initial data. This experiment complements the smooth-solution convergence study by illustrating how the scheme behaves near discontinuities, where localized numerical oscillations may arise.
\subsubsection{Broadwell Model}
\label{subsec:broadwell}
We next consider the one-dimensional Broadwell relaxation system
\cite{broadwell1964shock,PareschiRusso}, which serves as a benchmark
for assessing the performance of the proposed numerical scheme for
hyperbolic systems with nonlinear relaxation effects. Let
\[
Q=(\rho,m,z)^T
\]
denote the vector of state variables. The governing system is written
as
\begin{equation}\label{Broadwell_model}
Q_t+\mathrm{F}(Q)_x=\mathrm{S}(Q,\varepsilon),
\end{equation}
where the flux and source terms are given by
\[
\mathrm{F}(Q)=
\begin{pmatrix}
m\\
z\\
m
\end{pmatrix},
\qquad
\mathrm{S}(Q,\varepsilon)=
\frac{1}{2\varepsilon}
\begin{pmatrix}
0\\
0\\
\rho^2+m^2-2\rho z
\end{pmatrix}.
\]
Here, $\rho$, $m$, and $z$ denote the density, momentum, and auxiliary
variable, respectively, while $\varepsilon>0$ controls the relaxation
rate.
\begin{figure}[htb]
    \centering
    % First row
    \begin{subfigure}[t]{0.48\textwidth}
        \centering      \includegraphics[width=\linewidth,height=5cm,keepaspectratio]{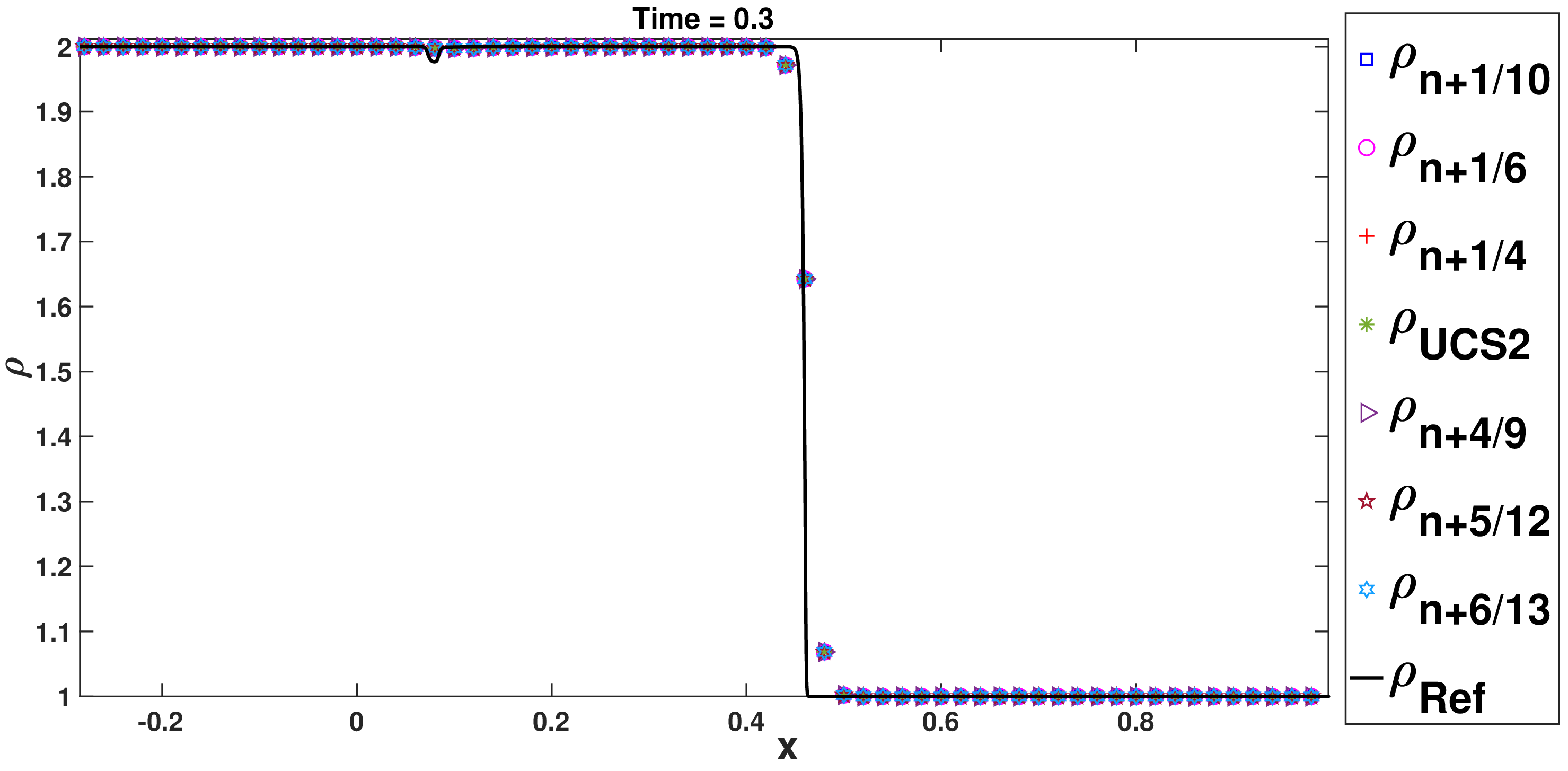}
        \caption{}
        \label{fig:broadwell-rho}
    \end{subfigure}
    \hfill
    \begin{subfigure}[t]{0.48\textwidth}
        \centering       \includegraphics[width=\linewidth,height=5cm,keepaspectratio]{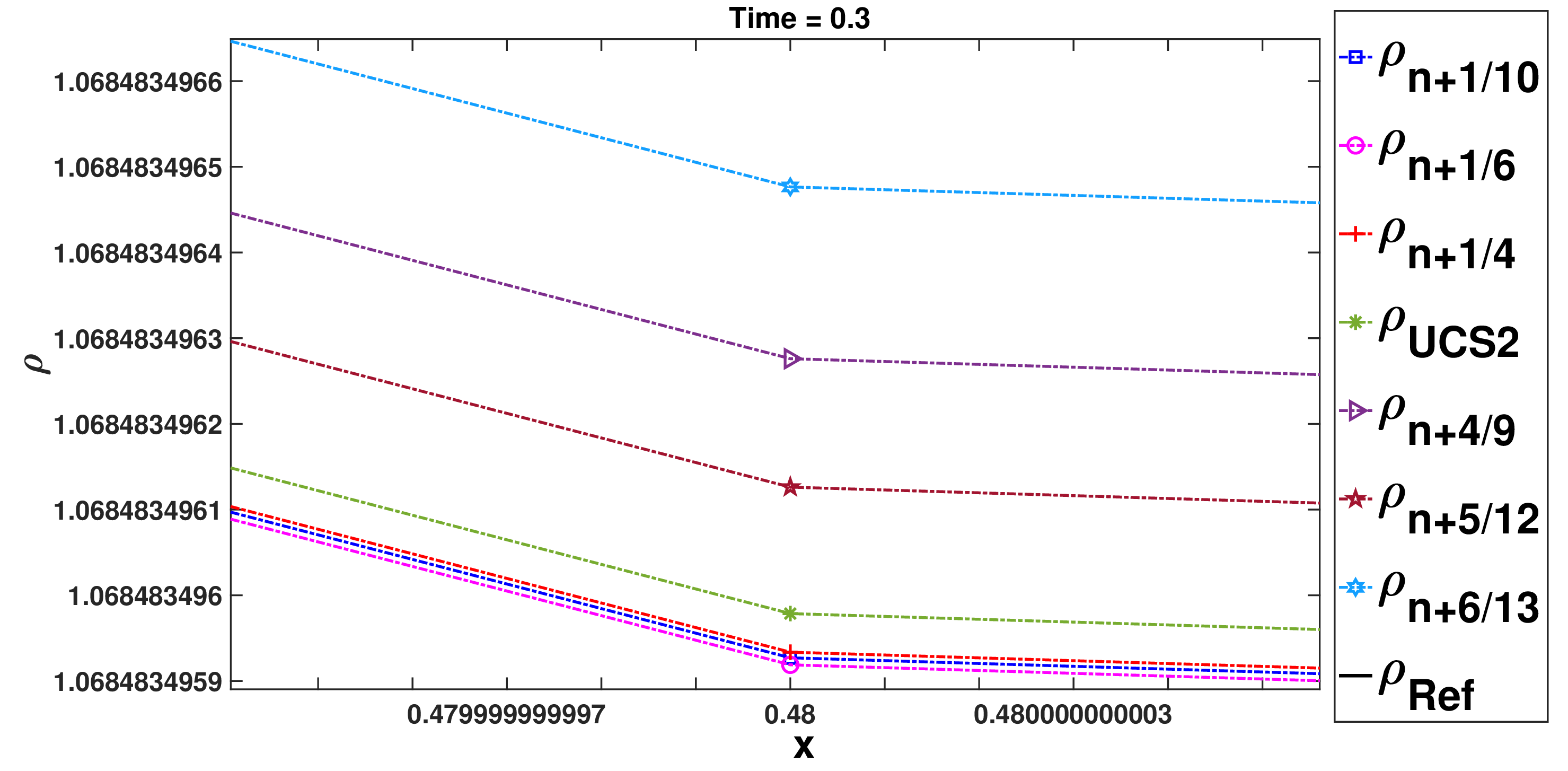}
        \caption{}
        \label{fig:broadwell-rho-zoom}
    \end{subfigure}
    \vspace{0.4cm}
    % Second row
    \begin{subfigure}[t]{0.48\textwidth}
        \centering       \includegraphics[width=\linewidth,height=5cm,keepaspectratio]{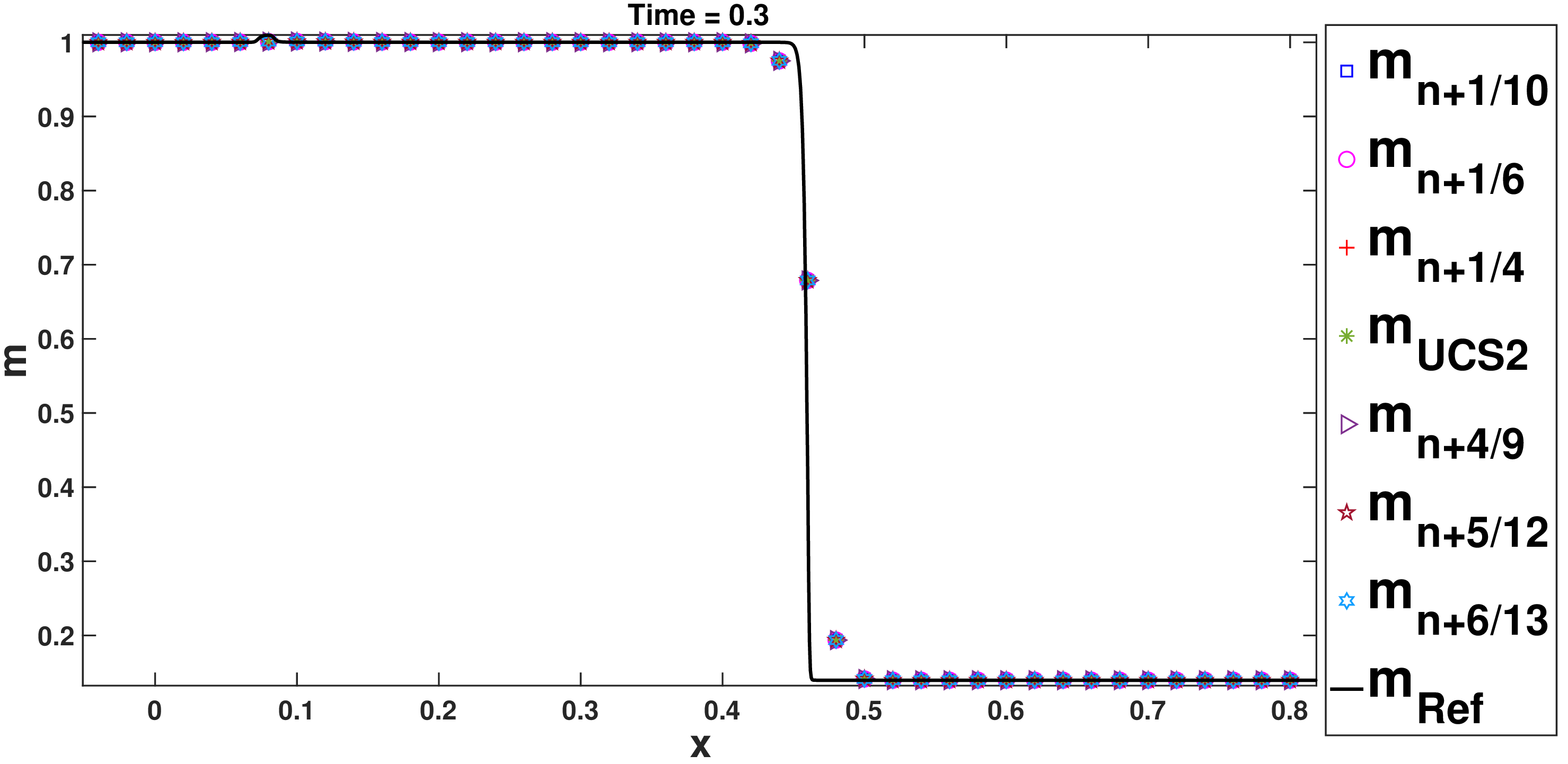}
        \caption{}
        \label{fig:broadwell-m}
    \end{subfigure}
    \hfill
    \begin{subfigure}[t]{0.48\textwidth}
        \centering       \includegraphics[width=\linewidth,height=5cm,keepaspectratio]{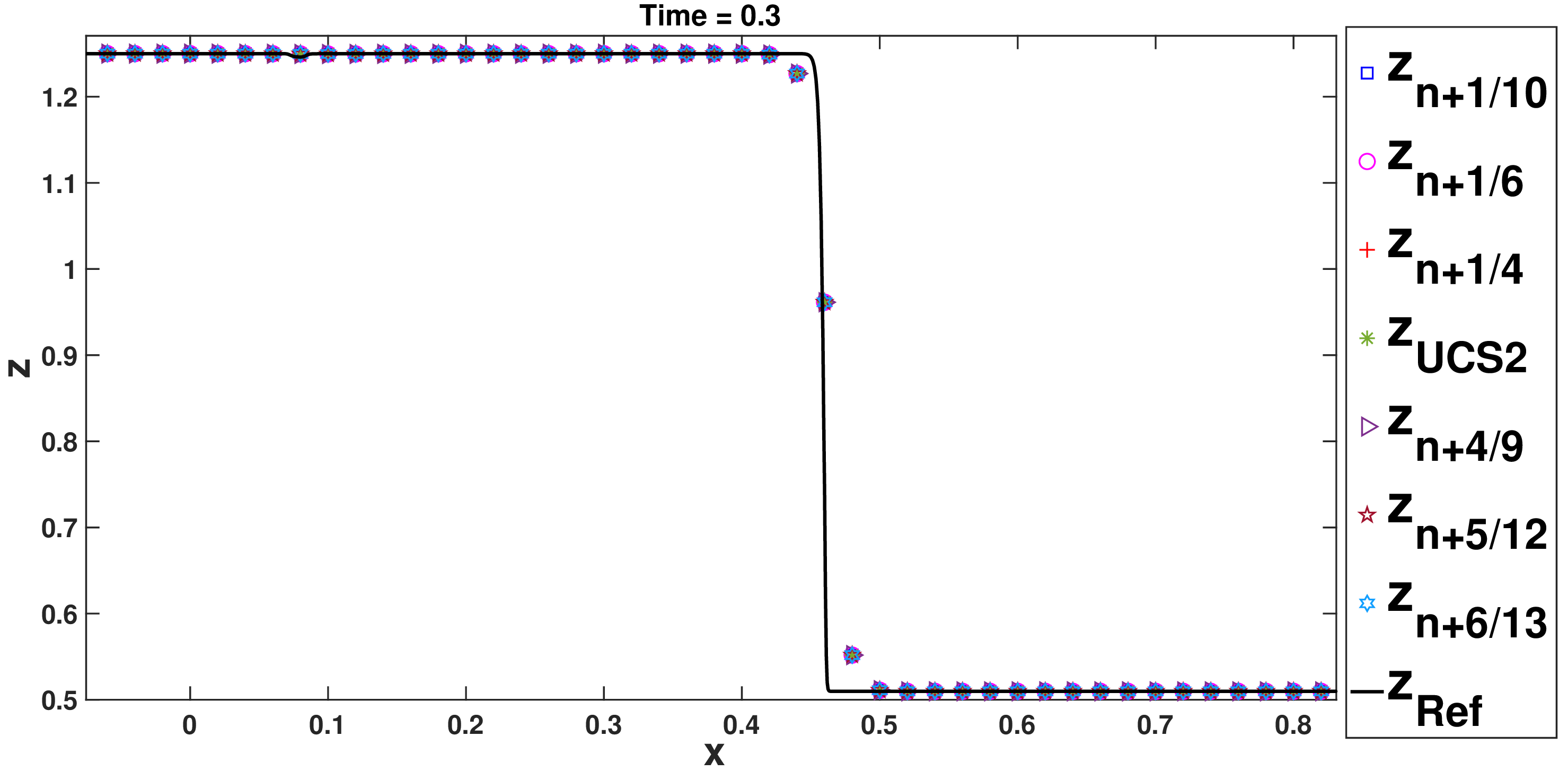}
        \caption{}
        \label{fig:broadwell-z}
    \end{subfigure}
    \caption{Non-smooth initial data for the Broadwell model with $\varepsilon = 10^{-10}$: (a) density $\rho$; (b) magnified view of the density profile; (c) momentum $m$; and (d) variable $z$.}
    \label{nonsmooth:broadwell:ep10}
\end{figure}

\medskip
\noindent\textbf{Test 1: Discontinuous Initial Data.}
We first consider piecewise-constant initial data to examine the ability of the numerical schemes to resolve discontinuities in the presence of nonlinear relaxation effects. The initial conditions are prescribed as
\begin{equation}\label{broadwell:discontinuous}
(\rho(x,0),m(x,0),z(x,0))=
\begin{cases}
(2,1,1), & x\leq 0.2,\\
(1,0.13962,1), & x>0.2.
\end{cases}
\end{equation}
These data introduce a jump discontinuity at $x=0.2$. The computational
domain is $[-1,1]$, with transmissive boundary conditions imposed at
both endpoints, and the simulations are performed up to the final time
$T=0.3$.

To examine the effect of the parameter $\alpha$, we consider the set of values specified in \eqref{alpha:para}.
The proposed scheme is compared with UCS2, while a fine-grid IMEX-RK2 solution is used as a reference. Three relaxation parameters, $\varepsilon=10^{-10}$, $10^{-3}$, and $1$, are considered to examine the numerical behavior in the fluid dynamic limit, thin regime, rarefied regimes, respectively.
\begin{figure}[htb]
    \centering
    % First row: Density
    \begin{subfigure}[t]{0.48\textwidth}
        \centering       \includegraphics[width=\linewidth,height=4.5cm,keepaspectratio]{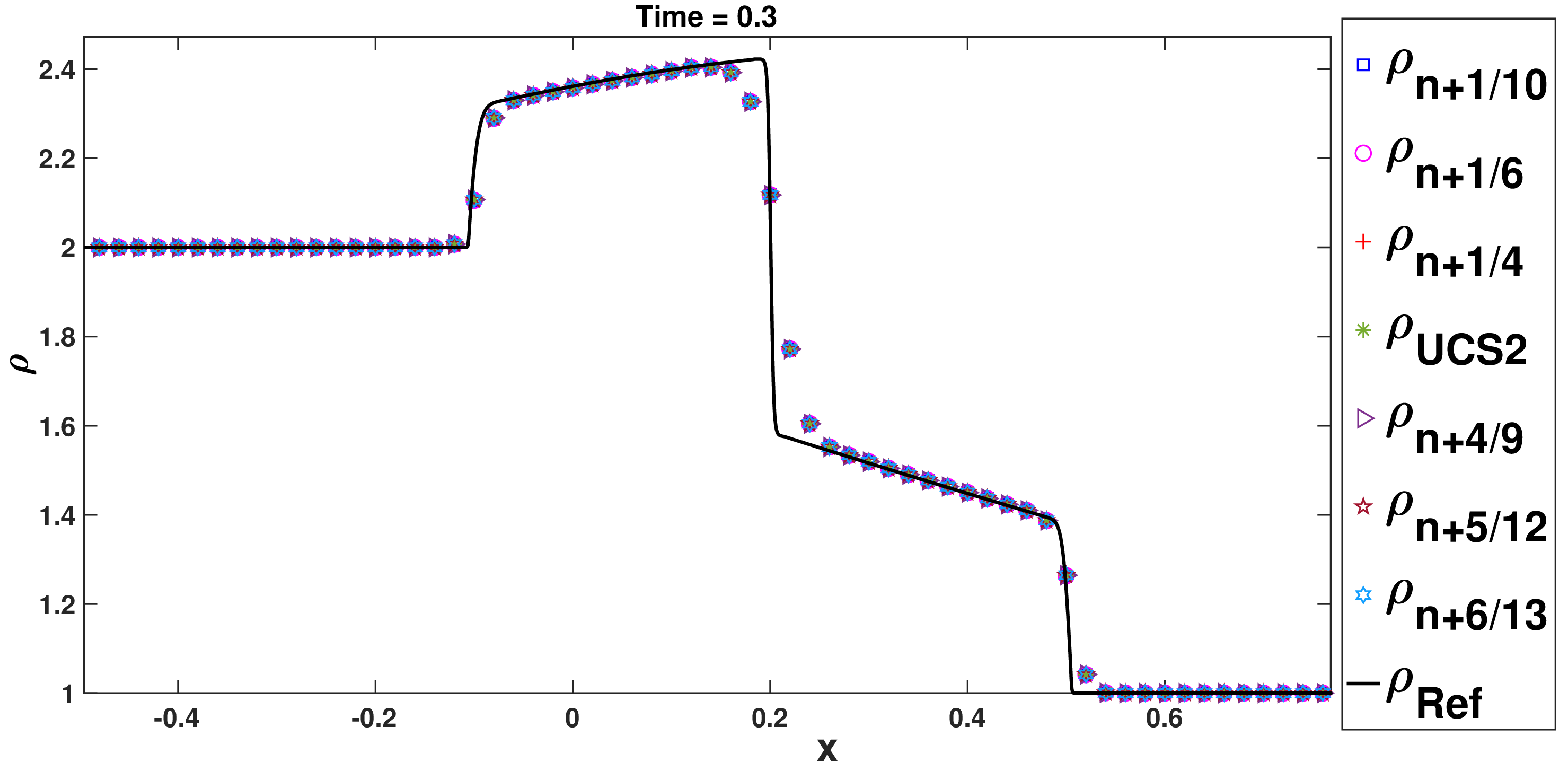}
        \caption{}
        \label{fig:broadwell-ep0-rho}
    \end{subfigure}
    \hfill
    \begin{subfigure}[t]{0.48\textwidth}
        \centering        \includegraphics[width=\linewidth,height=4.5cm,keepaspectratio]{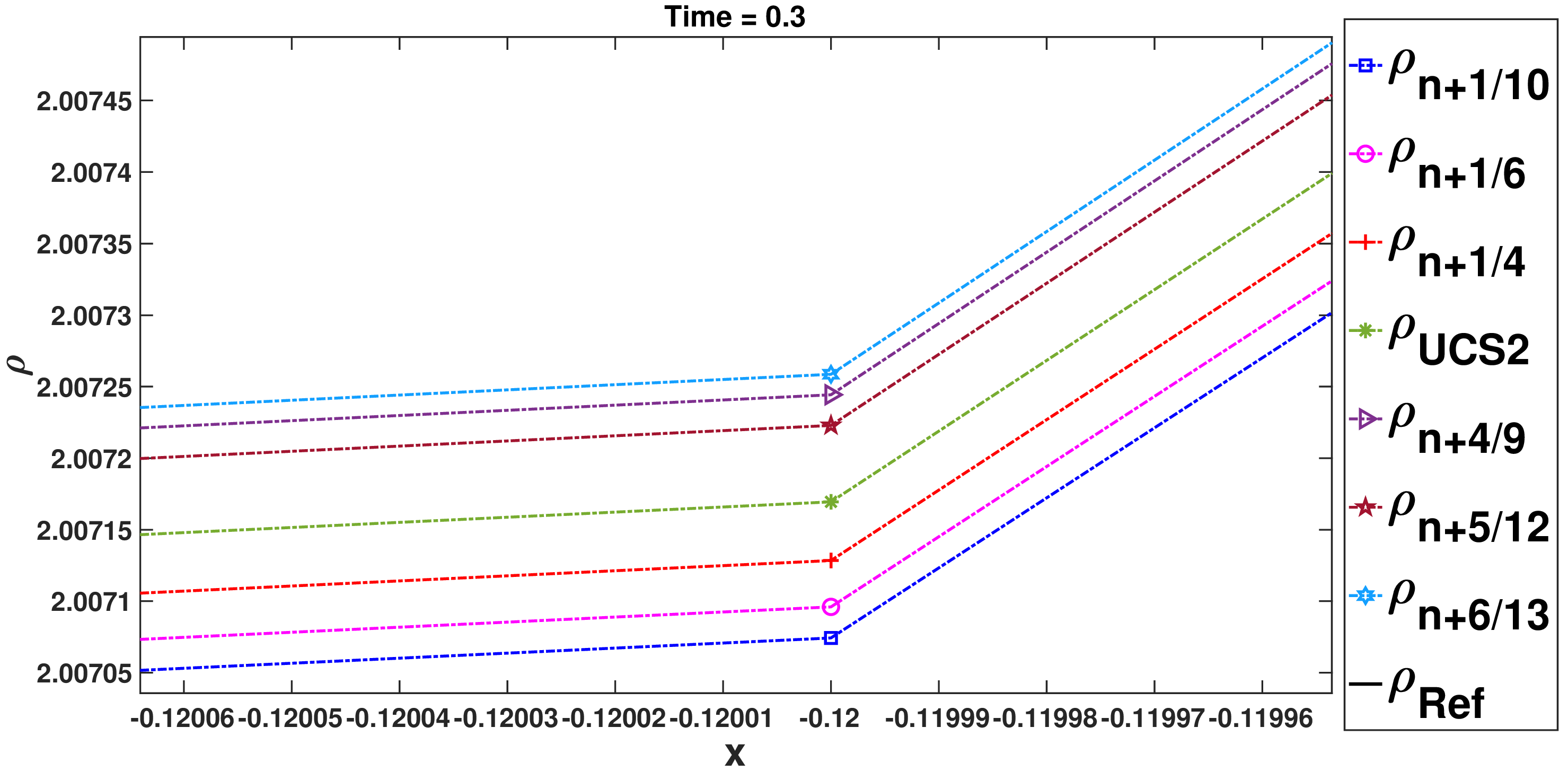}
        \caption{}
        \label{fig:broadwell-ep0-rho-zoom}
    \end{subfigure}

    \vspace{0.4cm}

    % Second row: Momentum
    \begin{subfigure}[t]{0.48\textwidth}
        \centering        \includegraphics[width=\linewidth,height=4.5cm,keepaspectratio]{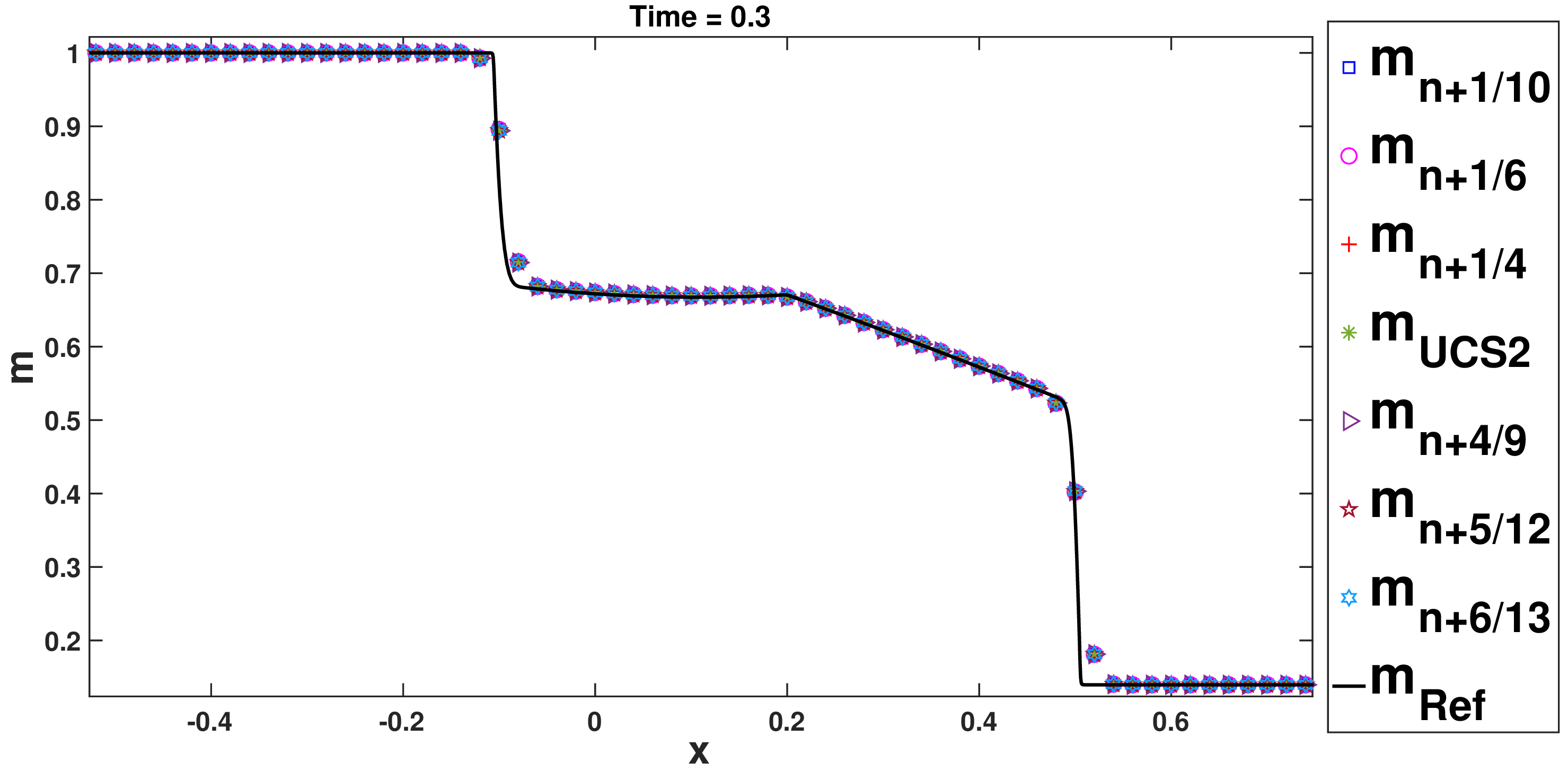}
        \caption{}
        \label{fig:broadwell-ep0-m}
    \end{subfigure}
    \hfill
    \begin{subfigure}[t]{0.48\textwidth}
        \centering        \includegraphics[width=\linewidth,height=4.5cm,keepaspectratio]{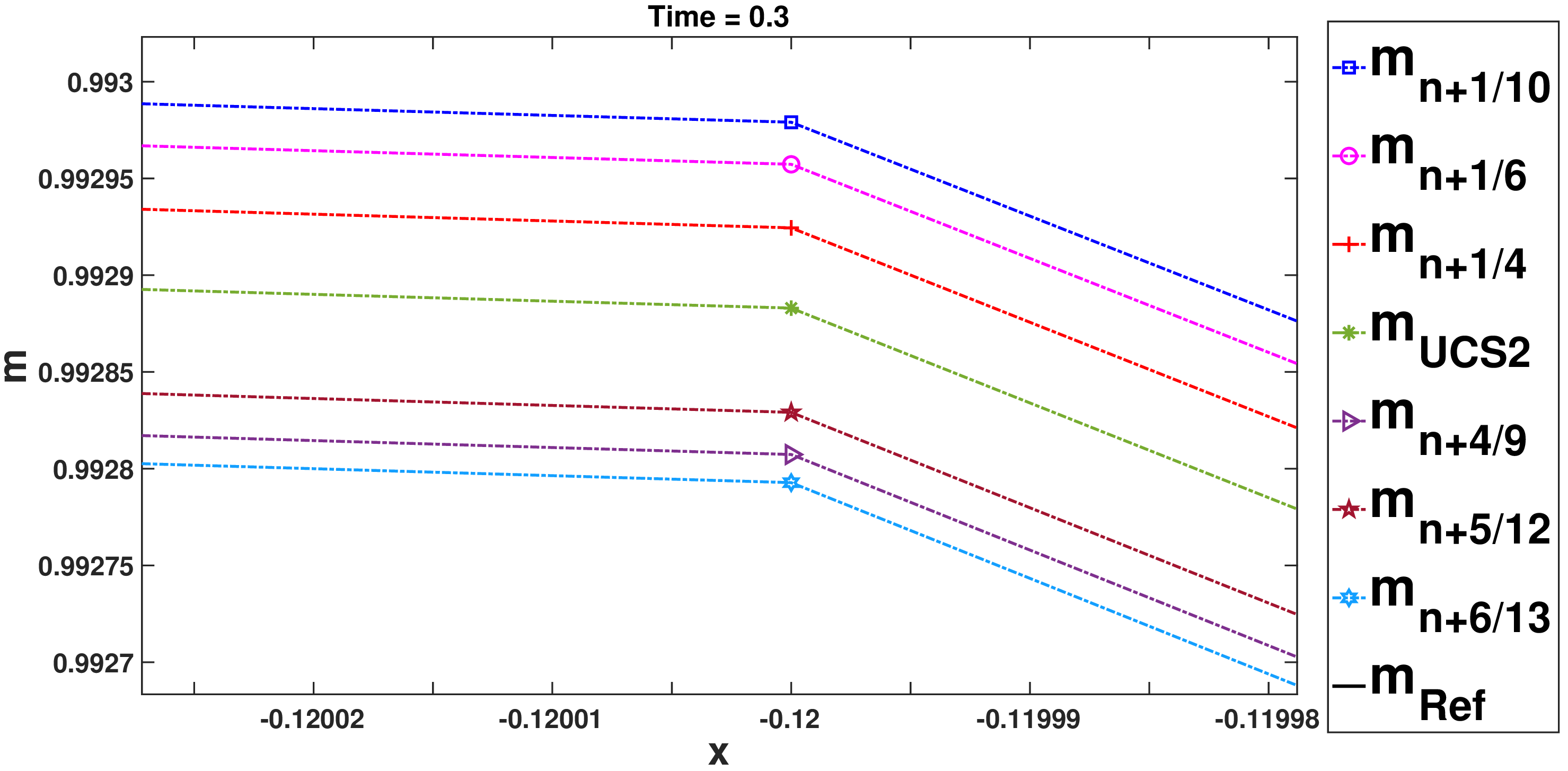}
        \caption{}
        \label{fig:broadwell-ep0-m-zoom}
    \end{subfigure}

    \vspace{0.4cm}

    % Third row: Variable z
    \begin{subfigure}[t]{0.48\textwidth}
        \centering        \includegraphics[width=\linewidth,height=4.5cm,keepaspectratio]{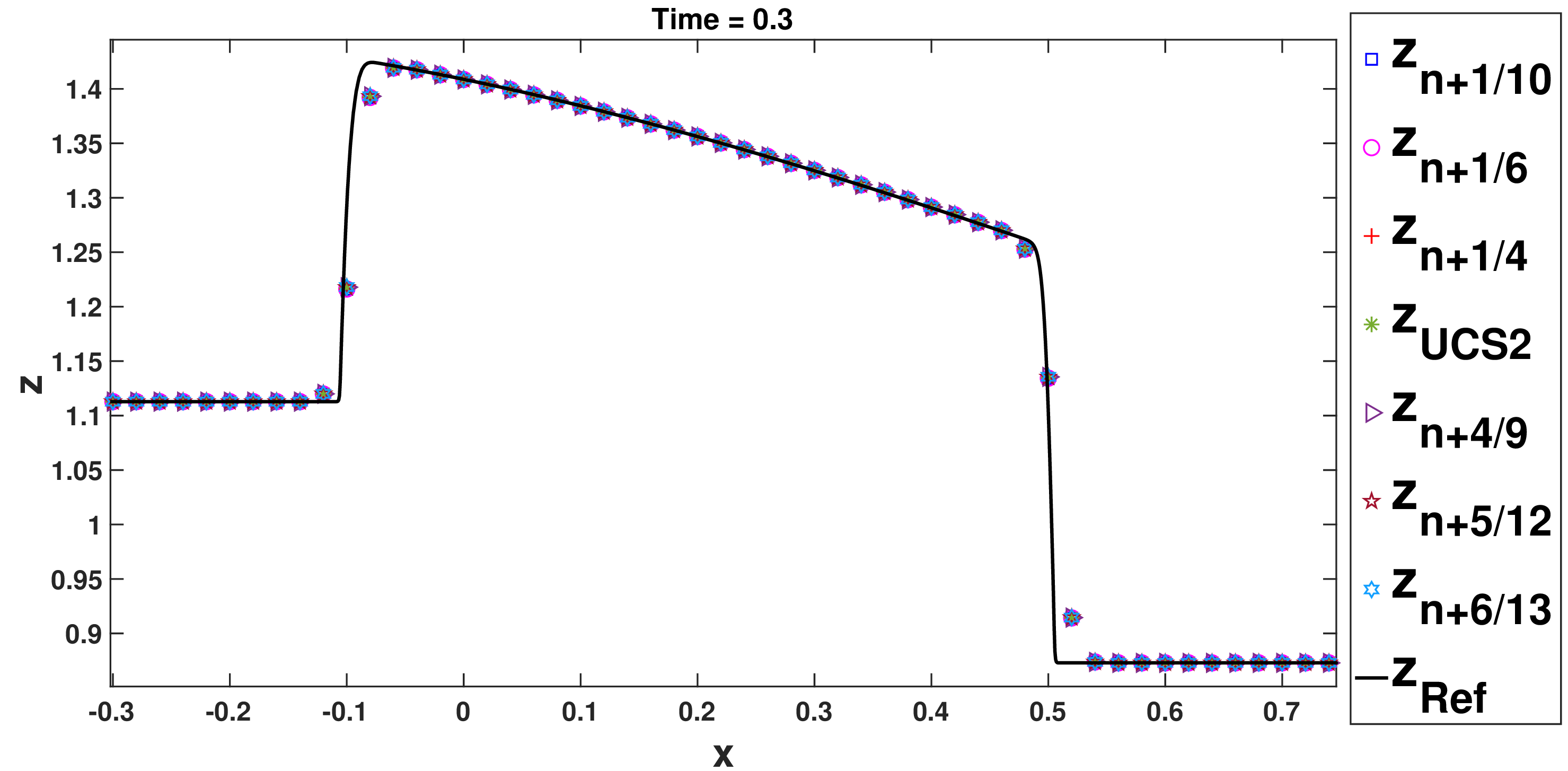}
        \caption{}
        \label{fig:broadwell-ep0-z}
    \end{subfigure}
    \hfill
    \begin{subfigure}[t]{0.48\textwidth}
        \centering        \includegraphics[width=\linewidth,height=4.5cm,keepaspectratio]{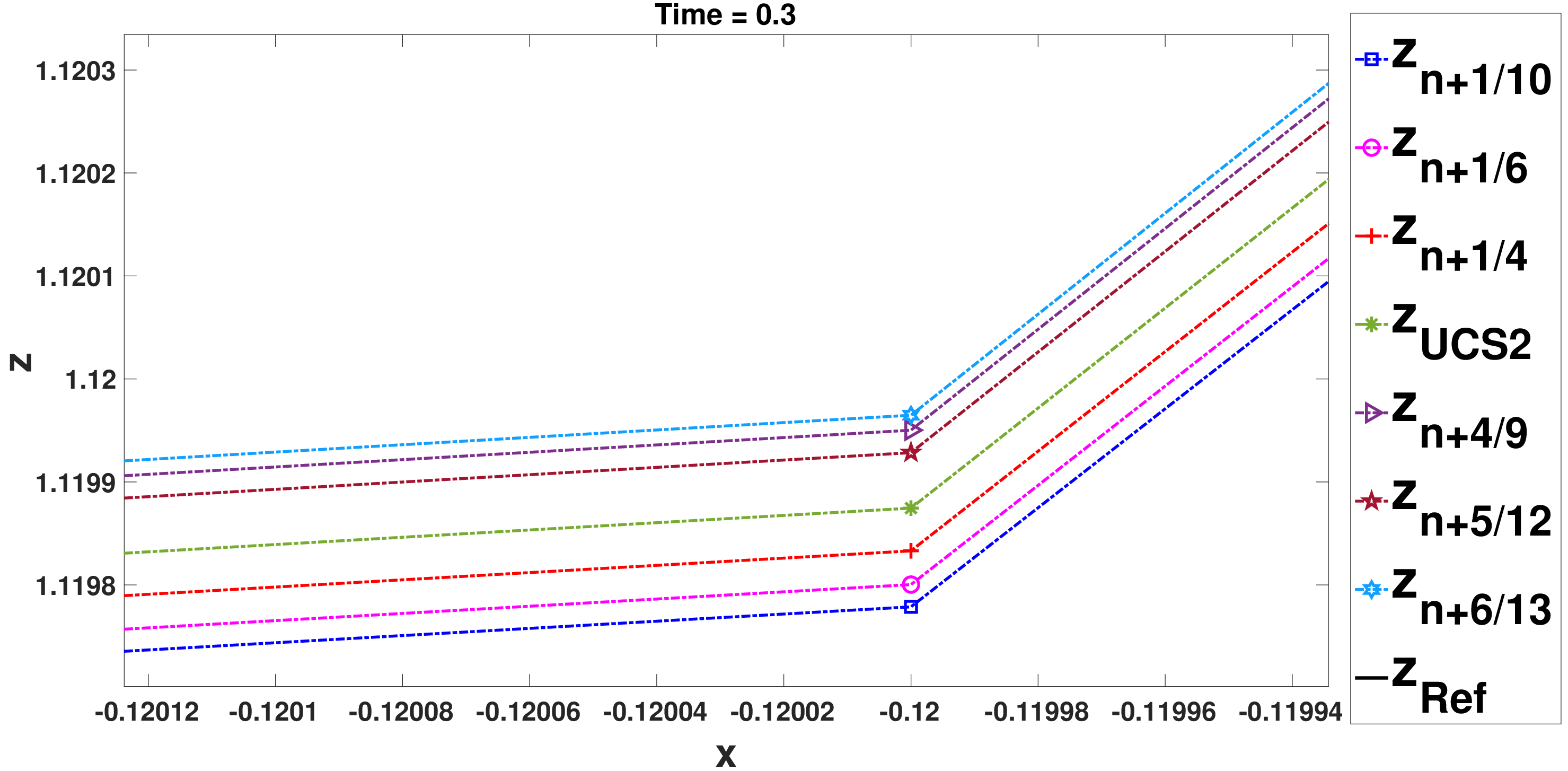}
        \caption{}
        \label{fig:broadwell-ep0-z-zoom}
    \end{subfigure}
    \caption{Numerical results for nonsmooth initial data in the Broadwell model with $\varepsilon = 1$: (a) density $\rho$; (b) magnified view of $\rho$; (c) momentum $m$; (d) magnified view of $m$; (e) variable $z$; and (f) magnified view of $z$.}
    \label{nonsmooth:broadwell:ep0}
\end{figure}
\begin{figure}[htb]
    \centering
    % First row
    \begin{subfigure}[t]{0.48\textwidth}
        \centering        \includegraphics[width=\linewidth,height=5cm,keepaspectratio]{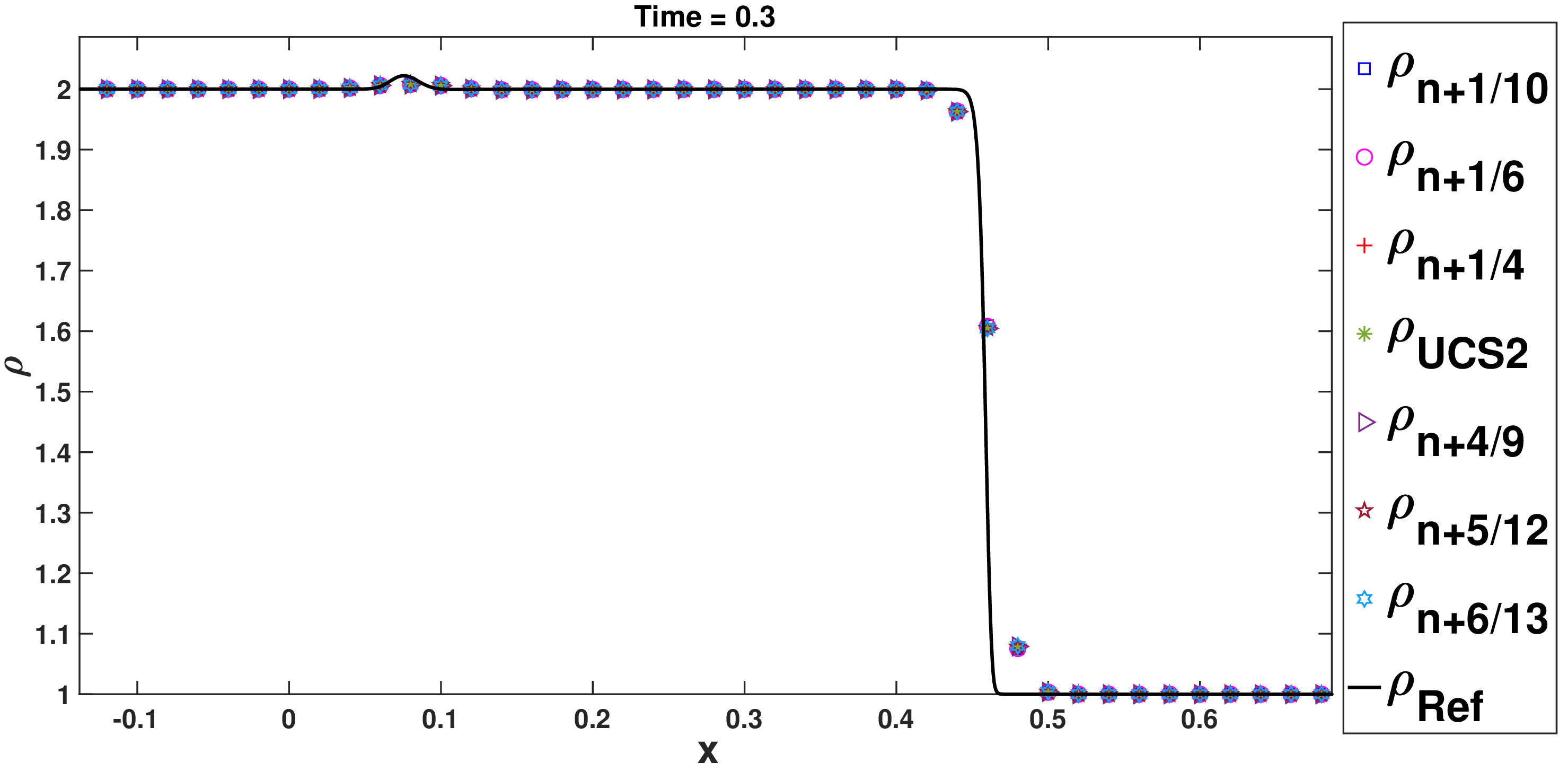}
        \caption{}
        \label{fig:broadwell-ep3-rho}
    \end{subfigure}
    \hfill
    \begin{subfigure}[t]{0.48\textwidth}
        \centering        \includegraphics[width=\linewidth,height=5cm,keepaspectratio]{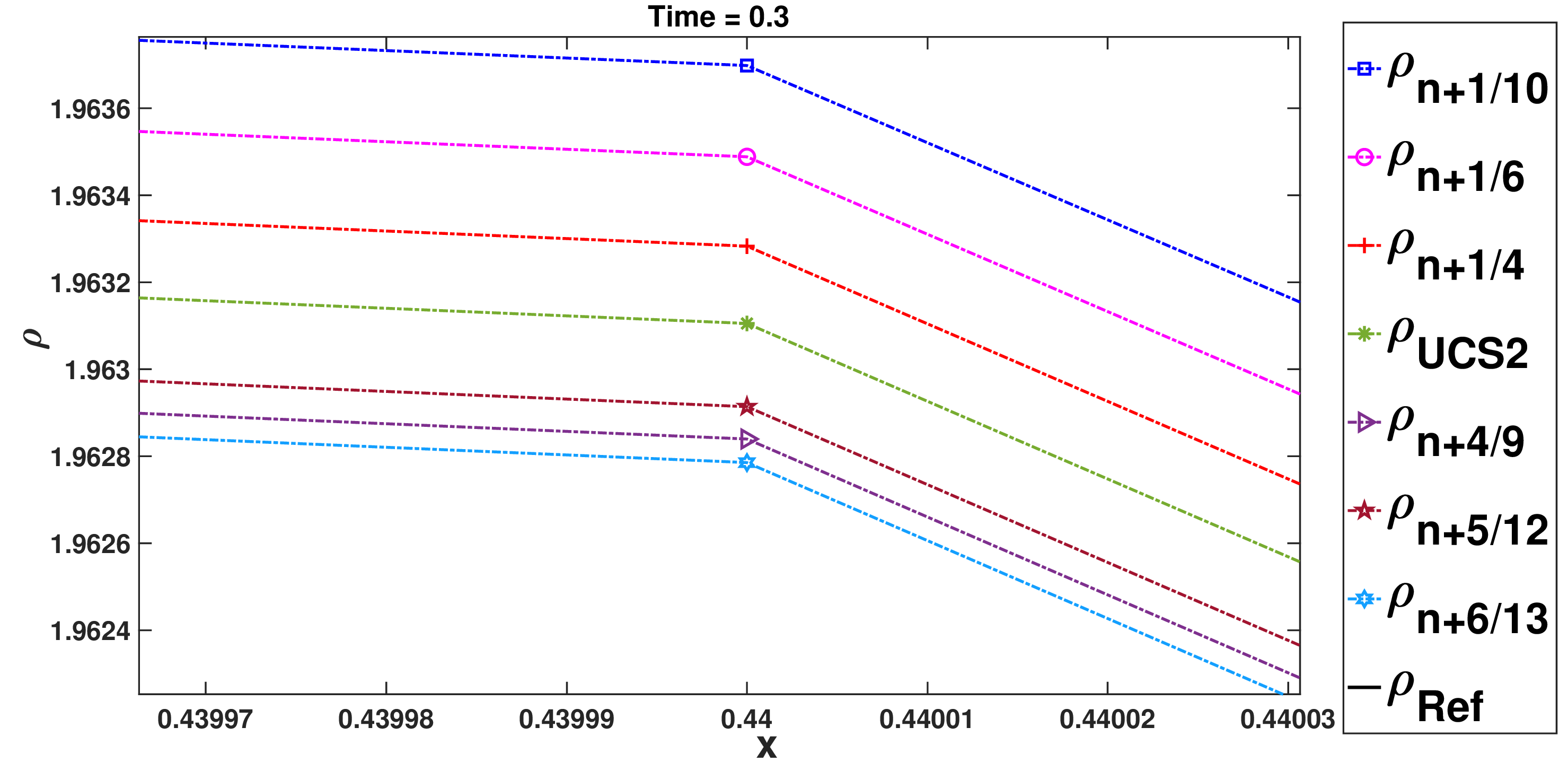}
        \caption{}
        \label{fig:broadwell-ep3-rho-zoom}
    \end{subfigure}
    \vspace{0.4cm}
    % Second row
    \begin{subfigure}[t]{0.48\textwidth}
        \centering        \includegraphics[width=\linewidth,height=5cm,keepaspectratio]{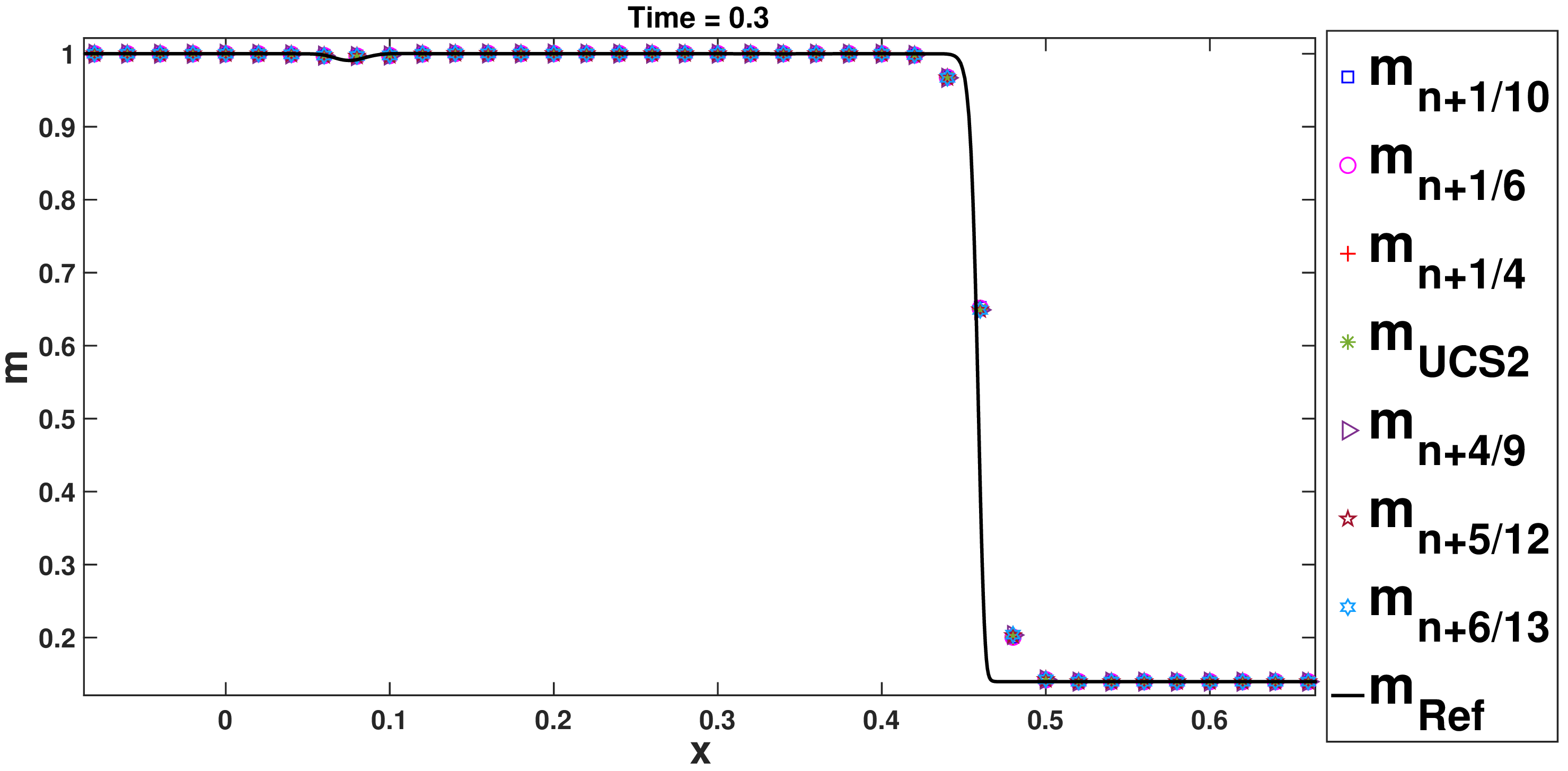}
        \caption{}
        \label{fig:broadwell-ep3-m}
    \end{subfigure}
    \hfill
    \begin{subfigure}[t]{0.48\textwidth}
        \centering        \includegraphics[width=\linewidth,height=5cm,keepaspectratio]{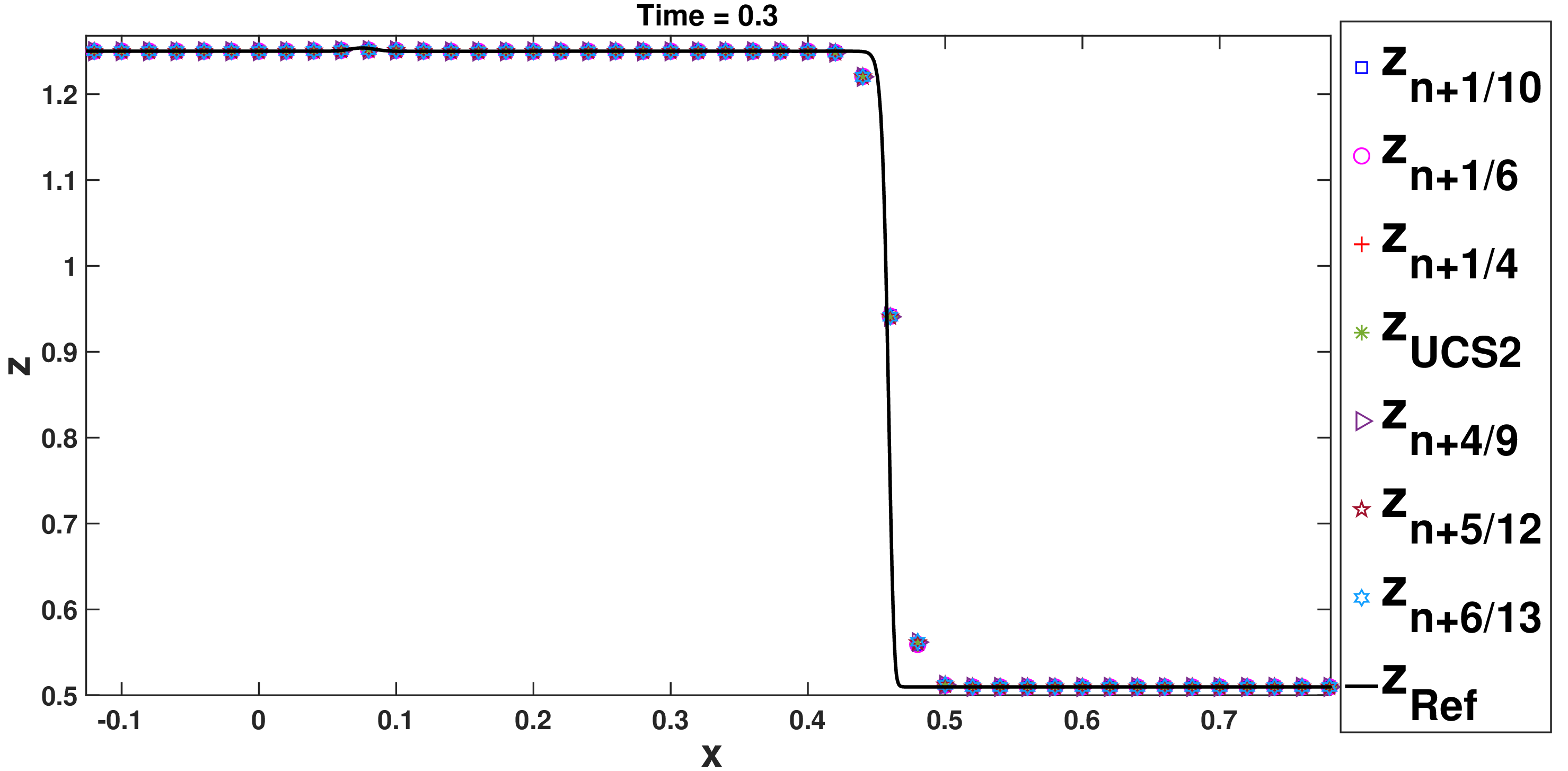}
        \caption{}
        \label{fig:broadwell-ep3-z}
    \end{subfigure}
    \caption{Numerical results for nonsmooth initial data in the Broadwell model with $\varepsilon = 10^{-3}$: (a) density $\rho$; (b) magnified view of $\rho$; (c) momentum $m$; and (d) variable $z$.}
    \label{nonsmooth:broadwell:ep3}
\end{figure}
The results for the stiff relaxation regime, $\varepsilon=10^{-10}$, are presented in Figure~\ref{nonsmooth:broadwell:ep10}. The solution profiles of $\rho$, $m$, and $z$ are shown along with enlarged views near the discontinuity to facilitate a detailed comparison of the local resolution of $\rho$ for different values of $\alpha$. For the weak-relaxation regime, $\varepsilon=1$, Figure~\ref{nonsmooth:broadwell:ep0} displays the full solution profiles of all three variables in the left panels and the corresponding zoomed-in views for different values of $\alpha$ in the right panels. The results for the intermediate relaxation regime, $\varepsilon=10^{-3}$, are presented in Figure~\ref{nonsmooth:broadwell:ep3}, which shows the full solution profiles of $\rho$, $m$, and $z$, together with an enlarged view of $\rho$ near the discontinuity. Taken together, these results illustrate the influence of the relaxation parameter on the numerical solutions and examine whether the effects of $\alpha$ persist across different relaxation regimes.

Overall, the numerical results demonstrate that the parameter $\alpha$ influences the resolution of discontinuities in the numerical solution. Among the values considered, $\alpha=1/10$ exhibits the closest agreement with the fine-grid IMEX-RK2 reference solution near the discontinuity, as observed from the computed profiles. The choices $\alpha=1/6$ and $\alpha=1/4$ also provide satisfactory agreement with the reference solution and perform comparably to UCS2 in the test cases considered. These findings highlight the role of $\alpha$ in controlling the numerical behavior of the proposed scheme for nonsmooth solutions across different relaxation regimes.
\subsubsection{Euler Equations with Heat Transfer}
\label{subsec:euler_heat}
We next consider the one-dimensional Euler system with heat transfer to investigate the performance of the proposed numerical scheme for compressible flows with thermal relaxation. In the absence of
gravitational forces, the governing equations are written in conservative form as
\begin{equation}\label{Euler_heat}
\begin{cases}
\rho_t+(\rho u)_x=0,\\[2mm]
(\rho u)_t+(\rho u^2+p)_x=0,\\[2mm]
(\rho E)_t+(\rho uE+pu)_x=-K\rho(T-T_0),
\end{cases}
\end{equation}
where $\rho$, $u$, $p$, and $T$ denote the density, velocity, pressure, and temperature, respectively. The total energy per unit mass is defined as $E=e+\frac{u^2}{2}$, where $e$ is the specific internal energy. The heat-transfer coefficient is given by $K=1/\varepsilon$, where $\varepsilon>0$ is the thermal relaxation parameter, and $T_0$ denotes the prescribed ambient temperature. For a $\gamma$-law gas, the constitutive relations are
\[
p=(\gamma-1)\rho e,\qquad
e=c_vT=\frac{\mathcal R T}{\gamma-1},
\]
where $\gamma$ is the ratio of specific heats, $c_v$ is the specific heat at constant volume, and $\mathcal R$ is the specific gas constant. Introducing the conservative state vector
\[
Q=(\rho,m, z)^T,\quad
m=\rho u,\quad z=\rho E,
\]
the system can be written compactly as
\begin{equation}\label{Euler_heat_compact}
Q_t+\mathrm F(Q)_x=\mathrm S(Q,\varepsilon),
\end{equation}
where the flux function and the heat-transfer source term are
\[
\mathrm F(Q)=
\begin{bmatrix}
\rho u\\
\left(1-\dfrac{\gamma -1}{2}\right)\dfrac{m^2}{\rho}+(\gamma -1)z \\
\dfrac{m}{\rho}\left( \gamma z - \left(\dfrac{\gamma -1}{2}\right)\dfrac{m^2}{\rho}\right)
\end{bmatrix},
\quad
\mathrm S(Q,\varepsilon)=
\begin{bmatrix}
0\\
0\\
\dfrac{1}{\varepsilon}\left(\rho T_0 -\dfrac{1}{c_v} \left(z - \dfrac{m^2}{2\rho}\right)\right)
\end{bmatrix}.
\]
\medskip

\noindent\textbf{Test: Discontinuous Initial Data.}
We consider a discontinuous initial-value problem to assess the ability of the proposed scheme to capture sharp transitions in the presence of thermal relaxation. The initial conditions for the density, velocity, and total energy per unit mass are prescribed as
\begin{equation}\label{discontinuous_euler}
(\rho(x,0),u(x,0),E(x,0))=
\begin{cases}
(1,0,1), & x\leq 0.5,\\
(1/8,0,1), & x>0.5.
\end{cases}
\end{equation}
These data introduce a jump discontinuity at $x=0.5$. The computational domain is discretized using $N=200$ uniformly distributed grid points, with transmissive boundary conditions imposed at both endpoints.
\begin{figure}[ht]
    \centering
    % First row
    \begin{subfigure}[t]{0.48\textwidth}
        \centering        \includegraphics[width=\linewidth,height=5cm,keepaspectratio]{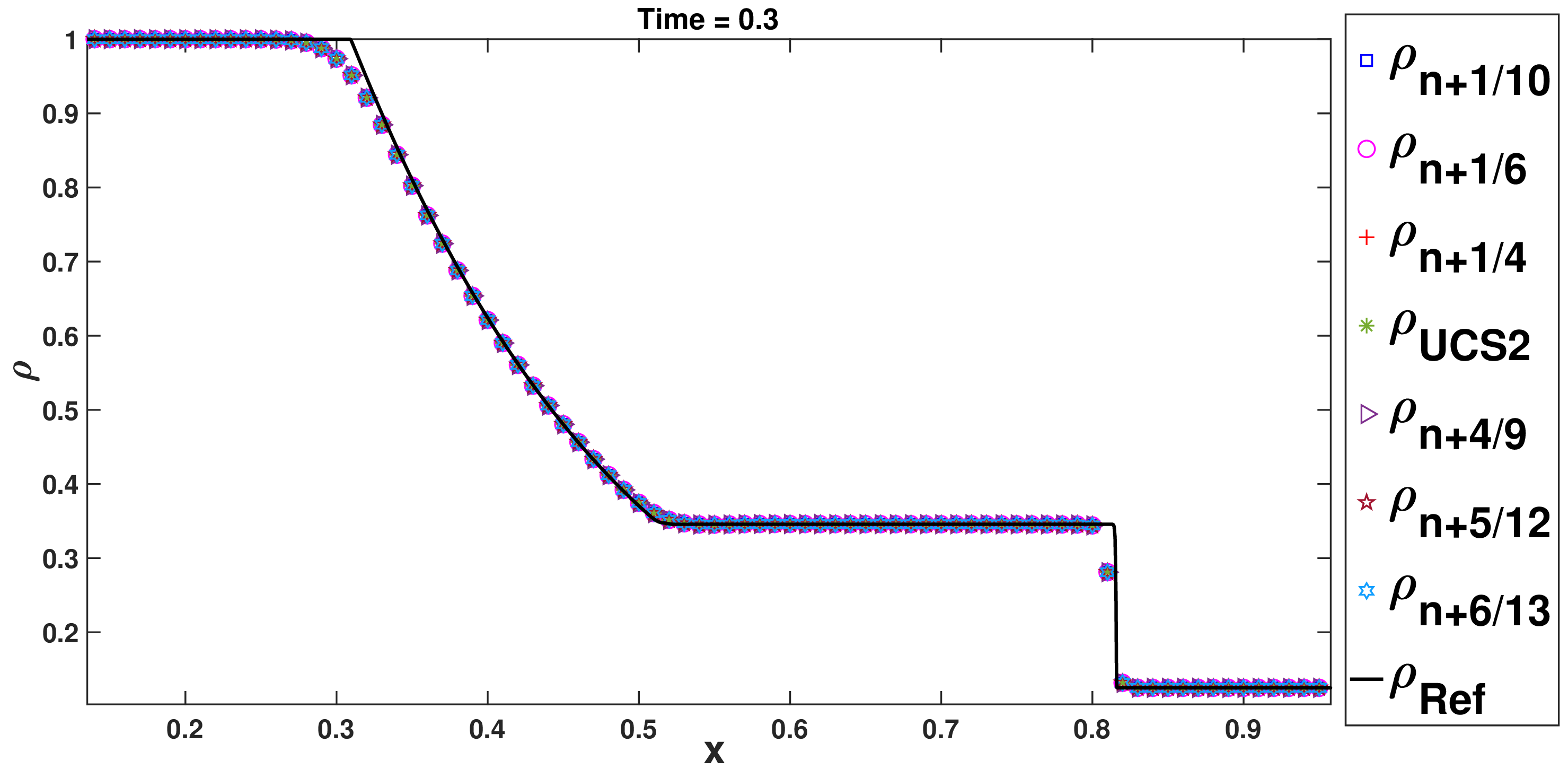}
        \caption{}
        \label{fig:euler-ep8-rho}
    \end{subfigure}
    \hfill
    \begin{subfigure}[t]{0.48\textwidth}
        \centering        \includegraphics[width=\linewidth,height=5cm,keepaspectratio]{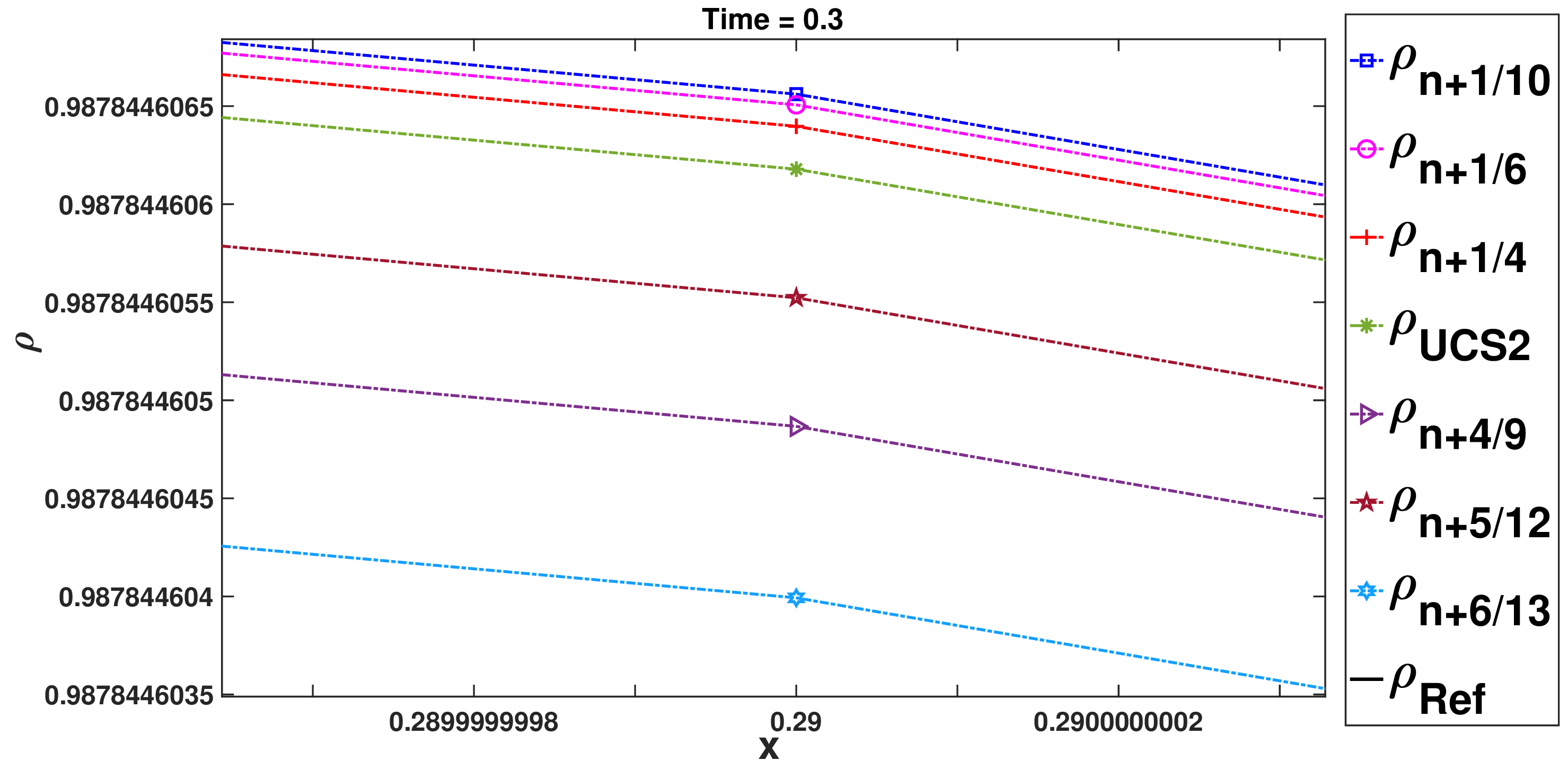}
        \caption{}
        \label{fig:euler-ep8-rho-zoom}
    \end{subfigure}
    
    \vspace{0.4cm}
    
    % Second row
    \begin{subfigure}[t]{0.48\textwidth}
        \centering        \includegraphics[width=\linewidth,height=5cm,keepaspectratio]{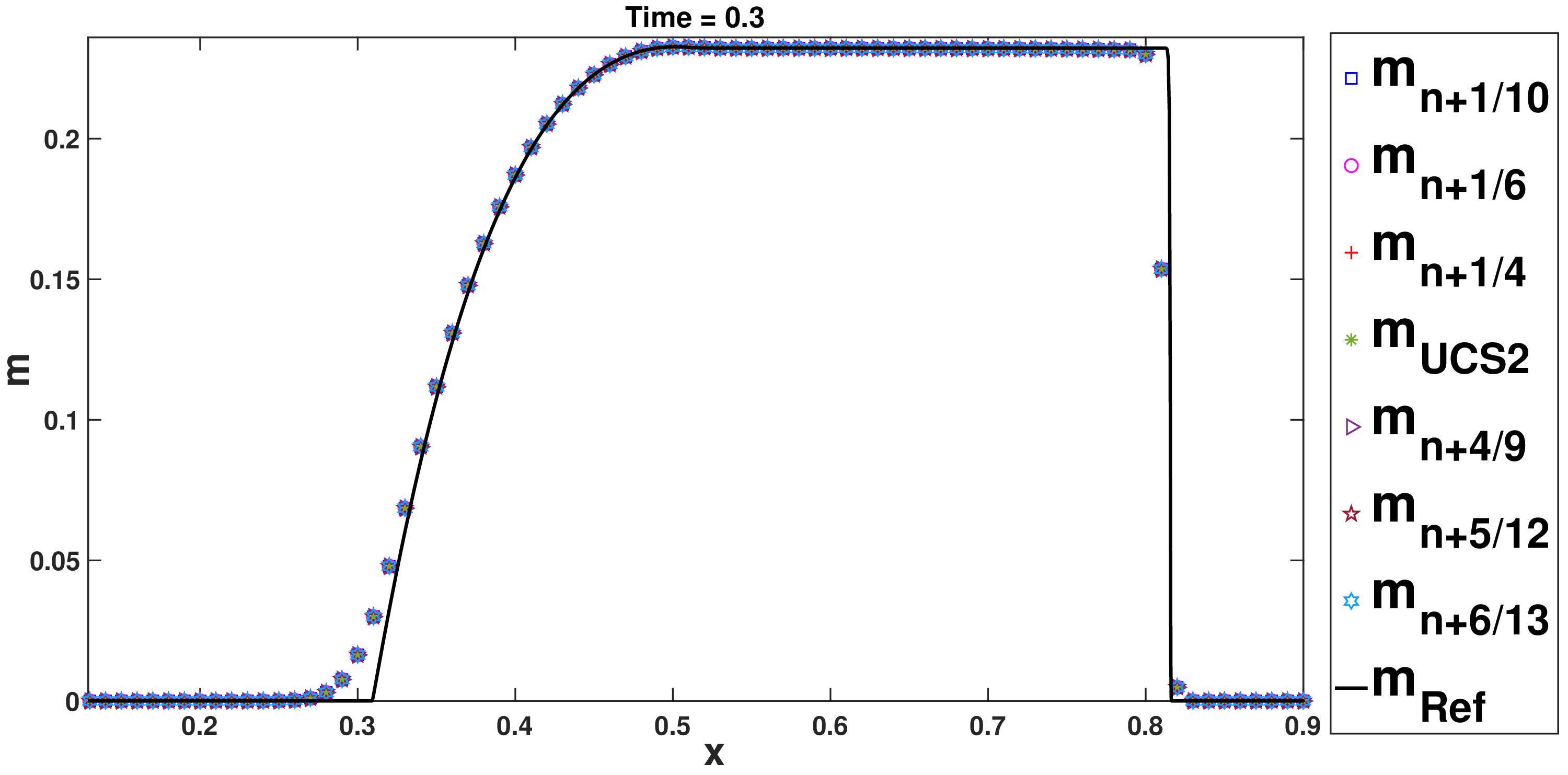}
        \caption{}
        \label{fig:euler-ep8-m}
    \end{subfigure}
    \hfill
    \begin{subfigure}[t]{0.48\textwidth}
        \centering        \includegraphics[width=\linewidth,height=5cm,keepaspectratio]{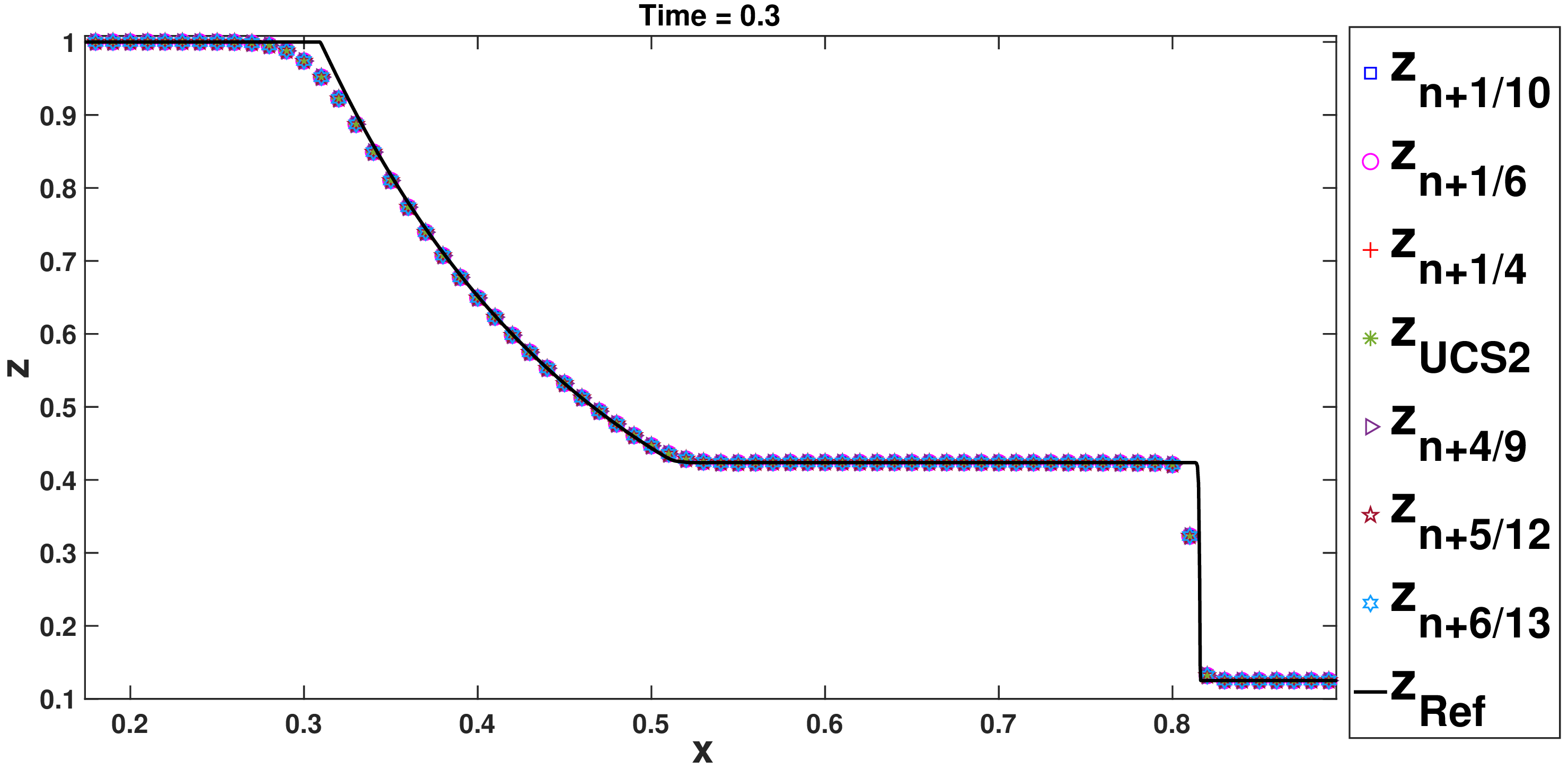}
        \caption{}
        \label{fig:euler-ep8-z}
    \end{subfigure}
    \caption{Numerical results for nonsmooth initial data in the Euler equations with heat transfer model, $\varepsilon = 10^{-8}$: (a) density $\rho$; (b) magnified view of $\rho$; (c) momentum $m$; and (d) variable $z$.}
    \label{nonsmooth_euler:ep8}
\end{figure}

To examine the effect of the parameter $\alpha$, we consider the same set of values as those used in the previous numerical experiments, as specified in \eqref{alpha:para}.
Figure~\ref{nonsmooth_euler:ep8} presents the numerical profiles of the density $\rho$, momentum $m=\rho u$, and total energy density $z=\rho E$ for the different values of $\alpha$. The results obtained using the proposed scheme are compared with those produced by UCS2, while a fine-grid IMEX-RK2 solution serves as a reference. An enlarged view of the density near the discontinuity is also provided to facilitate a closer examination of the local solution behavior.

The enlarged density profiles indicate that, among the parameter values considered, $\alpha=1/10$ provides the closest visual agreement with the fine-grid IMEX-RK2 reference solution near the discontinuity. Its numerical profile exhibits a sharper transition and smaller apparent deviations from the reference than those obtained with the other tested values of $\alpha$. The comparison with UCS2 also indicates that the proposed scheme with $\alpha=1/10$ resolves the discontinuity more sharply in this test. These observations illustrate the influence of $\alpha$ on the numerical behavior of the proposed scheme for the Euler system with thermal relaxation.

\subsection{Two-Dimensional Test Cases}
In this subsection, we investigate the performance of the proposed numerical scheme for two-dimensional test problems based on the compressible Euler equations. To construct a relaxation formulation suitable for the numerical treatment of these equations, we employ the Jin-Xin relaxation approach, which introduces auxiliary variables to approximate the directional Euler fluxes. This reformulation leads to a hyperbolic relaxation system containing source terms that drive the auxiliary variables toward their equilibrium values. In the limit as the relaxation parameter tends to zero, the relaxation system recovers the original two-dimensional Euler equations. We use this formulation to assess the ability of the proposed scheme to capture discontinuities and resolve complex wave structures in the stiff relaxation regime. The corresponding model and numerical configurations are described below.
\subsubsection{Two-Dimensional Jin-Xin Relaxation-Type Euler Model}\label{sec:eulerxinjin2d}
We consider the two-dimensional Jin-Xin relaxation system associated with the Euler equations, following the formulation in \cite{SW_friction}:
\begin{equation}
\begin{cases}
U_t + V_x + W_y = 0,\\[1mm]
V_t + A U_x = -\dfrac{1}{\varepsilon}
\left(V-F_1(U)\right),\\[2mm]
W_t + B U_y = -\dfrac{1}{\varepsilon}
\left(W-F_2(U)\right),
\end{cases}
\end{equation}
where $\varepsilon>0$ denotes the relaxation parameter, and $A$ and $B$ are diagonal matrices whose diagonal entries are $a_i$ and $b_i$, respectively. For the numerical experiments, we adopt the values
\[
a_1=a_2=a_3=a_4=9.75,
\qquad
b_1=b_2=b_3=b_4=7.62.
\]
The vector of conserved variables is defined as
\[
U=
\begin{bmatrix}
\rho\\
\rho u\\
\rho v\\
E
\end{bmatrix},
\]
where $\rho$ is the density, $u$ and $v$ denote the velocity components in the $x$- and $y$-directions, respectively, and $E$ represents the total energy. The corresponding Euler fluxes in the two coordinate directions are
\[
F_1(U)=
\begin{bmatrix}
\rho u\\
\rho u^2+p\\
\rho uv\\
(E+p)u
\end{bmatrix},
\qquad
F_2(U)=
\begin{bmatrix}
\rho v\\
\rho uv\\
\rho v^2+p\\
(E+p)v
\end{bmatrix}.
\]
The pressure is determined by the ideal-gas equation of state,
\[
p=(\gamma-1)
\left(E-\frac12\rho(u^2+v^2)\right),
\]
where $\gamma$ denotes the ratio of specific heats.

\begin{figure}[htb]
    \centering
    \begin{subfigure}[t]{0.325\textwidth}
        \centering        \includegraphics[width=\linewidth,height=5.5cm,trim=12.0cm 0cm 11.5cm 0cm,clip]{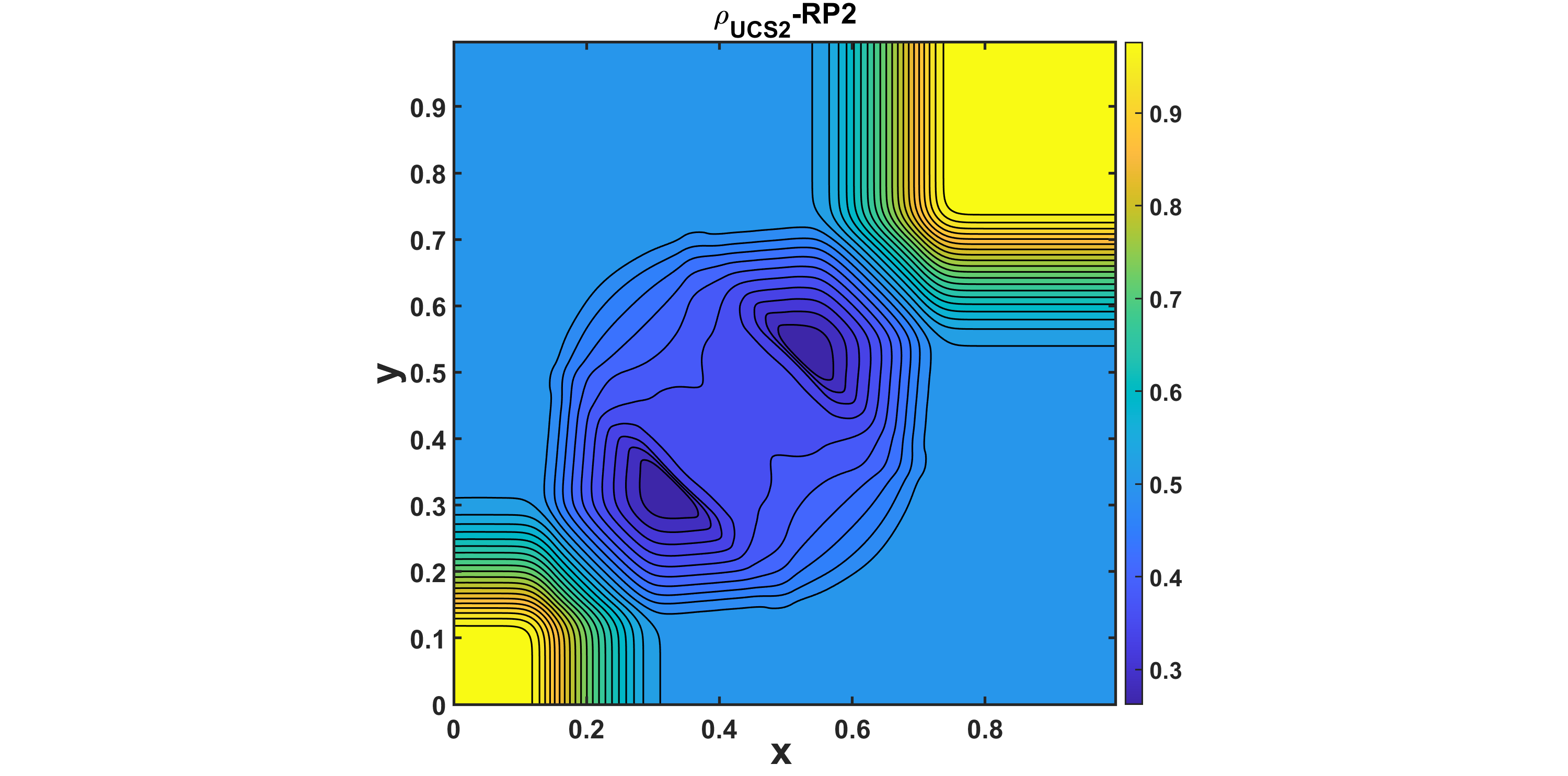}
        \caption{}
        \label{fig:rp2-ucs2}
    \end{subfigure}
    \hfill
    \begin{subfigure}[t]{0.325\textwidth}
        \centering        \includegraphics[width=\linewidth,height=5.5cm,trim=12.0cm 0cm 11.5cm 0cm,clip]{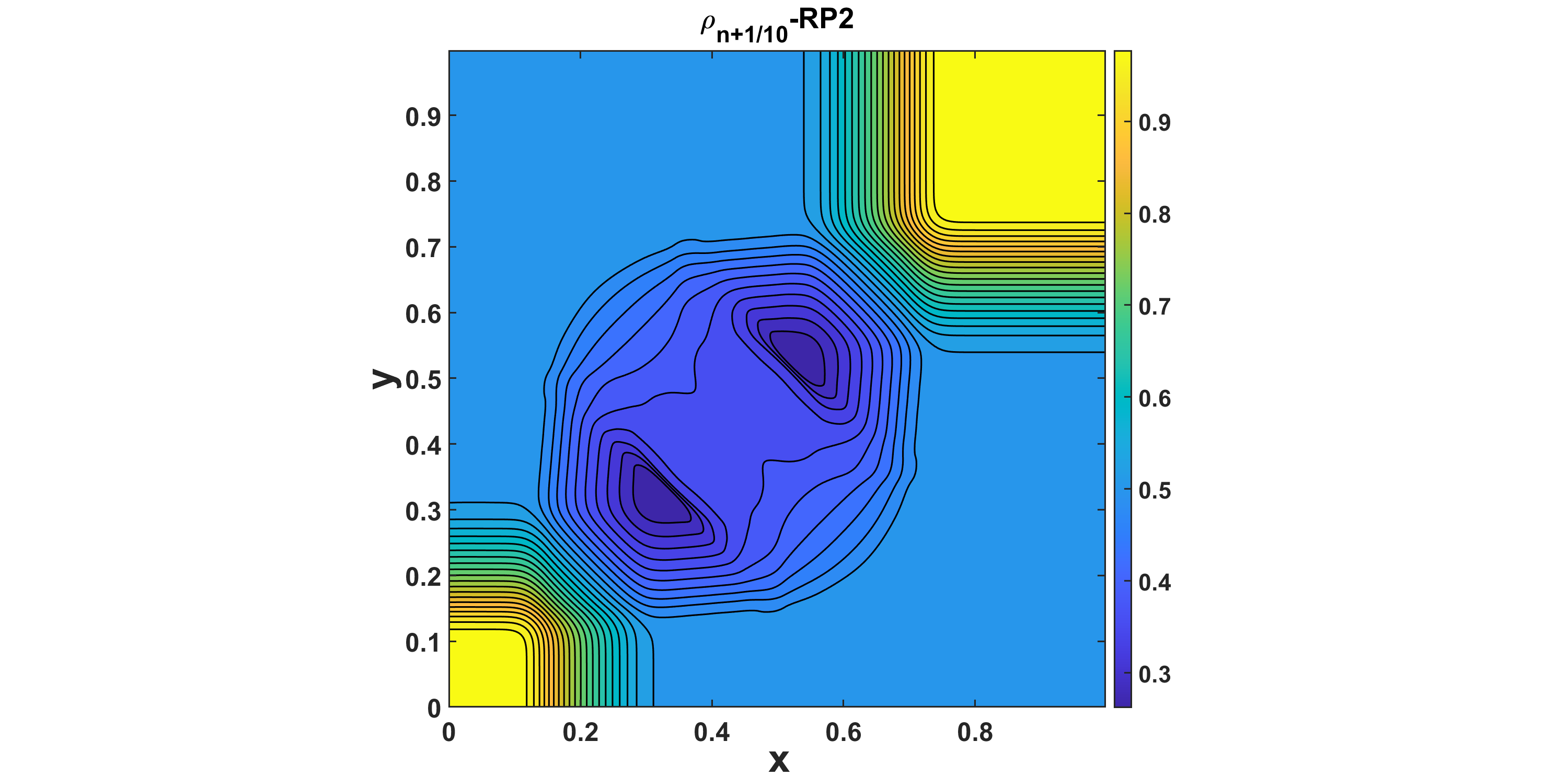}
        \caption{}
        \label{fig:rp2-proposed}
    \end{subfigure}
    \hfill
    \begin{subfigure}[t]{0.325\textwidth}
        \centering        \includegraphics[width=\linewidth,height=5.5cm,trim=12.0cm 0cm 11.5cm 0cm,clip]{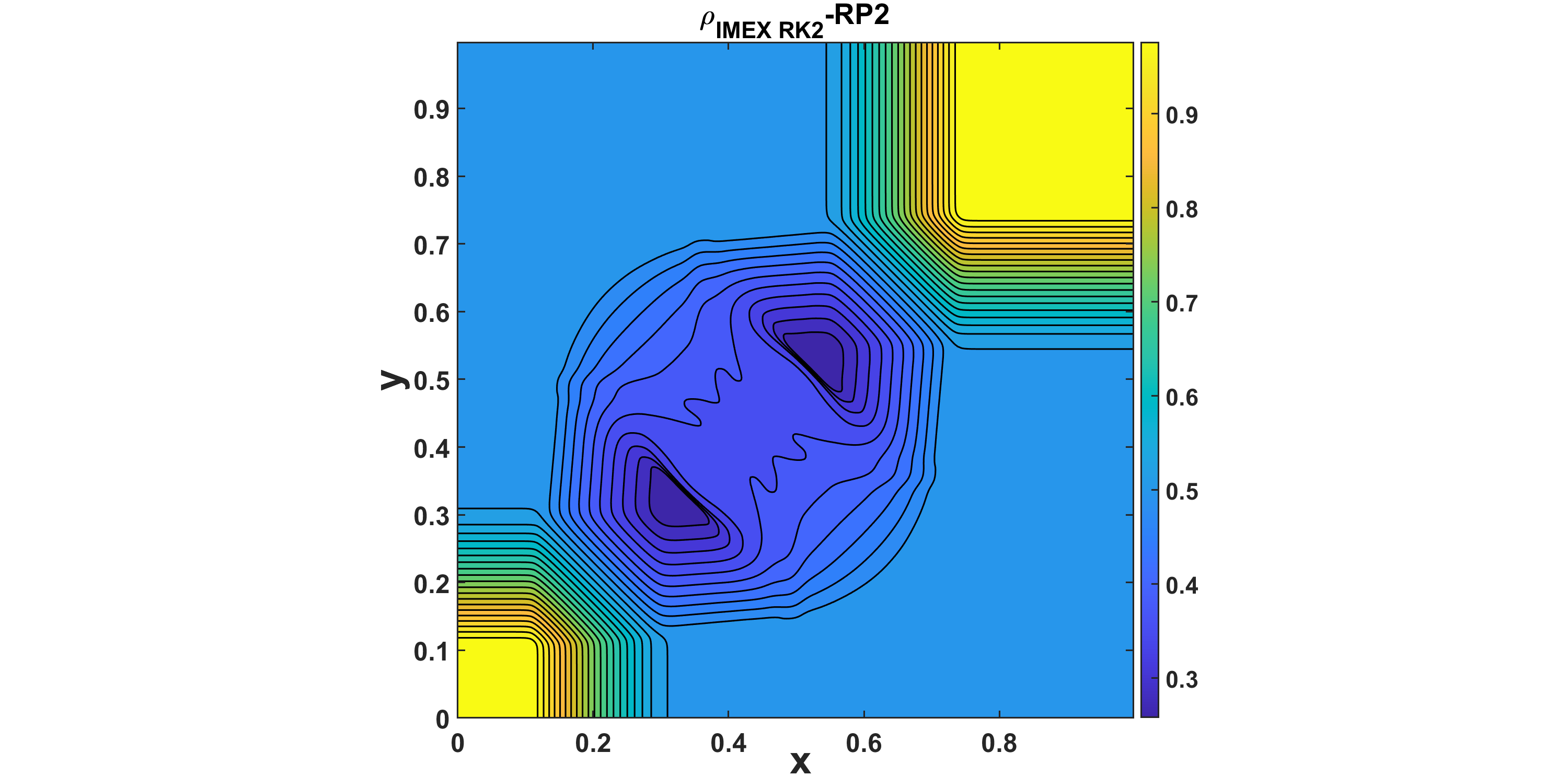}
        \caption{}
        \label{fig:rp2-imex}
    \end{subfigure}
    \caption{Numerical solutions for the density $\rho$ in Configuration 2, obtained using (a) the UCS2 scheme, (b) the proposed scheme with $\alpha=1/10$, and (c) the IMEX-RK2 scheme.}
    \label{/RP:euler:ep2}
\end{figure}

\begin{figure}[htb]
    \centering
    \begin{subfigure}[t]{0.325\textwidth}
        \centering        \includegraphics[width=\linewidth,height=5.5cm,trim=11.5cm 0cm 11.5cm 0cm,clip]{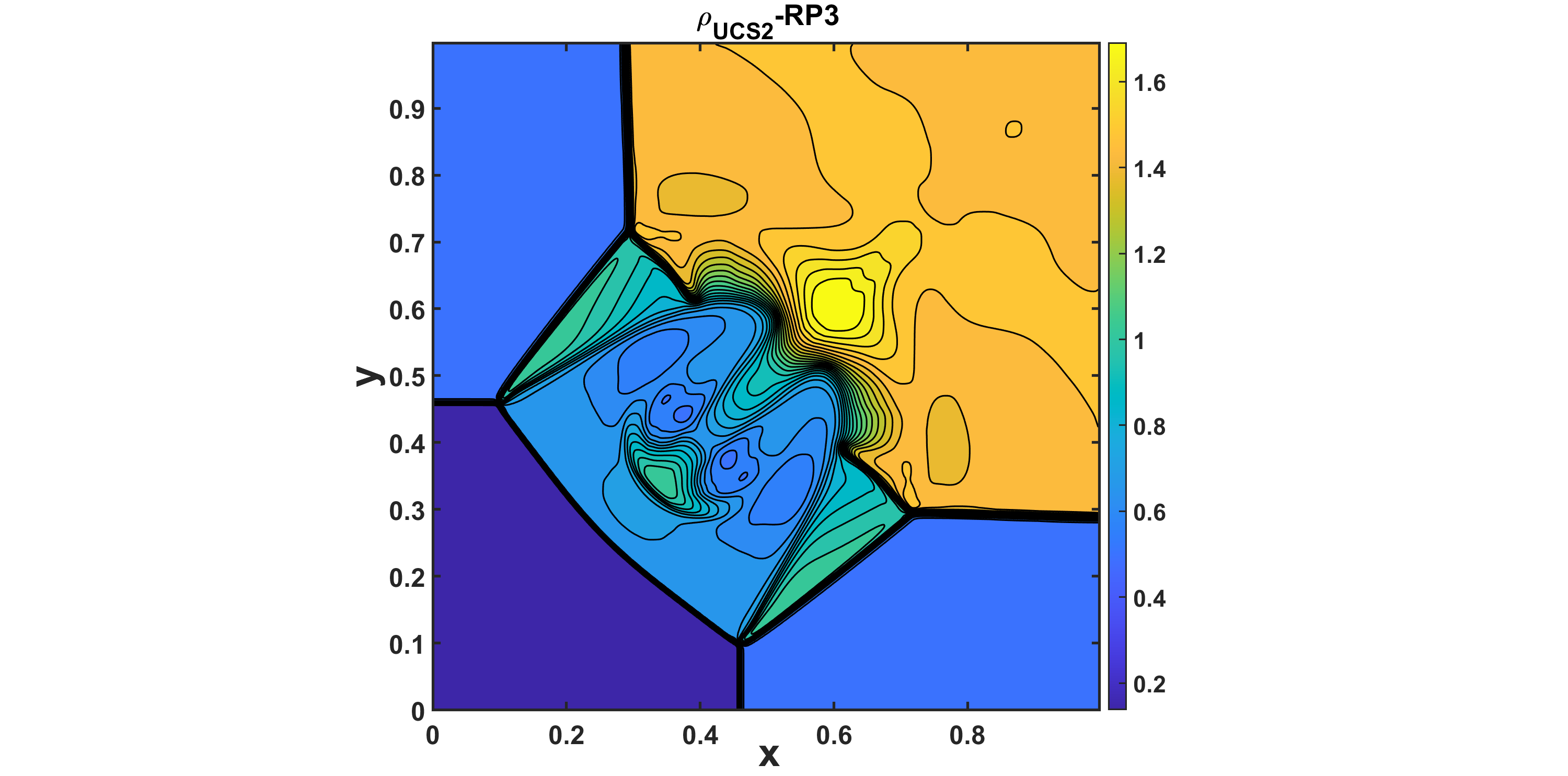}
        \caption{}
        \label{fig:rp3-ucs2}
    \end{subfigure}
    \hfill
    \begin{subfigure}[t]{0.325\textwidth}
        \centering        \includegraphics[width=\linewidth,height=5.5cm,trim=12.0cm 0cm 11.5cm 0cm,clip]{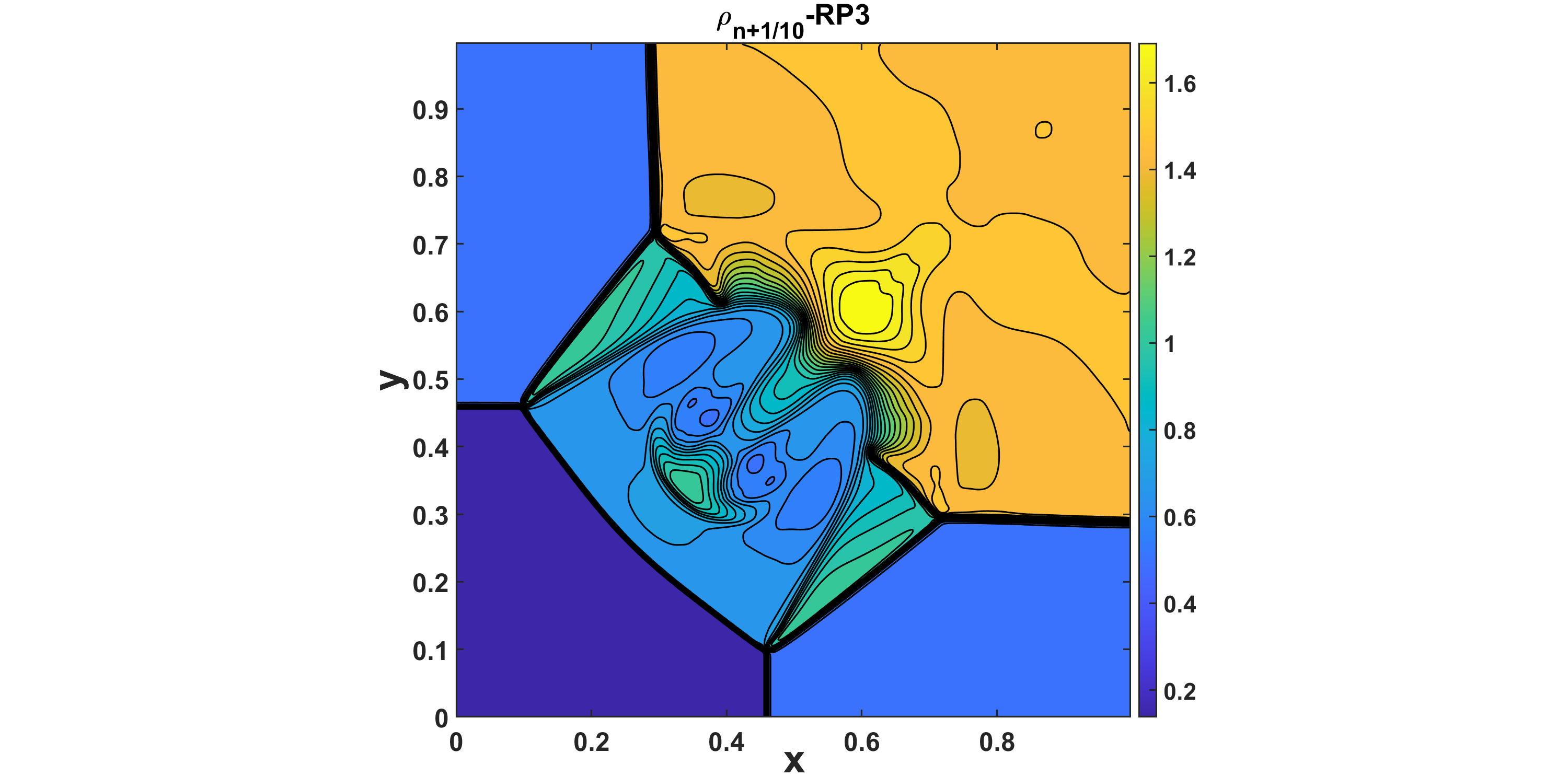}
        \caption{}
        \label{fig:rp3-proposed}
    \end{subfigure}
    \hfill
    \begin{subfigure}[t]{0.325\textwidth}
        \centering        \includegraphics[width=\linewidth,height=5.5cm,trim=12.0cm 0cm 10.7cm 0cm,clip]{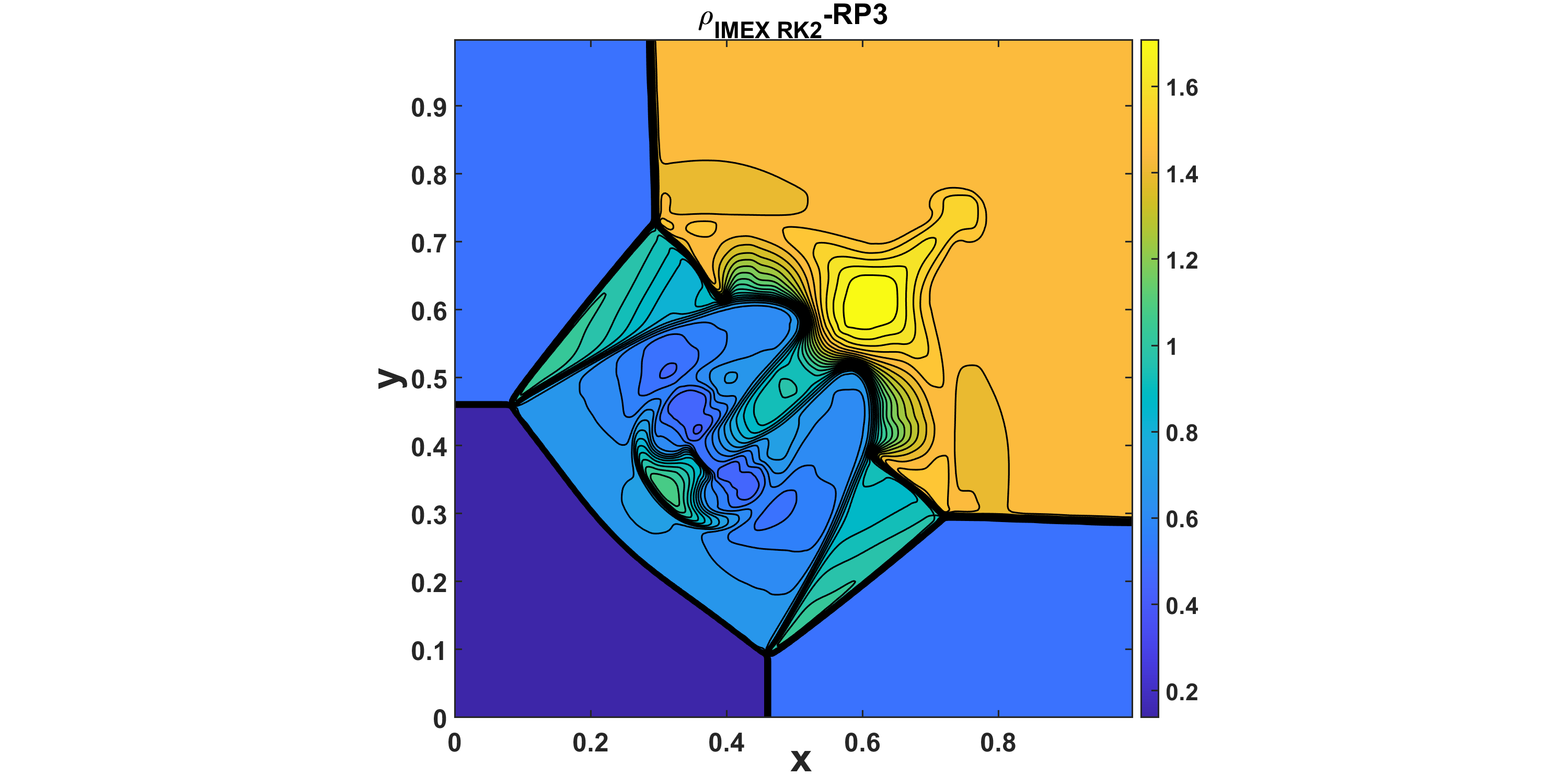}
        \caption{}
        \label{fig:rp3-imex}
    \end{subfigure}
    \caption{Comparison of the numerical density profiles $\rho$ for Configuration 3 obtained using (a) the UCS2 scheme, (b) the proposed scheme with $\alpha=1/10$, and (c) the IMEX-RK2 scheme.}
    \label{/RP:euler:ep3}
\end{figure}

\begin{figure}[htb]
    \centering
    \begin{subfigure}[t]{0.325\textwidth}
        \centering        \includegraphics[width=\linewidth,height=5.5cm,trim=12.0cm 0cm 10.7cm 0cm,clip]{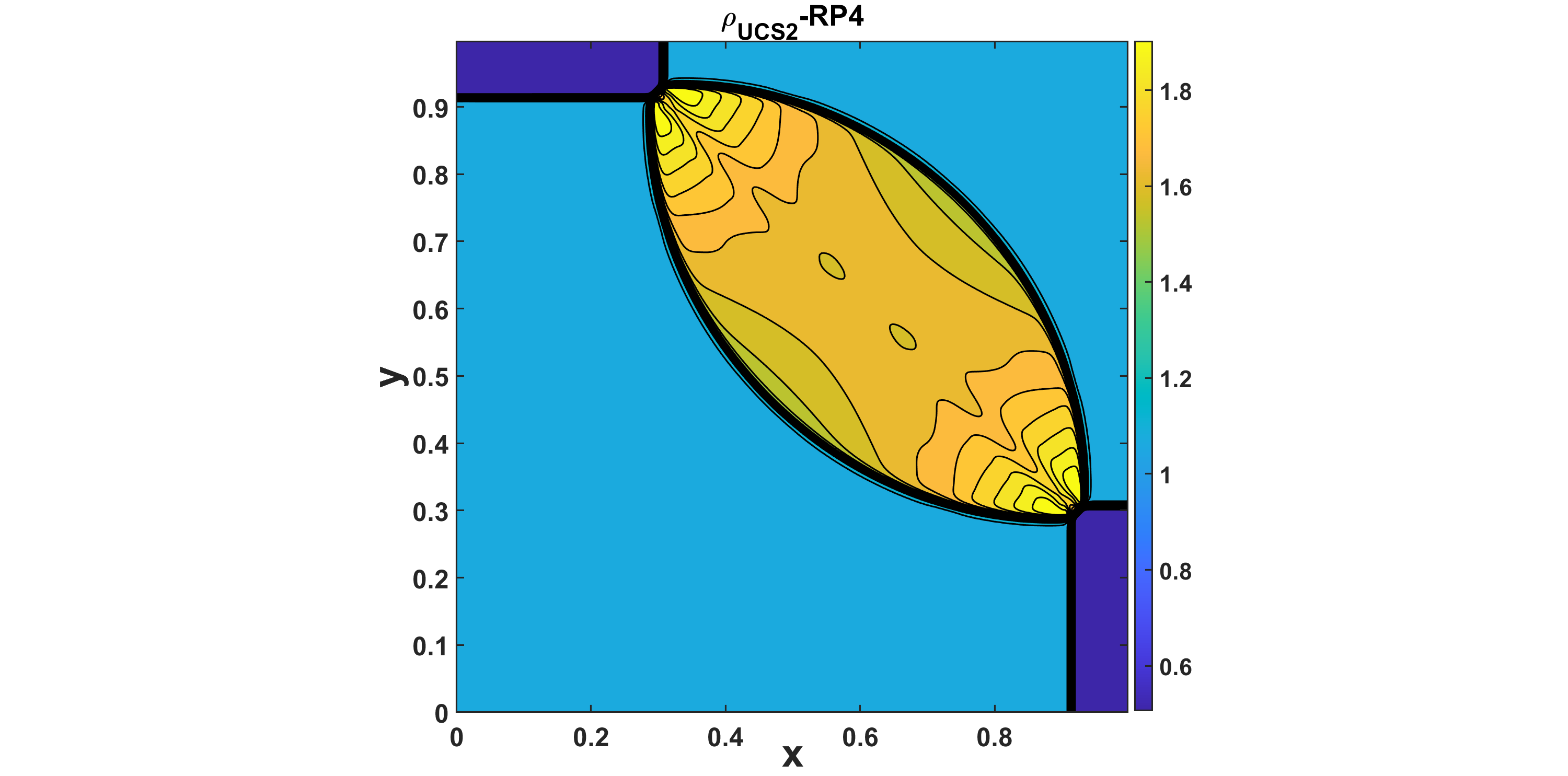}
        \caption{}
        \label{fig:rp4-ucs2}
    \end{subfigure}
    \hfill
    \begin{subfigure}[t]{0.325\textwidth}
        \centering        \includegraphics[width=\linewidth,height=5.5cm,trim=11.8cm 0cm 10.7cm 0cm,clip]{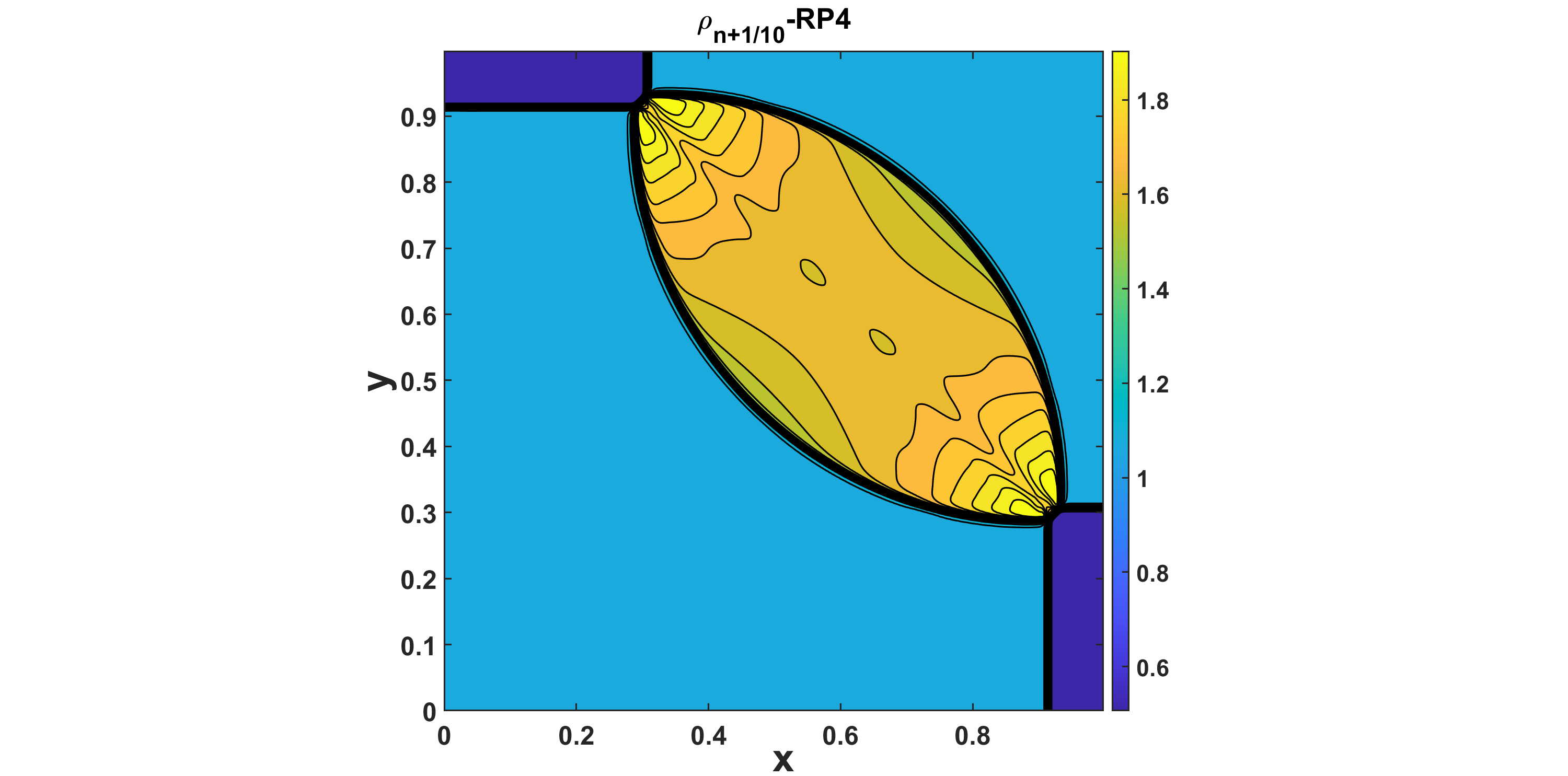}
        \caption{}
        \label{fig:rp4-proposed}
    \end{subfigure}
    \hfill
    \begin{subfigure}[t]{0.325\textwidth}
        \centering       \includegraphics[width=\linewidth,height=5.5cm,trim=12.0cm 0cm 10.7cm 0cm,clip]{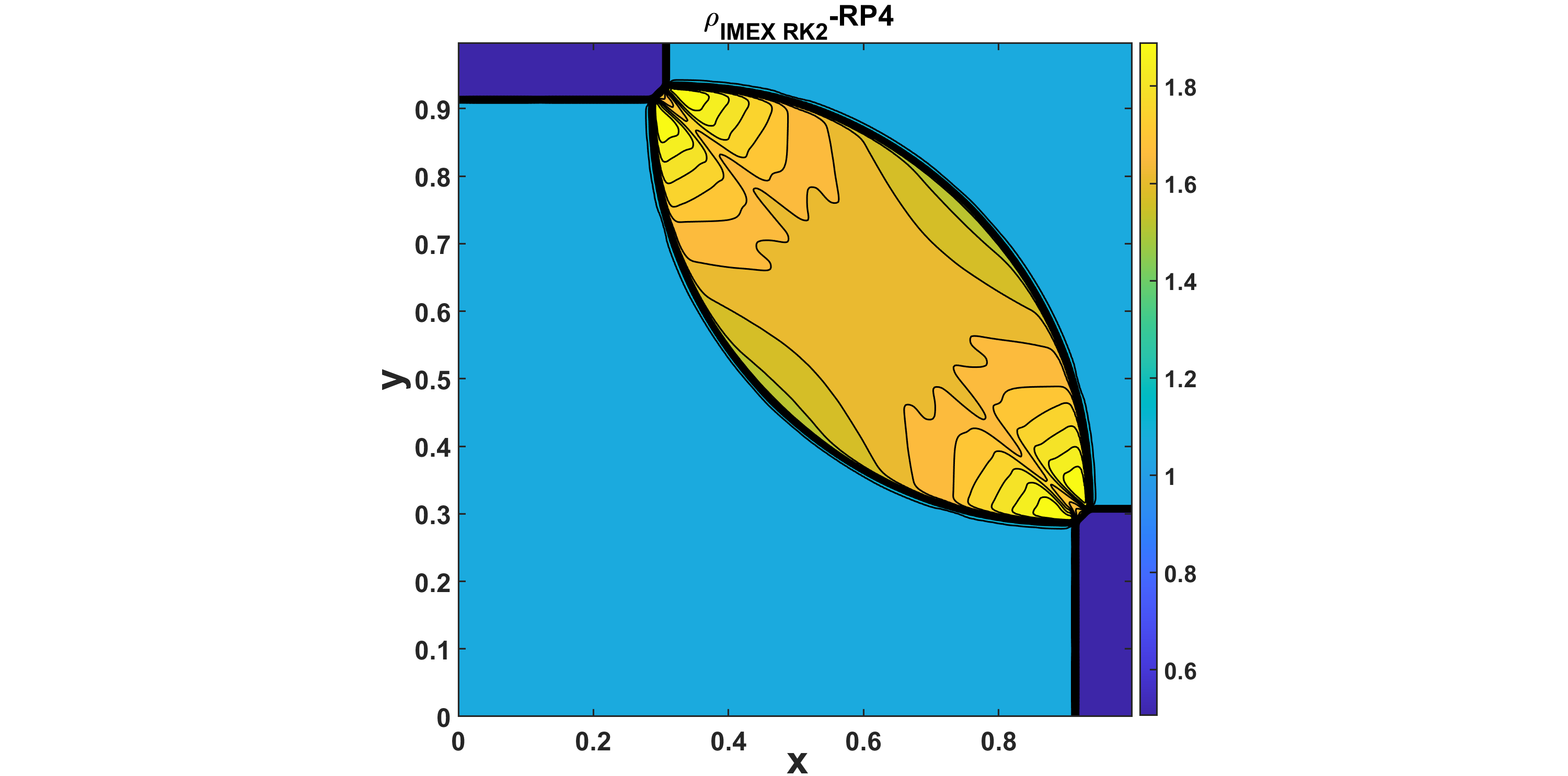}
        \caption{}
        \label{fig:rp4-imex}
    \end{subfigure}
    \caption{Comparison of the numerical density profiles $\rho$ for Configuration 4 obtained using (a) the UCS2 scheme, (b) the proposed scheme with $\alpha=1/10$, and (c) the IMEX-RK2 scheme.}
    \label{/RP:euler:ep4}
\end{figure}

\begin{figure}[htb]
    \centering
    \begin{subfigure}[t]{0.325\textwidth}
        \centering        \includegraphics[width=\linewidth,height=5.5cm,trim=12.0cm 0cm 10.7cm 0cm,clip]{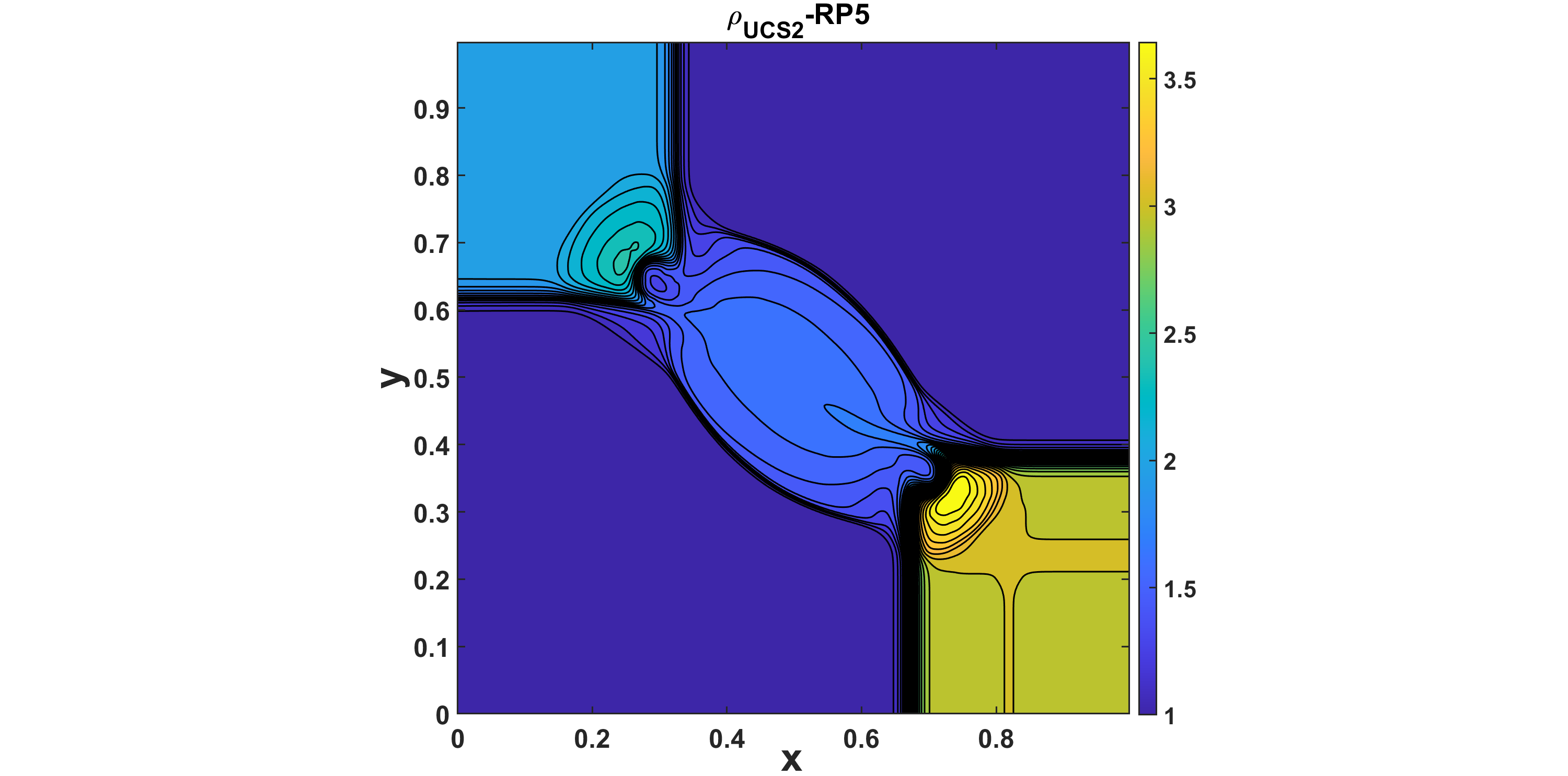}
        \caption{}
        \label{fig:rp5-ucs2}
    \end{subfigure}
    \hfill
    \begin{subfigure}[t]{0.325\textwidth}
        \centering        \includegraphics[width=\linewidth,height=5.5cm,trim=12.0cm 0cm 10.7cm 0cm,clip]{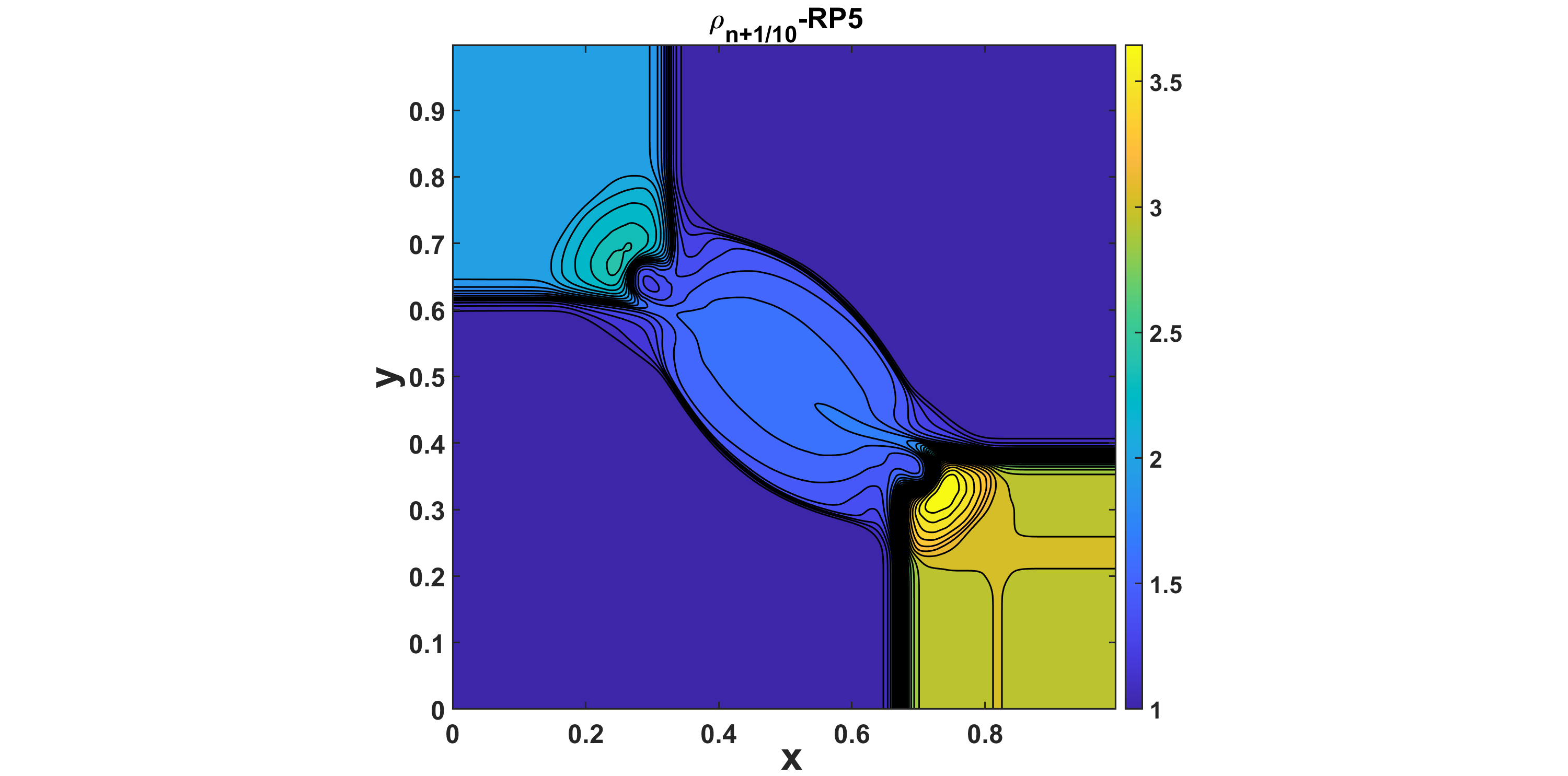}
        \caption{}
        \label{fig:rp5-proposed}
    \end{subfigure}
    \hfill
    \begin{subfigure}[t]{0.325\textwidth}
        \centering        \includegraphics[width=\linewidth,height=5.5cm,trim=12.0cm 0cm 10.7cm 0cm,clip]{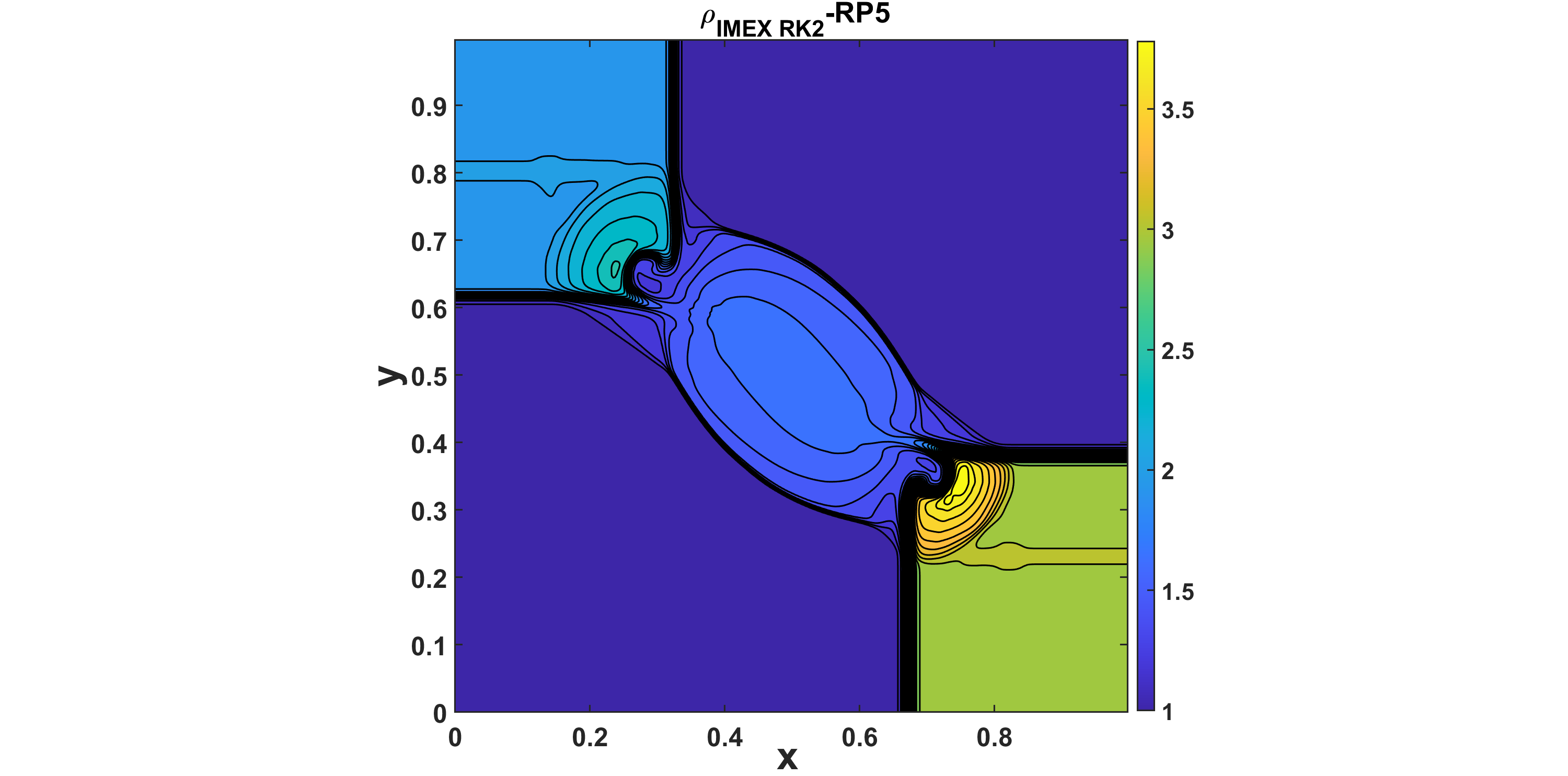}
        \caption{}
        \label{fig:rp5-imex}
    \end{subfigure}
    \caption{Comparison of the numerical density profiles $\rho$ for Configuration 5 obtained using (a) the UCS2 scheme, (b) the proposed scheme with $\alpha=1/10$, and (c) the IMEX-RK2 scheme.}
    \label{/RP:euler:ep5}
\end{figure}

\begin{figure}[htb]
    \centering
    \begin{subfigure}[t]{0.325\textwidth}
        \centering       \includegraphics[width=\linewidth,height=5.5cm,trim=12.5cm 0cm 10.7cm 0cm,clip]{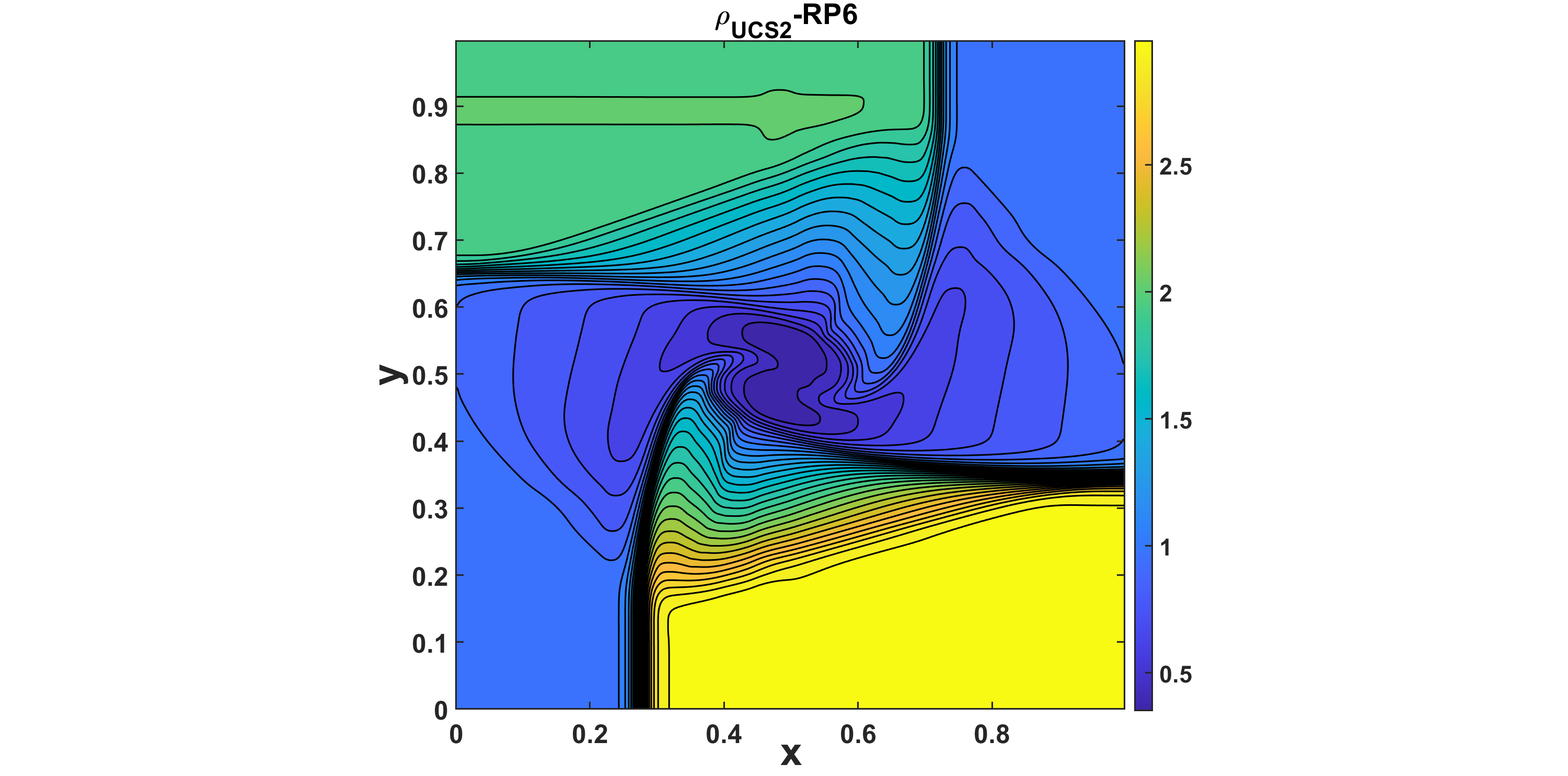}
        \caption{}
        \label{fig:rp6-ucs2}
    \end{subfigure}
    \hfill
    \begin{subfigure}[t]{0.325\textwidth}
        \centering        \includegraphics[width=\linewidth,height=5.5cm,trim=12.5cm 0cm 10.7cm 0cm,clip]{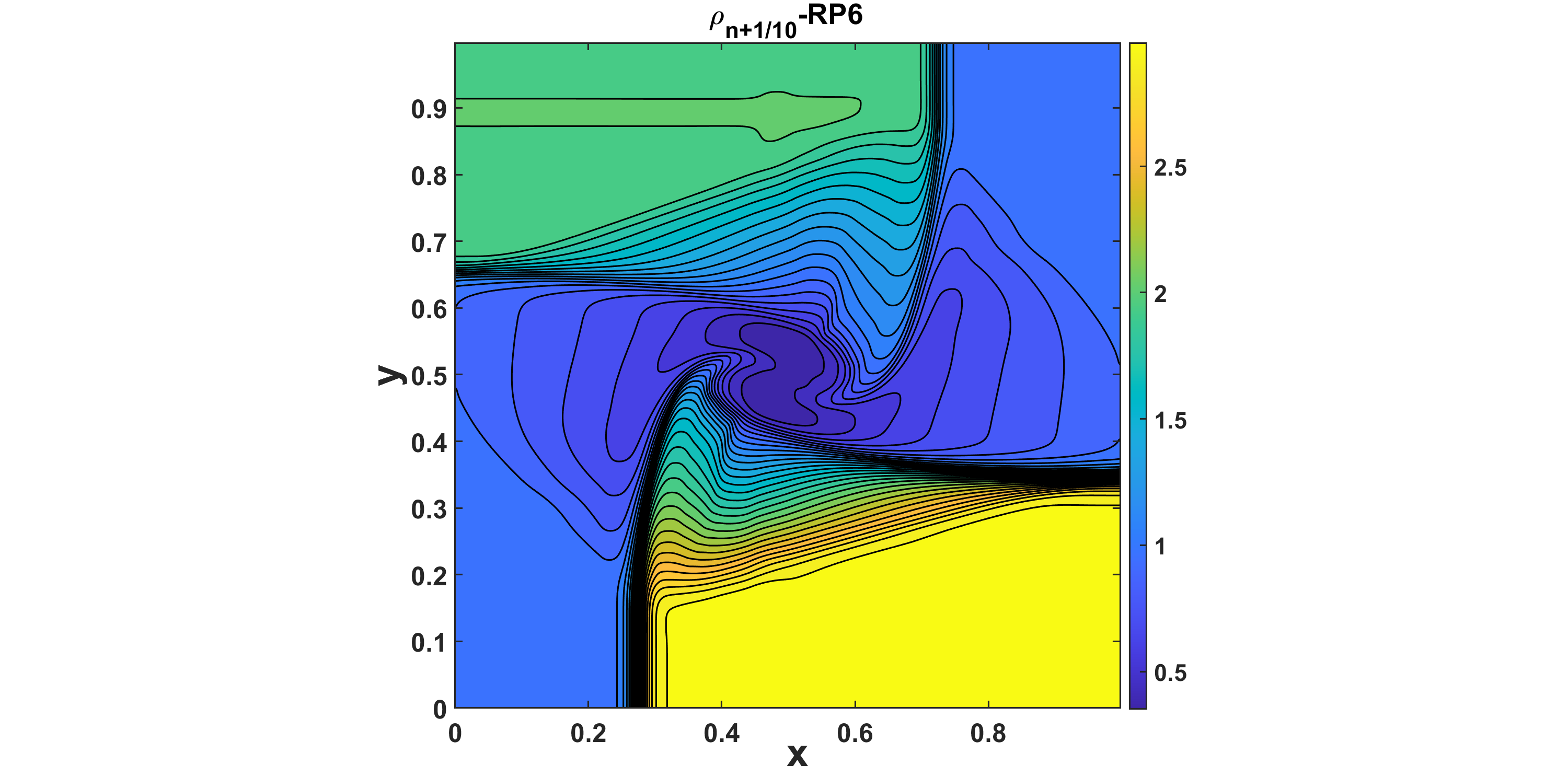}
        \caption{}
        \label{fig:rp6-proposed}
    \end{subfigure}
    \hfill
    \begin{subfigure}[t]{0.325\textwidth}
        \centering       \includegraphics[width=\linewidth,height=5.5cm,trim=12.5cm 0cm 10.7cm 0cm,clip]{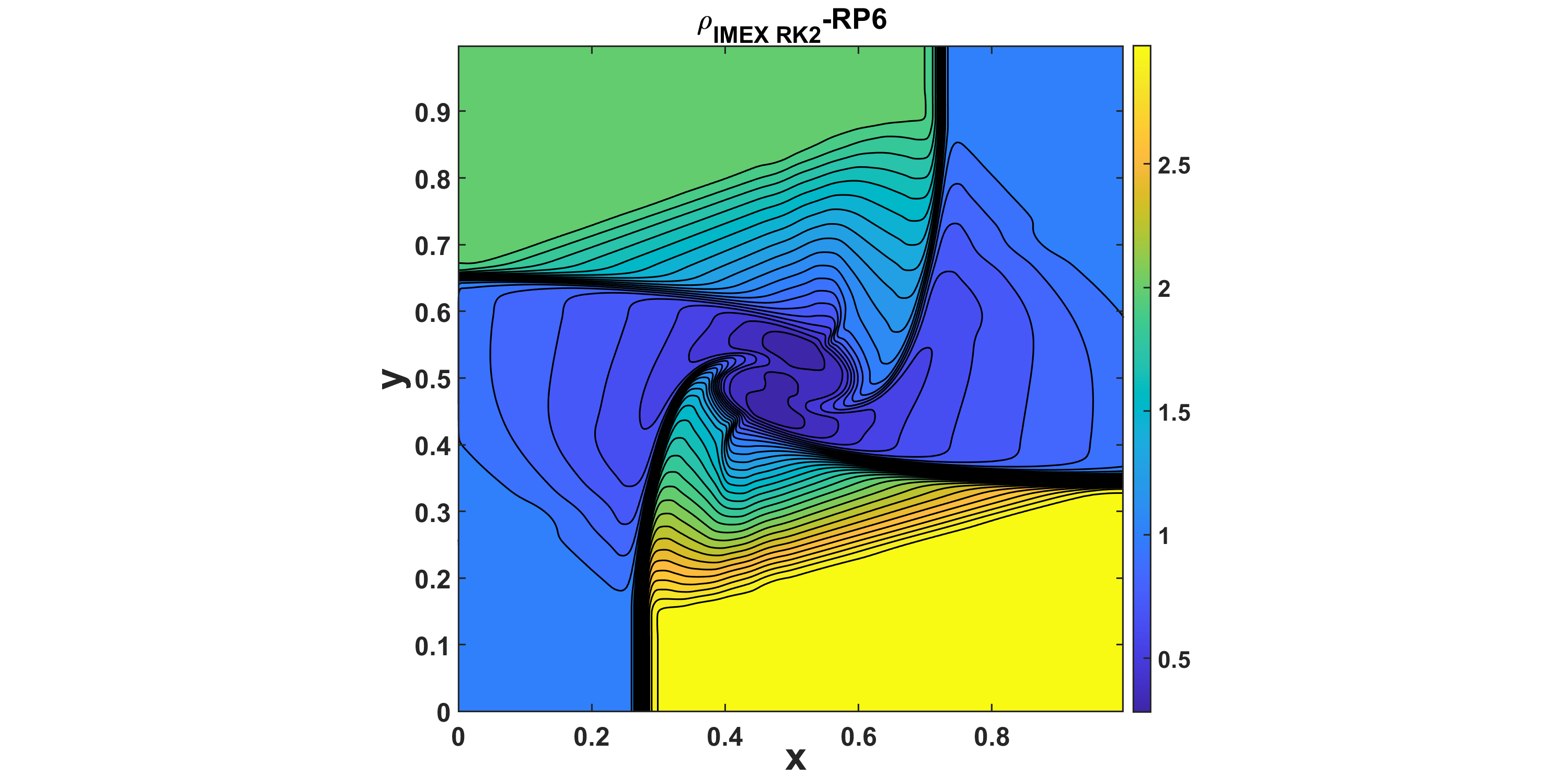}
        \caption{}
        \label{fig:rp6-imex}
    \end{subfigure}
    \caption{Comparison of the numerical density profiles $\rho$ for Configuration 6 obtained using (a) the UCS2 scheme, (b) the proposed scheme with $\alpha=1/10$, and (c) the IMEX-RK2 scheme.}
    \label{/RP:euler:ep6}
\end{figure}

\begin{figure}[htb]
    \centering
    \begin{subfigure}[t]{0.325\textwidth}
        \centering        \includegraphics[width=\linewidth,height=5.5cm,trim=12.5cm 0cm 10.7cm 0cm,clip]{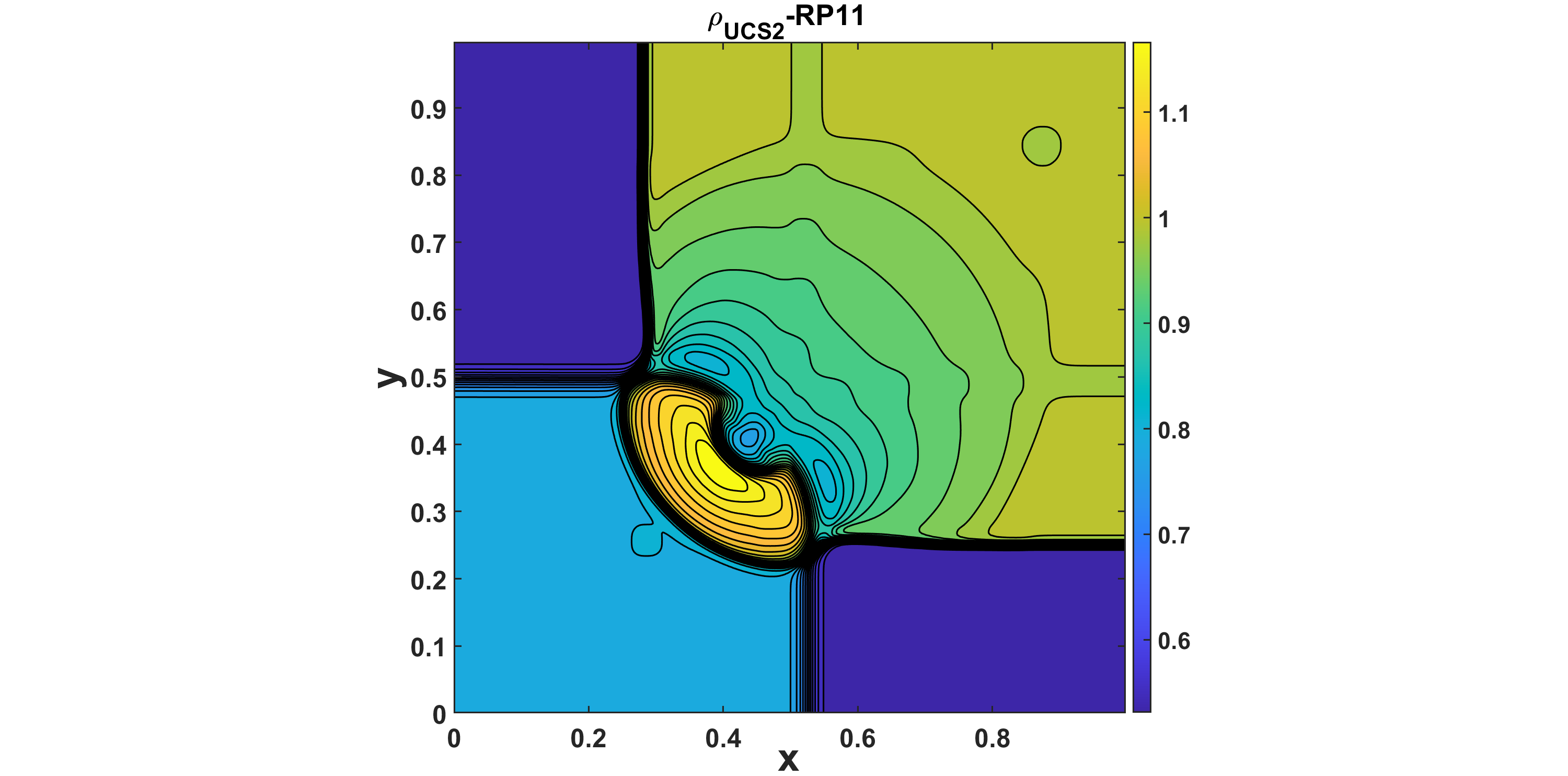}
        \caption{}
        \label{fig:rp11-ucs2}
    \end{subfigure}
    \hfill
    \begin{subfigure}[t]{0.325\textwidth}
        \centering        \includegraphics[width=\linewidth,height=5.5cm,trim=12.5cm 0cm 10.7cm 0cm,clip]{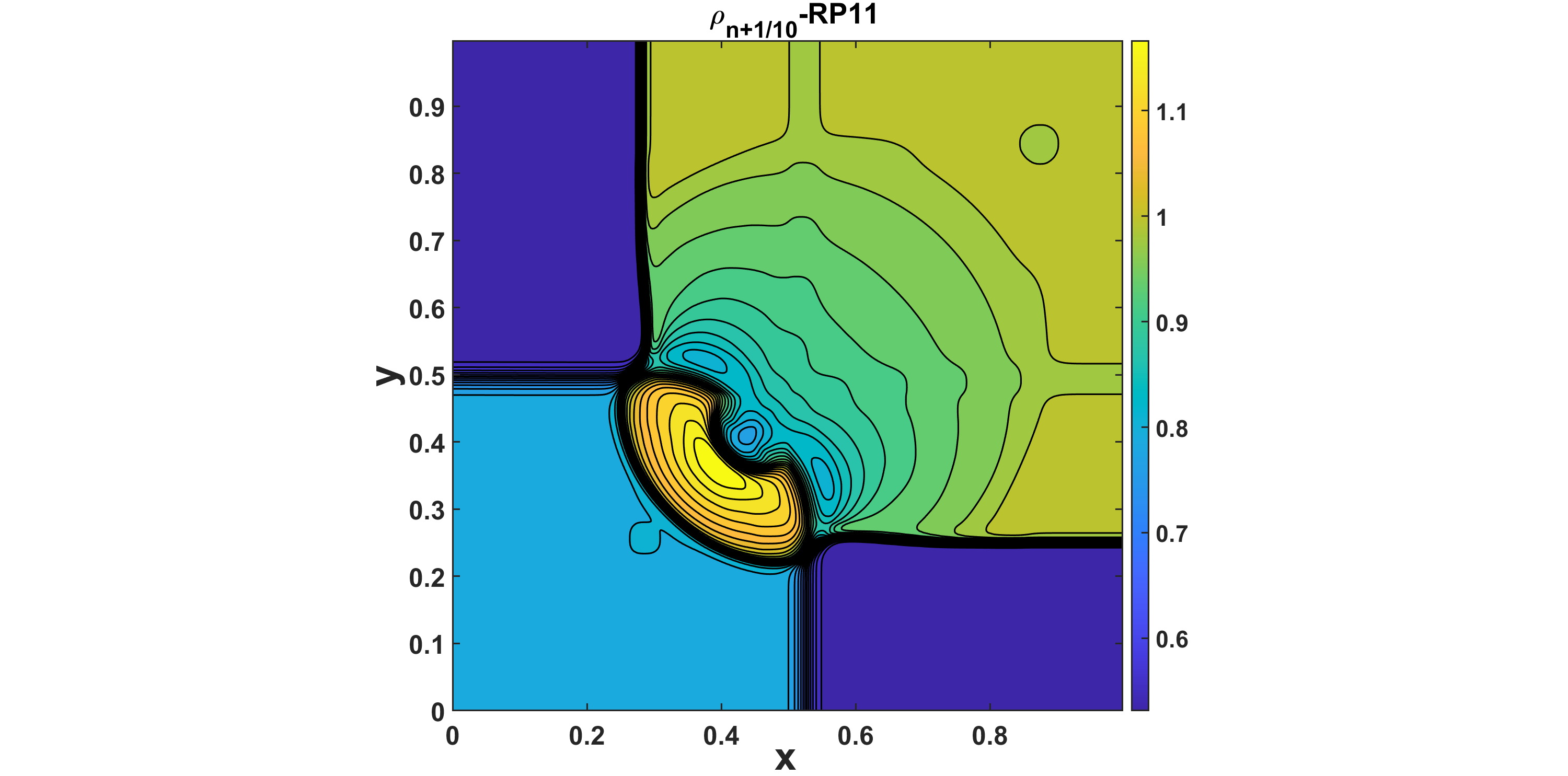}
        \caption{}
        \label{fig:rp11-general}
    \end{subfigure}
    \hfill
    \begin{subfigure}[t]{0.325\textwidth}
        \centering        \includegraphics[width=\linewidth,height=5.5cm,trim=12.5cm 0cm 10.7cm 0cm,clip]{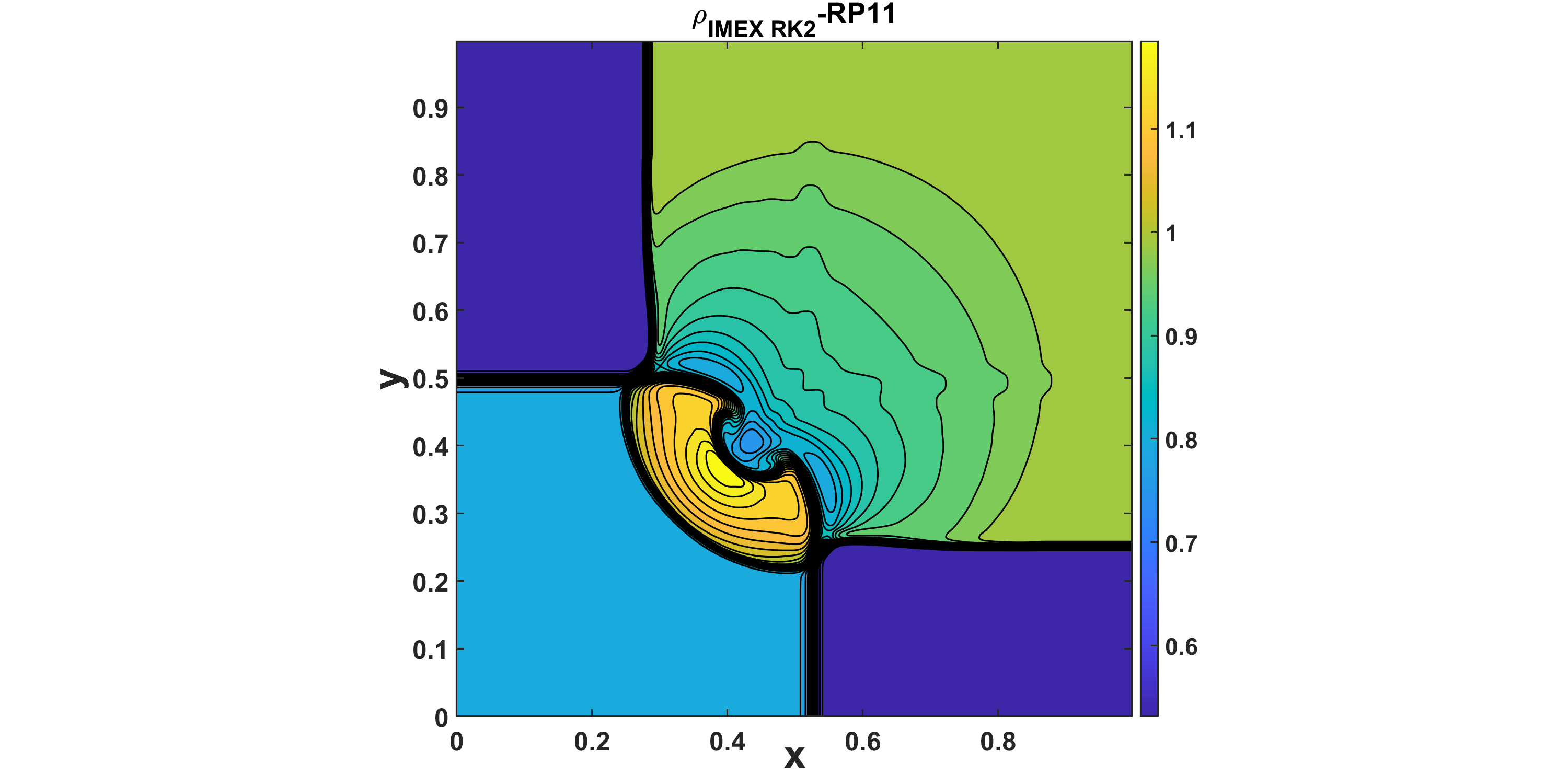}
        \caption{}
        \label{fig:rp11-imex}
    \end{subfigure}
    \caption{Comparison of the numerical density profiles $\rho$ for Configuration 11 obtained using (a) the UCS2 scheme, (b) the proposed scheme with $\alpha=1/10$, and (c) the IMEX-RK2 scheme.}
    \label{/RP:euler:ep11}
\end{figure}

\begin{figure}[htb]
    \centering
    \begin{subfigure}[t]{0.325\textwidth}
        \centering        \includegraphics[width=\linewidth,height=5.5cm,trim=12.5cm 0cm 10.7cm 0cm,clip]{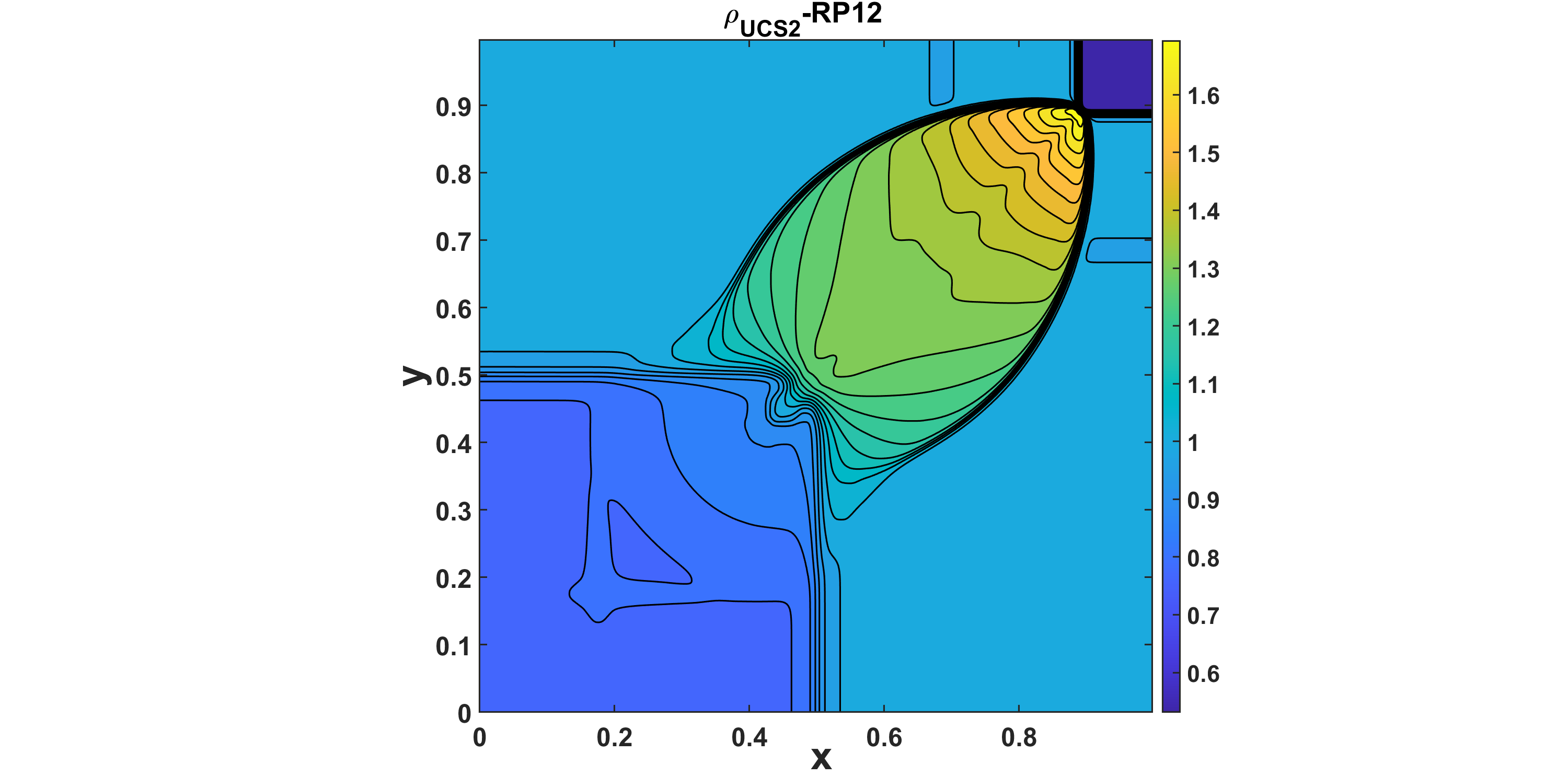}
        \caption{}
        \label{fig:rp12-ucs2}
    \end{subfigure}
    \hfill
    \begin{subfigure}[t]{0.325\textwidth}
        \centering        \includegraphics[width=\linewidth,height=5.5cm,trim=12.5cm 0cm 10.7cm 0cm,clip]{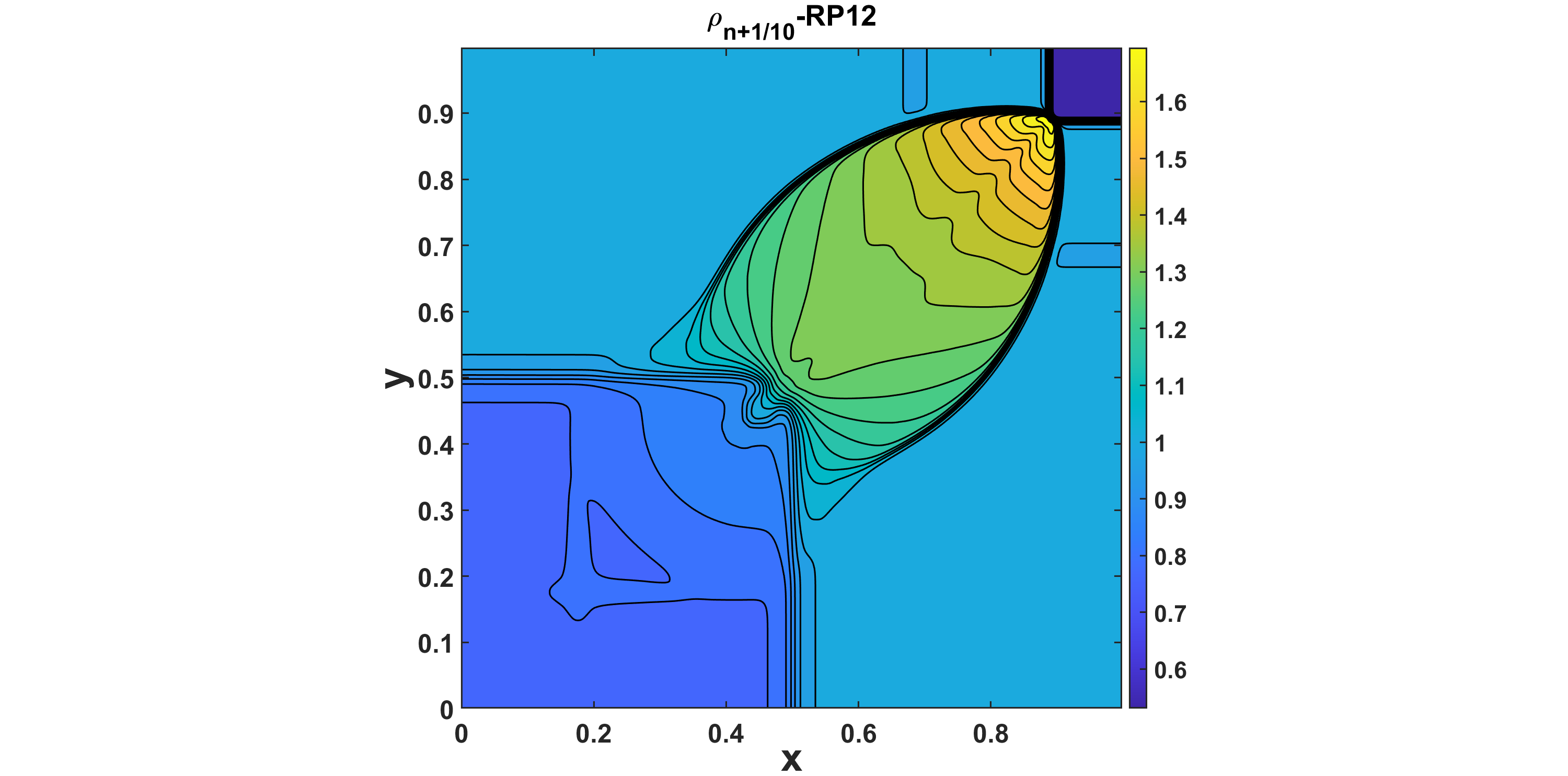}
        \caption{}
        \label{fig:rp12-proposed}
    \end{subfigure}
    \hfill
    \begin{subfigure}[t]{0.325\textwidth}
        \centering       \includegraphics[width=\linewidth,height=5.5cm,trim=12.5cm 0cm 10.7cm 0cm,clip]{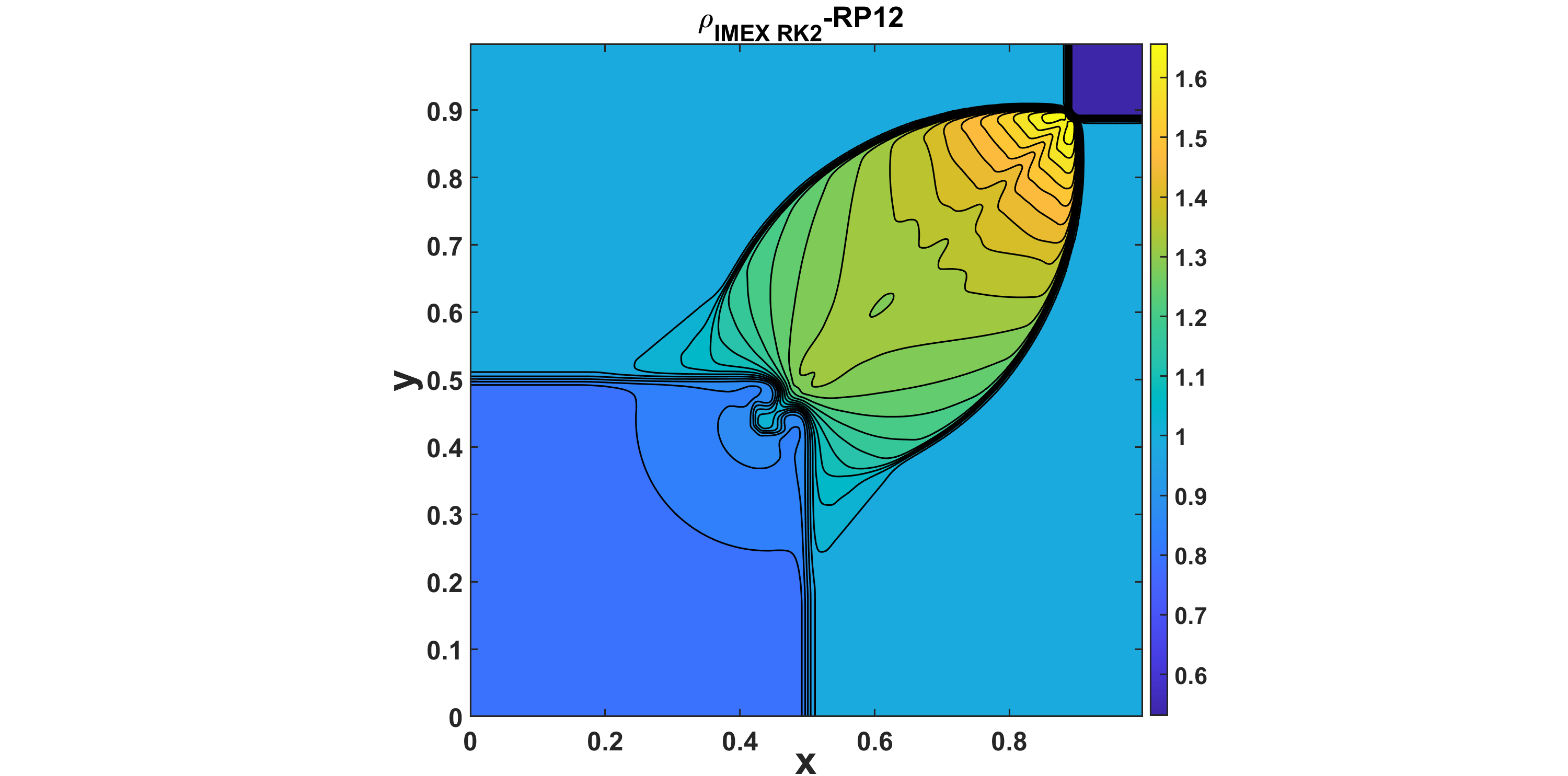}
        \caption{}
        \label{fig:rp12-imex}
    \end{subfigure}
    \caption{Comparison of the numerical density profiles $\rho$ for discontinuous initial data in Configuration 12 obtained using (a) the UCS2 scheme, (b) the proposed scheme with $\alpha=1/10$, and (c) the IMEX-RK2 scheme.}
    \label{/RP:euler:ep12}
\end{figure}

\begin{figure}[htb]
    \centering
    \begin{subfigure}[t]{0.325\textwidth}
        \centering        \includegraphics[width=\linewidth,height=5.5cm,trim=12.0cm 0cm 10.7cm 0cm,clip]{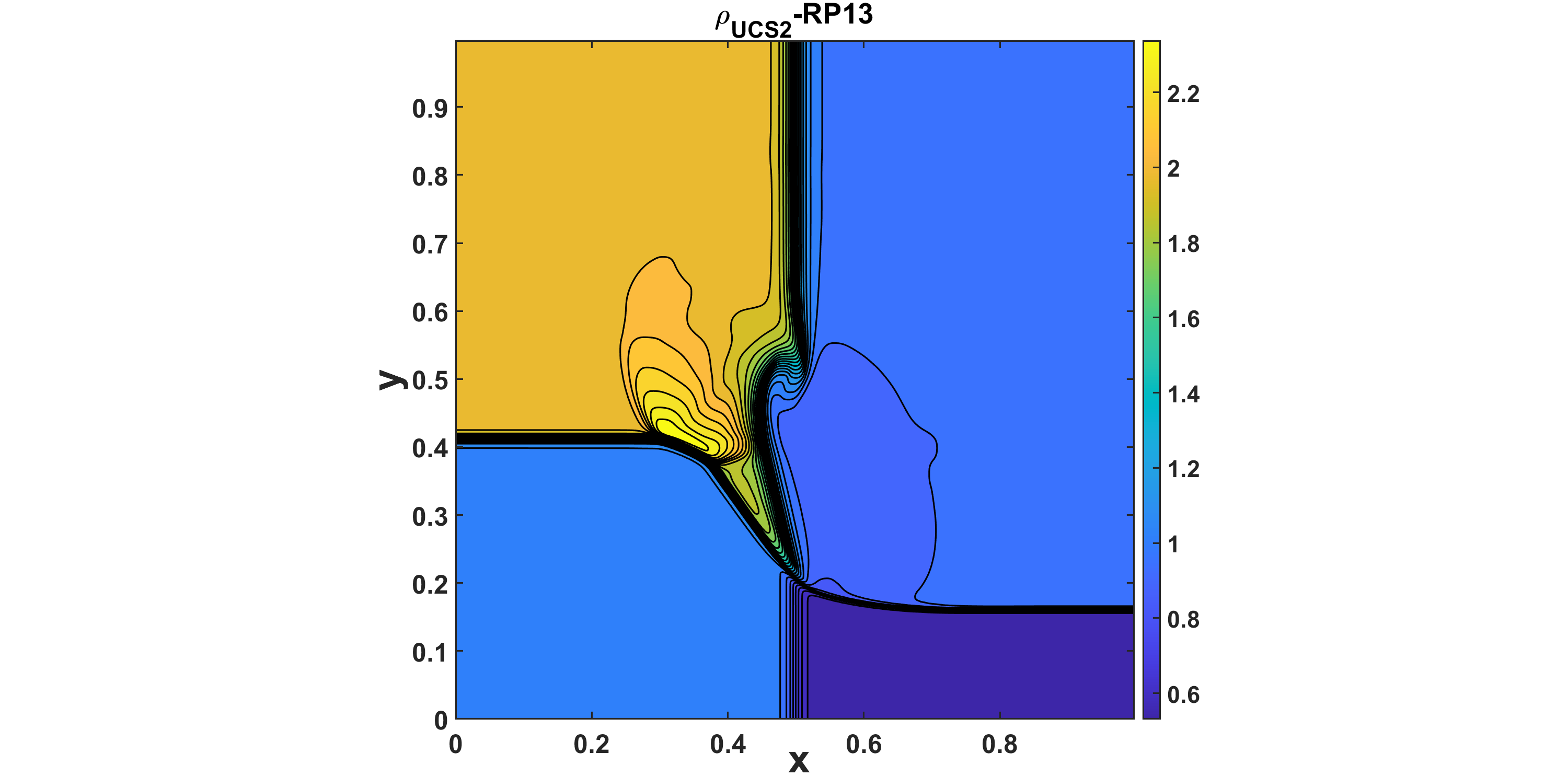}
        \caption{}
        \label{fig:rp13-ucs2}
    \end{subfigure}
    \hfill
    \begin{subfigure}[t]{0.325\textwidth}
        \centering        \includegraphics[width=\linewidth,height=5.5cm,trim=12.5cm 0cm 10.7cm 0cm,clip]{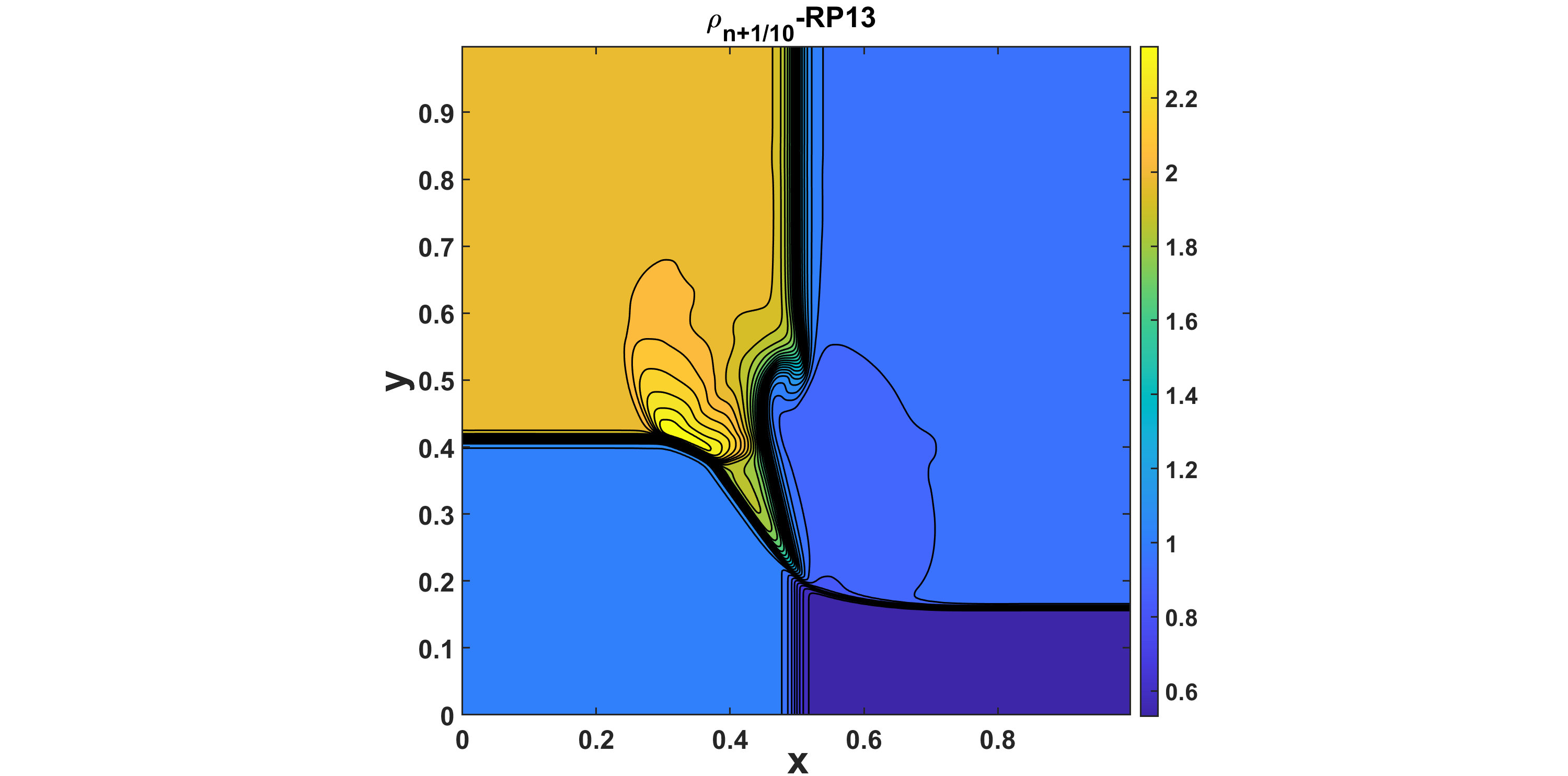}
        \caption{}
        \label{fig:rp13-proposed}
    \end{subfigure}
    \hfill
    \begin{subfigure}[t]{0.325\textwidth}
        \centering        \includegraphics[width=\linewidth,height=5.5cm,trim=12.0cm 0cm 10.7cm 0cm,clip]{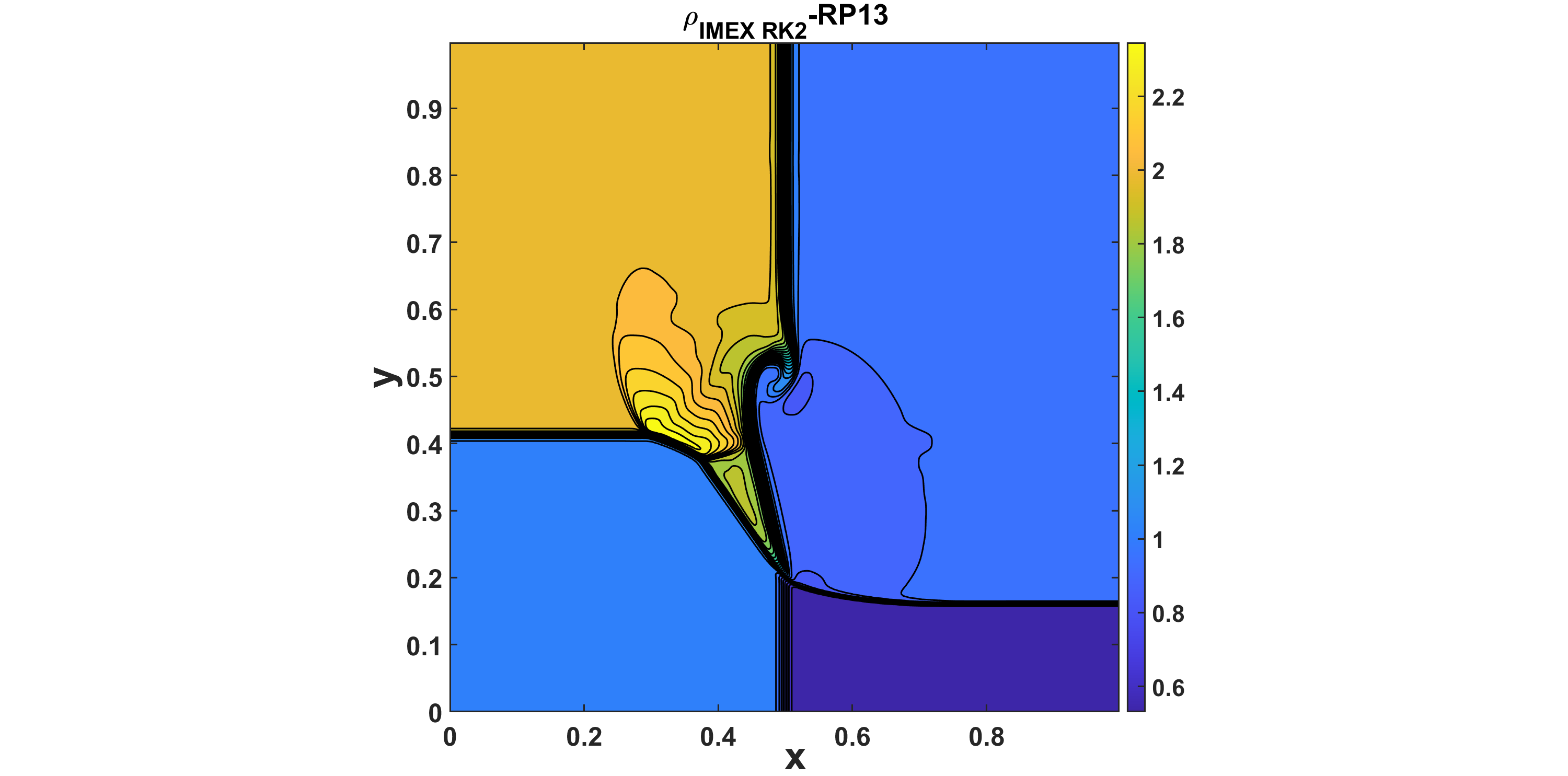}
        \caption{}
        \label{fig:rp13-imex}
    \end{subfigure}
    \caption{Comparison of the numerical density profiles $\rho$ for Configuration 13 obtained using (a) the UCS2 scheme, (b) the proposed scheme with $\alpha=1/10$, and (c) the IMEX-RK2 scheme.}
    \label{/RP:euler:ep13}
\end{figure}

\begin{figure}[htb]
    \begin{subfigure}[b]{0.325\textwidth}
         \centering      
         \includegraphics[scale=0.215, trim=12.0cm 0cm 10.7cm 0cm, clip]{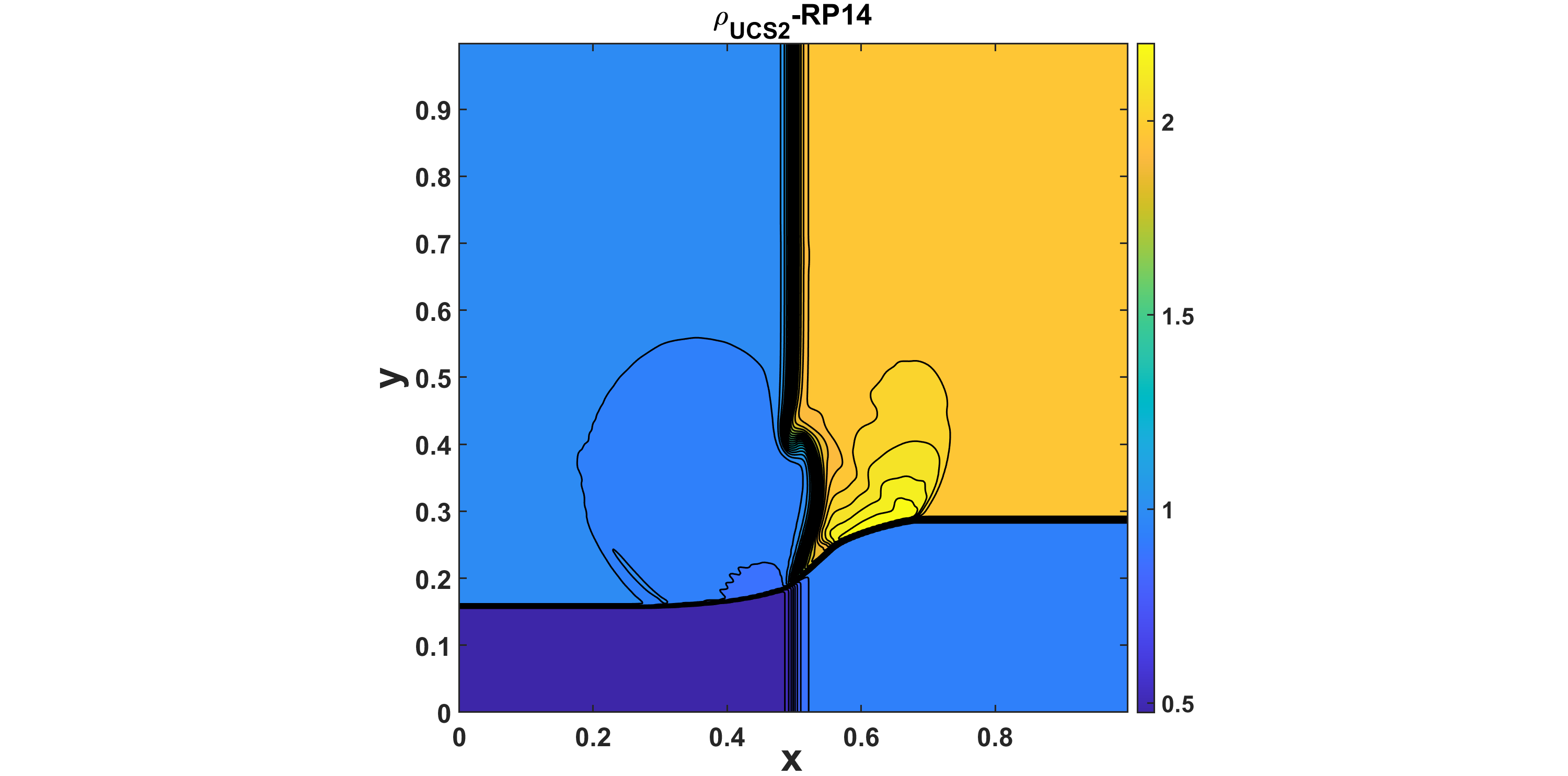}
          \caption{}
        \label{fig:rp14-ucs}
     \end{subfigure}
     \hfill
    \begin{subfigure}[b]{0.325\textwidth}
         \centering      
         \includegraphics[scale=0.215, trim=11.5cm 0cm 7.0cm 0cm, clip]{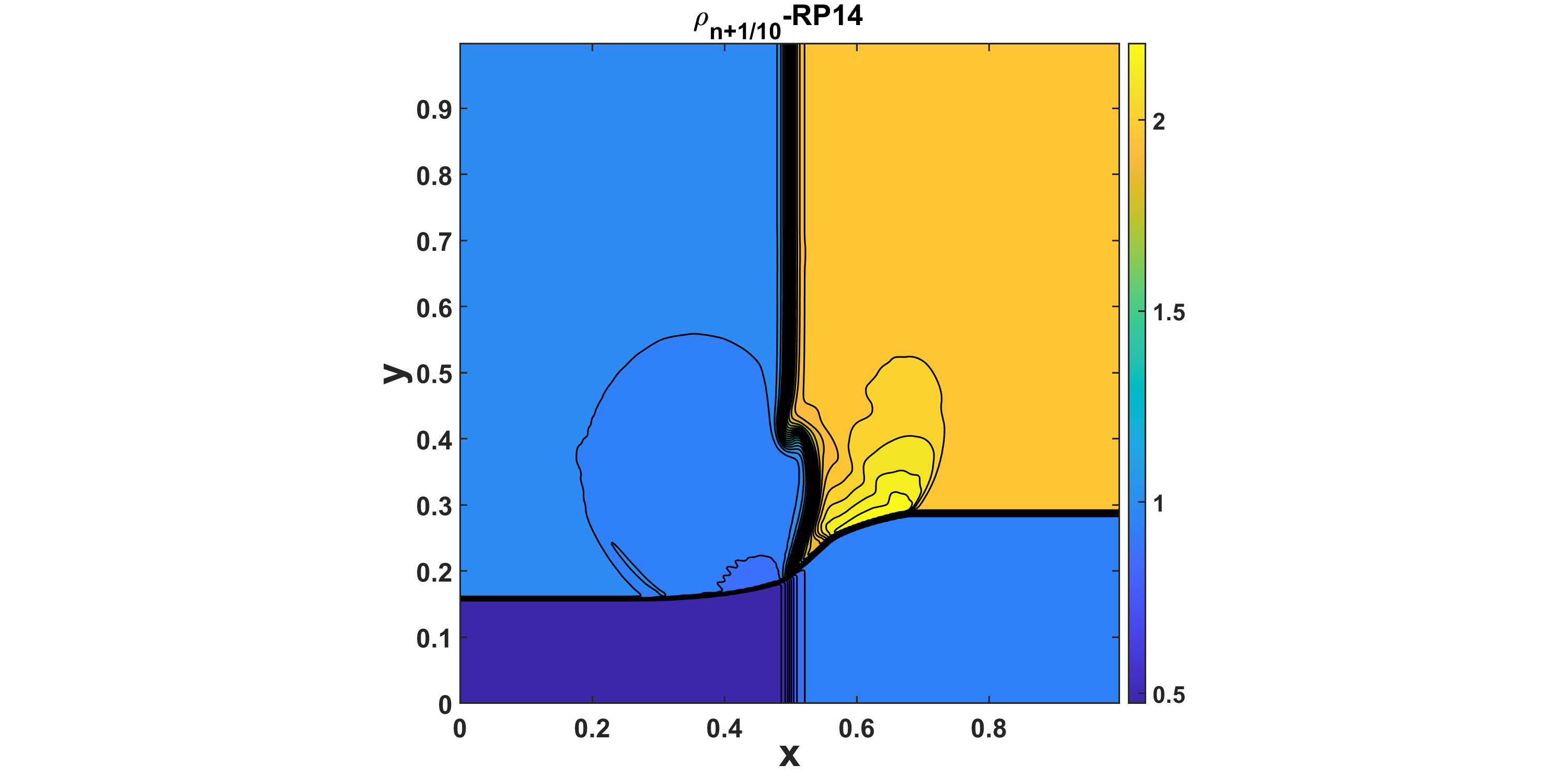}
          \caption{}
        \label{fig:rp14-n110}
     \end{subfigure}
      \hfill
    \begin{subfigure}[b]{0.325\textwidth}
         \centering      
         \includegraphics[scale=0.215, trim=12.0cm 0cm 10.7cm 0cm, clip]{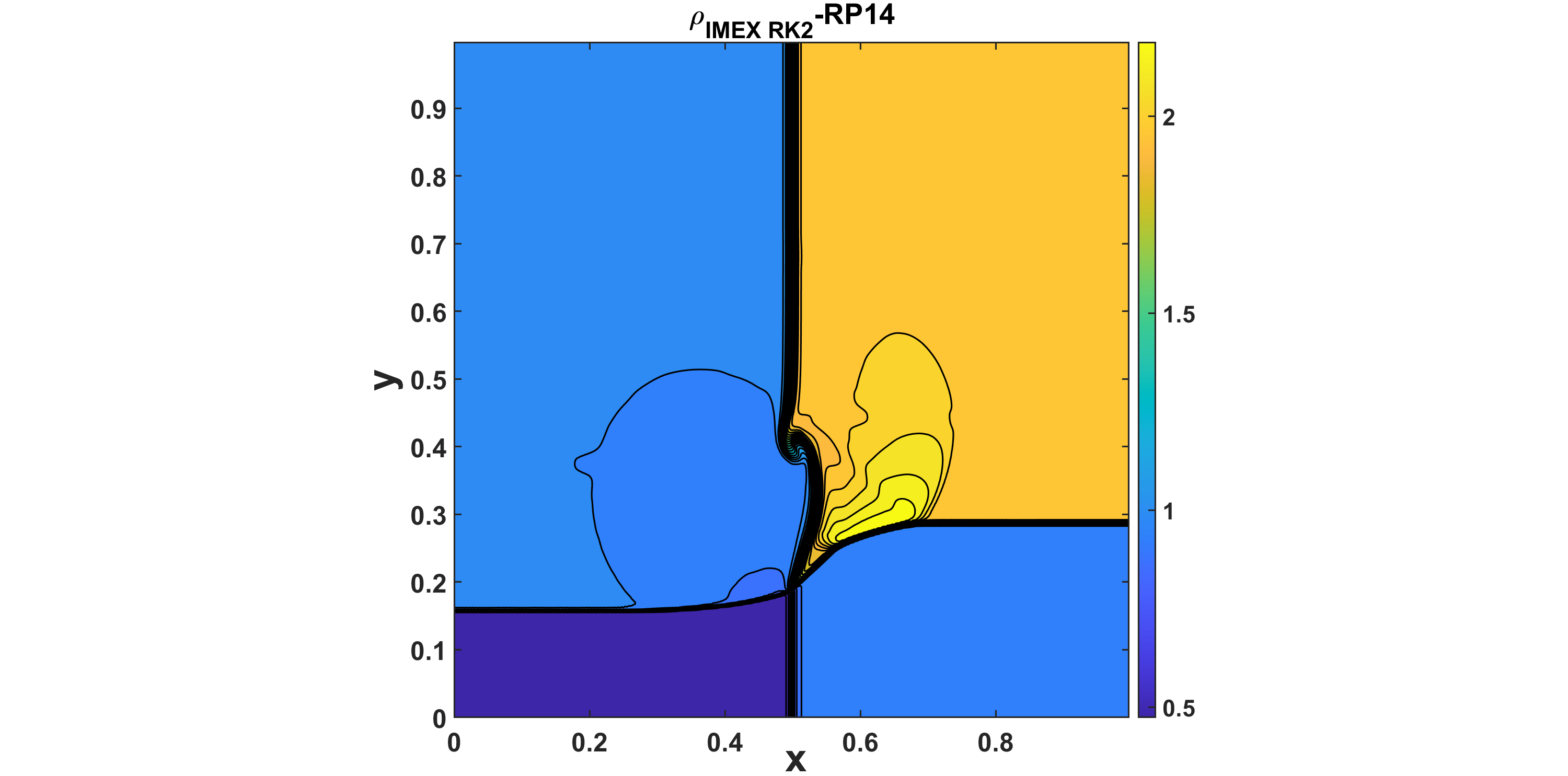}
          \caption{}
        \label{fig:rp14-imex}
     \end{subfigure}
     \caption{Comparison of the numerical density profiles $\rho$ for Configuration 14 obtained using (a) the UCS2 scheme, (b) the proposed scheme with $\alpha=1/10$, and (c) the IMEX-RK2 scheme.}
\label{/RP:euler:ep14}
\end{figure}

\begin{figure}[htb]
    \centering
    \begin{subfigure}[t]{0.325\textwidth}
        \centering        \includegraphics[width=\linewidth,height=5.5cm,trim=12.0cm 0cm 10.7cm 0cm,clip]{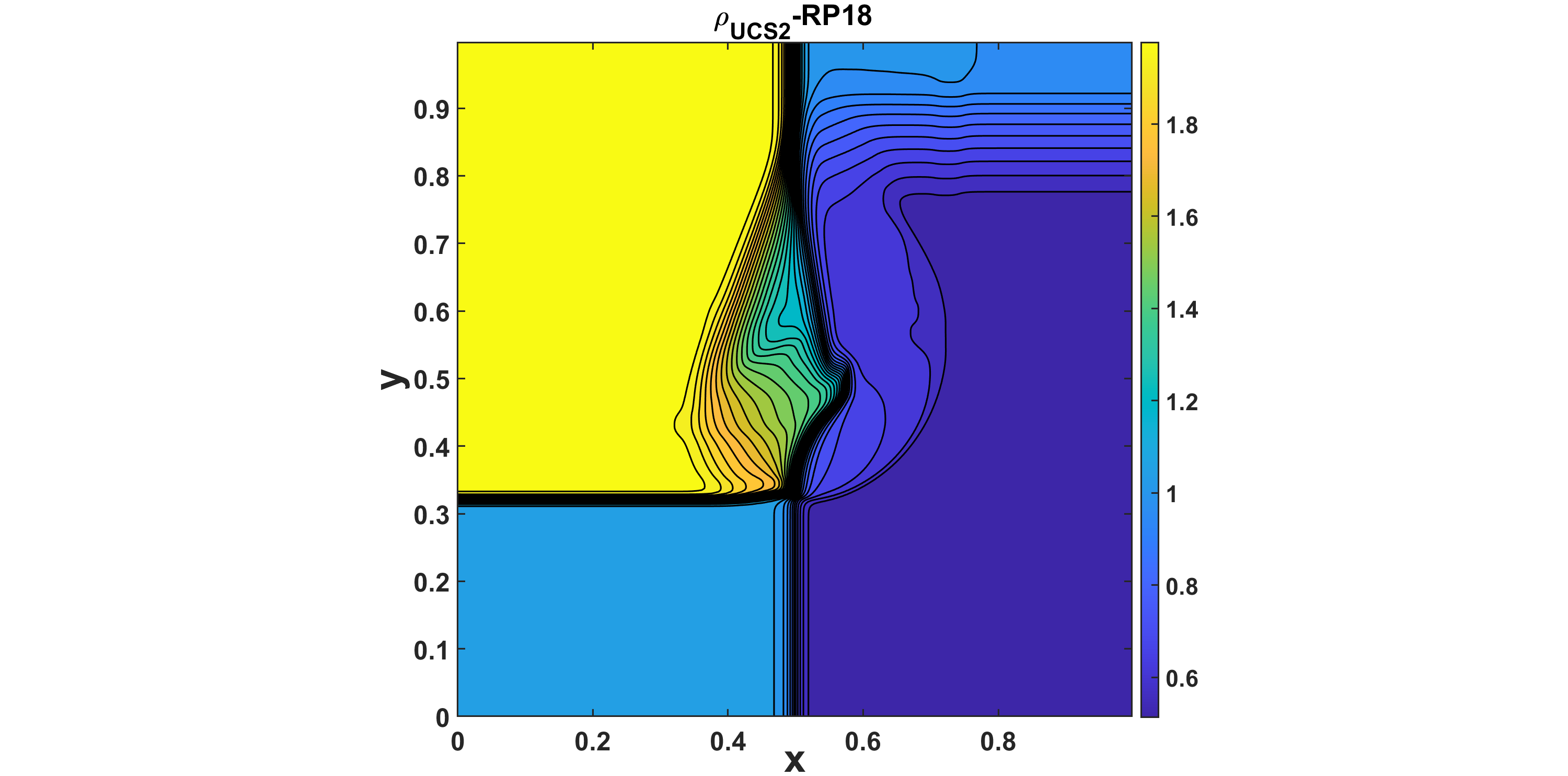}
        \caption{}
        \label{fig:rp18-ucs2}
    \end{subfigure}
    \hfill
    \begin{subfigure}[t]{0.325\textwidth}
        \centering        \includegraphics[width=\linewidth,height=5.5cm,trim=11.5cm 0cm 10.7cm 0cm,clip]{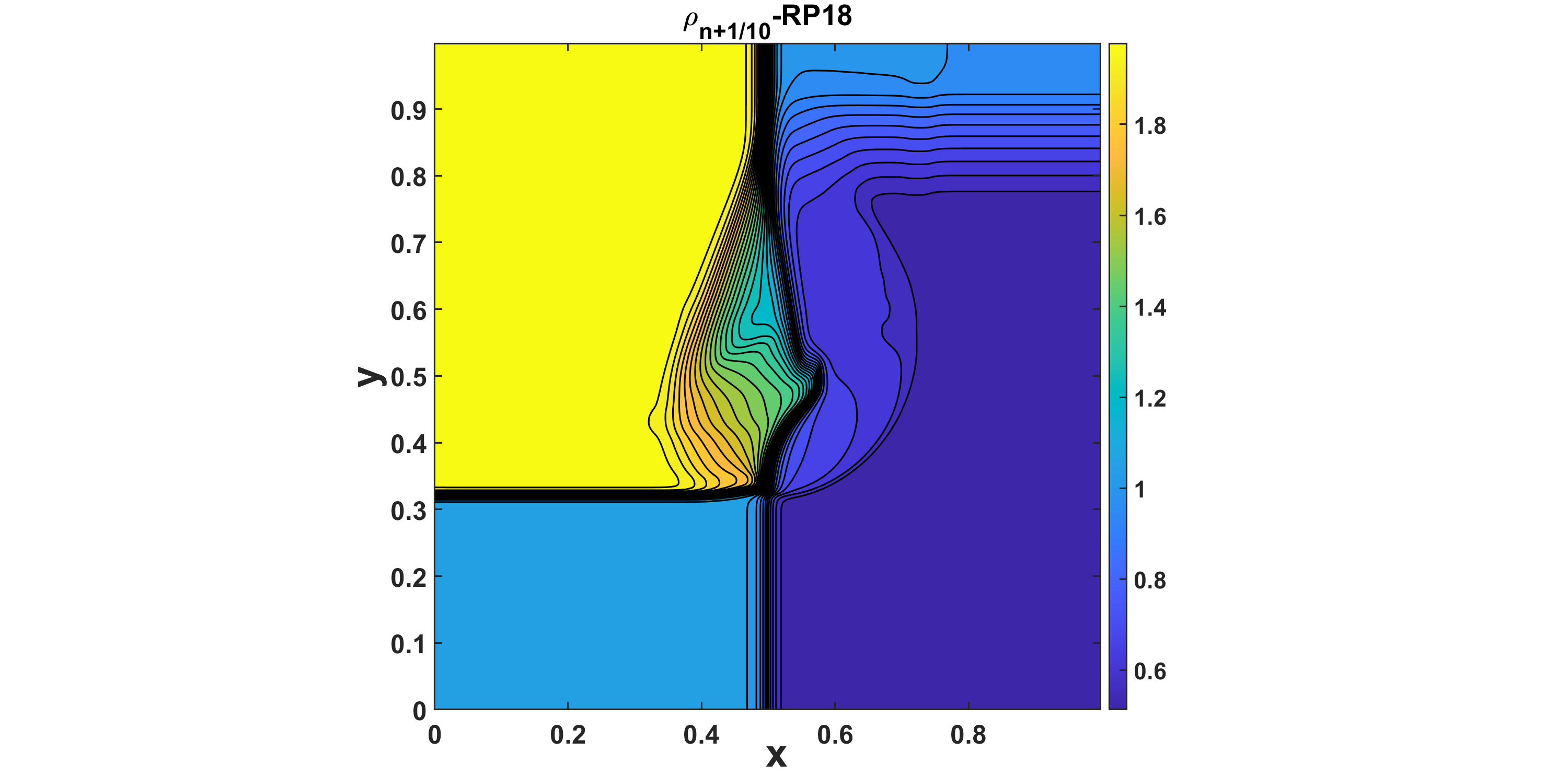}
        \caption{}
        \label{fig:rp18-proposed}
    \end{subfigure}
    \hfill
    \begin{subfigure}[t]{0.325\textwidth}
        \centering        \includegraphics[width=\linewidth,height=5.5cm,trim=12.0cm 0cm 10.7cm 0cm,clip]{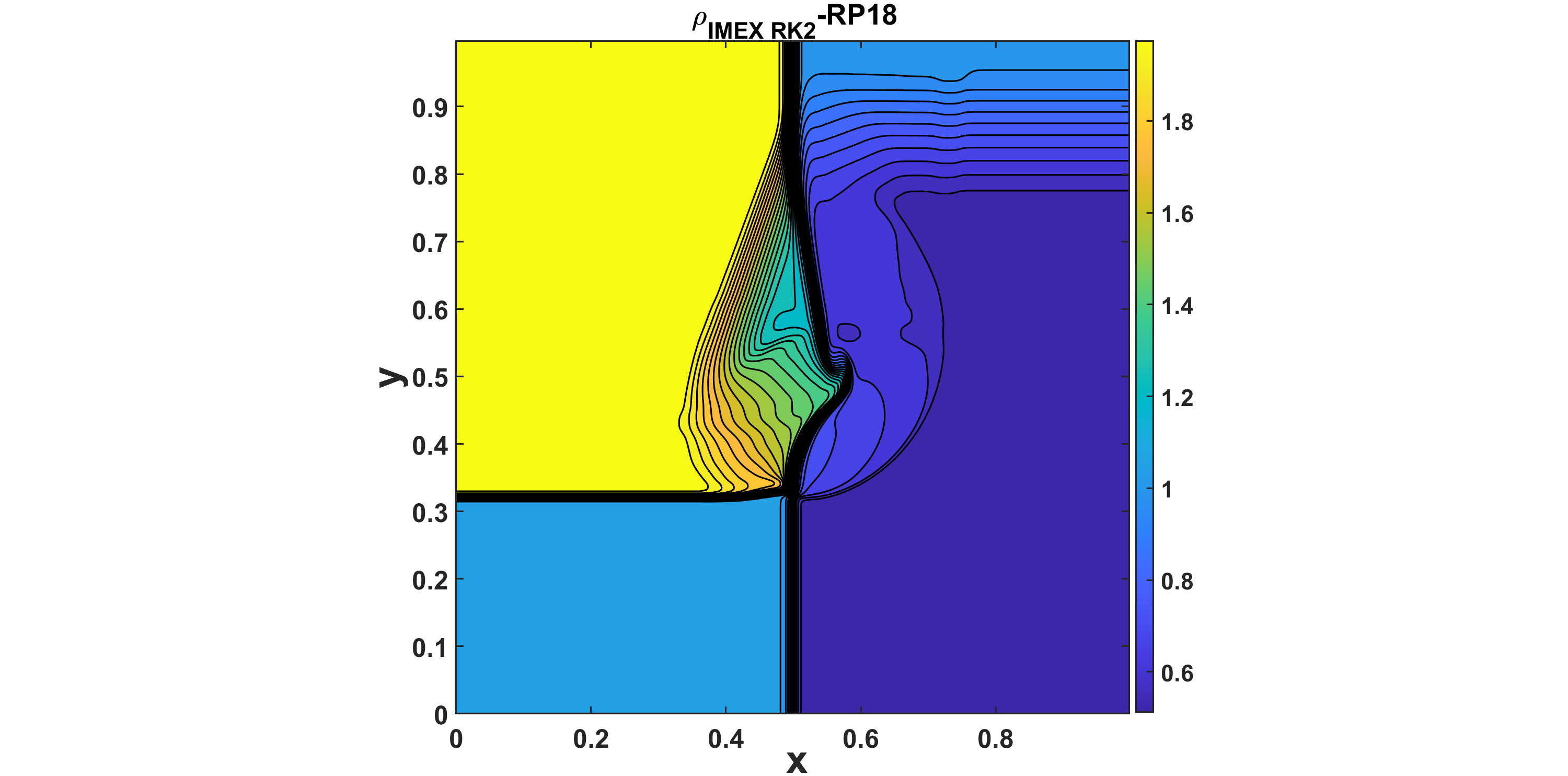}
        \caption{}
        \label{fig:rp18-imex}
    \end{subfigure}
    \caption{Comparison of the numerical density profiles $\rho$ for Configuration 18 obtained using (a) the UCS2 scheme, (b) the proposed scheme with $\alpha=1/10$, and (c) the IMEX-RK2 scheme.}
    \label{/RP:euler:ep18}
\end{figure}

Having assessed the accuracy of the proposed scheme for smooth, well-prepared initial data, we now investigate its performance for discontinuous initial data in two spatial dimensions. To this end, ten two-dimensional Riemann problems for the Euler equations are considered, covering a variety of wave structures, including shock waves, rarefaction waves, contact discontinuities, and their interactions. The corresponding density contour plots are presented in Figures~\ref{/RP:euler:ep2}-\ref{/RP:euler:ep18}. In the configuration labels, \(S\), \(R\), and \(C\) denote shock waves, rarefaction waves, and contact discontinuities, respectively. The proposed generalized scheme with \(\alpha=1/10\) is compared with UCS2, while IMEX-RK2 is used as a reference method for assessing the computed flow structures. The initial data and final time for each configuration are specified below.

All simulations are performed on a uniform Cartesian grid with \(N_x=N_y=400\) over the domain \([0,1]\times[0,1]\). The CFL number is fixed at \(0.9\), transmissive boundary conditions are imposed on all boundaries, and the relaxation parameter is set to \(\varepsilon=10^{-12}\). Unless otherwise specified, the relaxation variables are initialized on the equilibrium manifold according to
\[
V(x,y,0)=F_1(U(x,y,0)),\qquad
W(x,y,0)=F_2(U(x,y,0)).
\]
All contour plots are generated using $30$ contour levels. The final time is specified individually for each configuration.

\medskip
\noindent\textbf{Configuration 2 (\(R_1\)-\(R_2\)-\(R_3\)-\(R_4\)):}
The initial data for the primitive variables \((\rho,u,v,p)\) are prescribed as
\begin{equation}
(\rho,u,v,p)(x,y,0)=
\begin{cases}
(1.0,\,0,\,0,\,1.0), & x\geq 0.5,\; y\geq 0.5,\\
(0.5197,\,-0.7259,\,0,\,0.4), & x<0.5,\; y\geq 0.5,\\
(1.0,\,-0.7259,\,-0.7259,\,1.0), & x<0.5,\; y<0.5,\\
(0.5197,\,0,\,-0.7259,\,0.4), & x\geq 0.5,\; y<0.5.
\end{cases}
\end{equation}
The computation is carried out up to the final time \(T=0.2\).

This configuration involves four rarefaction waves arising from the discontinuities between the four initial states. The initial discontinuities are located along \(x=0.5\) and \(y=0.5\). The velocity components are negative in the left and lower halves, respectively, while the upper-right quadrant is initially at rest. The density and pressure also vary between the quadrants, providing a test of the numerical schemes' ability to capture expansion waves and represent the resulting multidimensional wave interactions. Figure~\ref{/RP:euler:ep2} compares the density contours obtained using the proposed scheme with \(\alpha=1/10\), UCS2, and IMEX-RK2.

\medskip
\noindent\textbf{Configuration 3 (\(S_1\)-\(S_2\)-\(S_3\)-\(S_4\)):} The initial data for the primitive variables \((\rho,u,v,p)\) are prescribed as
\begin{equation}
(\rho,u,v,p)(x,y,0)=
\begin{cases}
(1.5,\,0,\,0,\,1.5), & x\geq 0.8,\; y\geq 0.8,\\
(0.5323,\,1.2060,\,0,\,0.3), & x<0.8,\; y\geq 0.8,\\
(0.1380,\,1.2060,\,1.2060,\,0.0290), & x<0.8,\; y<0.8,\\
(0.5323,\,0,\,1.2060,\,0.3), & x\geq 0.8,\; y<0.8.
\end{cases}
\end{equation}
The computation is carried out up to the final time \(T=0.8\).

This configuration involves four shock waves, denoted by \(S_1\), \(S_2\), \(S_3\), and \(S_4\), arising from the discontinuities between the four initial states. The initial discontinuities are located at \(x=0.8\) and \(y=0.8\). The density and pressure in the lower-left quadrant are substantially smaller than those in the other three quadrants, while both velocity components are nonzero in this region. This problem thus serves as a stringent benchmark for assessing the capability of the numerical schemes to accurately resolve strong shock waves and their interactions in two-dimensional flows. Figure~\ref{/RP:euler:ep3} compares the density contours obtained using the proposed scheme with \(\alpha=1/10\), UCS2, and IMEX-RK2. In addition, the density distributions along the line \(y=0.8\) are examined to assess the resolution of the shock structures and the accuracy of the numerical solutions.

\noindent\textbf{Configuration 4 (\(S_1\)-\(S_2\)-\(S_3\)-\(S_4\)):} The initial data for the primitive variables \((\rho,u,v,p)\) are prescribed as
\begin{equation}
(\rho,u,v,p)(x,y,0)=
\begin{cases}
(1.1,\,0,\,0,\,1.1), & x\geq 0.5,\; y\geq 0.5,\\
(0.5026,\,0.8939,\,0,\,0.35), & x<0.5,\; y\geq 0.5,\\
(1.1,\,0.8939,\,0.8939,\,1.1), & x<0.5,\; y<0.5,\\
(0.5065,\,0,\,0.8939,\,0.35), & x\geq 0.5,\; y<0.5.
\end{cases}
\end{equation}
The computation is carried out up to the final time \(T=0.25\).

This configuration involves four shock waves, denoted by \(S_1\), \(S_2\), \(S_3\), and \(S_4\), arising from the discontinuities between the four initial states. The initial discontinuities are located along \(x=0.5\) and \(y=0.5\). The upper-right and lower-left quadrants have density and pressure values of \(1.1\), whereas the upper-left and lower-right quadrants have lower densities and pressures. The velocity components also vary between the quadrants, with both components nonzero in the lower-left quadrant. This problem provides a benchmark for evaluating the accuracy and robustness of the numerical schemes in capturing shock waves and their complex interactions in two-dimensional flows. The density contours obtained with the proposed scheme for \(\alpha=1/10\) are compared with the corresponding results from UCS2 and IMEX-RK2 in Figure~\ref{/RP:euler:ep4}.
\medskip

\noindent\textbf{Configuration 5 (\(C_1\)-\(C_2\)-\(C_3\)-\(C_4\)):} The initial data for the primitive variables \((\rho,u,v,p)\) are prescribed as
\begin{equation}
(\rho,u,v,p)(x,y,0)=
\begin{cases}
(1.0,\,-0.75,\,-0.5,\,1.0), & x\geq 0.5,\; y\geq 0.5,\\
(2.0,\,-0.75,\,0.5,\,1.0), & x<0.5,\; y\geq 0.5,\\
(1.0,\,0.75,\,0.5,\,1.0), & x<0.5,\; y<0.5,\\
(3.0,\,0.75,\,-0.5,\,1.0), & x\geq 0.5,\; y<0.5.
\end{cases}
\end{equation}
The computation is carried out up to the final time \(T=0.23\).

This configuration involves four contact waves, \(C_1\), \(C_2\), \(C_3\), and \(C_4\). The initial discontinuities are located along \(x=0.5\) and \(y=0.5\). The pressure is uniform across all four quadrants, \(p=1.0\), while the density and velocity components vary between the initial states. The horizontal velocity changes sign between the left and right halves of the domain, while the vertical velocity changes sign between the upper and lower halves. This test problem evaluates the capability of the numerical schemes to accurately capture contact discontinuities while maintaining sharp density interfaces with minimal numerical diffusion. The density contours produced by the proposed scheme with \(\alpha=1/10\), UCS2, and IMEX-RK2 are presented in Figure~\ref{/RP:euler:ep5}.\\
\noindent\textbf{Configuration 6 (\(C_1\)-\(C_2\)-\(C_3\)-\(C_4\)):} The initial data for the primitive variables \((\rho,u,v,p)\) are prescribed as
\begin{equation}
(\rho,u,v,p)(x,y,0)=
\begin{cases}
(1.0,\,0.75,\,-0.5,\,1.0), & x\geq 0.5,\; y\geq 0.5,\\
(2.0,\,0.75,\,0.5,\,1.0), & x<0.5,\; y\geq 0.5,\\
(1.0,\,-0.75,\,0.5,\,1.0), & x<0.5,\; y<0.5,\\
(3.0,\,-0.75,\,-0.5,\,1.0), & x\geq 0.5,\; y<0.5.
\end{cases}
\end{equation}
The computation is carried out up to the final time \(T=0.3\).
This configuration involves four contact waves, \(C_1\), \(C_2\), \(C_3\), and \(C_4\), arising from the discontinuities between the four initial states. The pressure is uniform throughout the domain, with \(p=1.0\), while the density varies between the quadrants. The velocity components have opposite signs in the left and right halves for the horizontal velocity and in the upper and lower halves for the vertical velocity. Compared with Configuration 5, the horizontal velocity signs are reversed, providing a complementary test of contact-wave propagation and interaction. Figure~\ref{/RP:euler:ep6} compares the density contours obtained using the proposed scheme with \(\alpha=1/10\), UCS2, and IMEX-RK2.\\

\noindent\textbf{Configuration 11 (\(C_1\)-\(C_2\)-\(S_1\)-\(S_2\)):} 
The initial data for the primitive variables \((\rho,u,v,p)\) are prescribed as
\begin{equation}
(\rho,u,v,p)(x,y,0)=
\begin{cases}
(1.0,\,0.1,\,0,\,1.0), & x\geq 0.5,\; y\geq 0.5,\\
(0.5313,\,0.8276,\,0,\,0.4), & x<0.5,\; y\geq 0.5,\\
(0.8,\,0.1,\,0,\,0.4), & x<0.5,\; y<0.5,\\
(0.5313,\,0.1,\,0.7276,\,0.4), & x\geq 0.5,\; y<0.5.
\end{cases}
\end{equation}
The computation is carried out up to the final time \(T=0.3\).

This configuration involves two contact waves, \(C_1\) and \(C_2\), and two shock waves, \(S_1\) and \(S_2\). The initial discontinuities are located along \(x=0.5\) and \(y=0.5\). The horizontal velocity is \(u=0.1\) in three quadrants, while the upper-left quadrant has a larger horizontal velocity of \(u=0.8276\). The vertical velocity is zero in the upper half of the domain and in the lower-left quadrant, but takes the value \(v=0.7276\) in the lower-right quadrant. The pressure is \(1.0\) in the upper-right quadrant and \(0.4\) elsewhere. These variations make the problem a mixed-wave test for assessing the schemes' ability to resolve contact discontinuities and capture shock structures. A comparative representation of the density contours obtained from the proposed scheme with \(\alpha=1/10\) and the UCS2 and IMEX-RK2 methods is provided in Figure~\ref{/RP:euler:ep11}.\\
\medskip
\noindent\textbf{Configuration 12 (\(C_1\)-\(C_2\)-\(S_1\)-\(S_2\)):} The initial data for the primitive variables \((\rho,u,v,p)\) are prescribed as
\[
(\rho,u,v,p)(x,y,0)=
\begin{cases}
(0.5313,\,0,\,0,\,0.4), & x\geq 0.5,\; y\geq 0.5,\\
(1.0,\,0.7276,\,0,\,1.0), & x<0.5,\; y\geq 0.5,\\
(0.8,\,0,\,0,\,1.0), & x<0.5,\; y<0.5,\\
(1.0,\,0,\,0.7276,\,1.0), & x\geq 0.5,\; y<0.5.
\end{cases}
\]
The computation is carried out up to the final time \(T=0.25\).
This configuration involves two contact waves, \(C_1\) and \(C_2\), and two shock waves, \(S_1\) and \(S_2\). The initial discontinuities are located along \(x=0.5\) and \(y=0.5\). The horizontal velocity is nonzero only in the upper-left quadrant, while the vertical velocity is nonzero only in the lower-right quadrant. The pressure is \(0.4\) in the upper-right quadrant and \(1.0\) in the remaining three quadrants, with density also varying among the initial states. This configuration tests the numerical schemes' ability to capture a mixed-wave pattern involving contact waves and shocks in a two-dimensional flow. Figure~\ref{/RP:euler:ep12} illustrates the density contour distributions generated by the proposed scheme with \(\alpha=1/10\), alongside the results from UCS2 and IMEX-RK2.\\

\medskip
\noindent\textbf{Configuration 13 (\(C_1\)-\(C_2\)-\(S_1\)-\(S_2\)):} The initial data for the primitive variables \((\rho,u,v,p)\) are prescribed as
\begin{equation}
(\rho,u,v,p)(x,y,0)=
\begin{cases}
(1.0,\,0,\,-0.3,\,1.0), & x\geq 0.5,\; y\geq 0.5,\\
(2.0,\,0,\,0.3,\,1.0), & x<0.5,\; y\geq 0.5,\\
(1.0625,\,0,\,0.8145,\,0.4), & x<0.5,\; y<0.5,\\
(0.5313,\,0,\,0.4276,\,0.4), & x\geq 0.5,\; y<0.5.
\end{cases}
\end{equation}
The computation is carried out up to the final time \(T=0.3\).
This configuration involves two contact discontinuities, \(C_1\) and \(C_2\), and two shock waves, \(S_1\) and \(S_2\). The initial discontinuities are located along \(x=0.5\) and \(y=0.5\). The horizontal velocity is zero in all four quadrants, whereas the transverse velocity, density, and pressure vary between the initial states. This problem therefore provides a test of the numerical schemes for a two-dimensional interaction involving both contact discontinuities and shock waves. Figure~\ref{/RP:euler:ep13} compares the density contours obtained using the proposed scheme with \(\alpha=1/10\), UCS2, and IMEX-RK2.\\
\medskip

\noindent\textbf{Configuration 14 (\(C_1\)-\(C_2\)-\(S_1\)-\(S_2\)):} The initial data for the primitive variables \((\rho,u,v,p)\) are prescribed as
\begin{equation}
(\rho,u,v,p)(x,y,0)=
\begin{cases}
(2.0,\,0,\,-0.5606,\,8.0), & x\geq 0.5,\; y\geq 0.5,\\
(1.0,\,0,\,-1.2172,\,8.0), & x<0.5,\; y\geq 0.5,\\
(0.4736,\,0,\,1.2172,\,2.6667), & x<0.5,\; y<0.5,\\
(0.9474,\,0,\,1.1606,\,2.6667), & x\geq 0.5,\; y<0.5.
\end{cases}
\end{equation}
The computation is carried out up to the final time \(T=0.1\).

This configuration involves two contact discontinuities, \(C_1\) and \(C_2\), and two shock waves, \(S_1\) and \(S_2\). The initial discontinuities are located along \(x=0.5\) and \(y=0.5\). The horizontal velocity is zero in all four quadrants, while the density and transverse velocity vary between the initial states. In addition, the pressure takes different values in the upper and lower halves of the domain. This problem tests the ability of the numerical schemes to resolve contact discontinuities and shock waves in a two-dimensional flow with distinct initial pressure levels. Figure~\ref{/RP:euler:ep14} compares the density contours obtained using the proposed scheme with \(\alpha=1/10\), UCS2, and IMEX-RK2.

\medskip
\noindent\textbf{Configuration 18 (\(C_1\)-\(C_2\)-\(S_1\)-\(R_1\)):} The initial data for the primitive variables \((\rho,u,v,p)\) are prescribed as
\[
(\rho,u,v,p)(x,y,0)=
\begin{cases}
(1.0,\,0,\,1.0,\,1.0), & x\geq 0.5,\; y\geq 0.5,\\
(2.0,\,0,\,-0.3,\,1.0), & x<0.5,\; y\geq 0.5,\\
(1.0625,\,0,\,0.2145,\,0.4), & x<0.5,\; y<0.5,\\
(0.5197,\,0,\,0.2741,\,0.4), & x\geq 0.5,\; y<0.5.
\end{cases}
\]
The computation is carried out up to the final time \(T=0.2\).
This configuration involves two contact discontinuities, \(C_1\) and \(C_2\), one shock wave, \(S_1\), and one rarefaction wave, \(R_1\). The initial discontinuities are located along \(x=0.5\) and \(y=0.5\). As in the preceding configurations, the horizontal velocity is zero in all four quadrants, whereas the transverse velocity and density vary among the initial states. The pressure is \(1.0\) in the upper half of the domain and \(0.4\) in the lower half. The coexistence of contact discontinuities, a shock wave, and a rarefaction wave makes this configuration a test of the schemes' ability to represent different wave structures within the same two-dimensional Riemann problem. Figure~\ref{/RP:euler:ep18} compares the density contours obtained using the proposed scheme with \(\alpha=1/10\), UCS2, and IMEX-RK2.

The density contour plots provide a basis for assessing the performance of the proposed generalized scheme compared to UCS2 and IMEX-RK2 in the ten configurations. The comparison focuses on the resolution of shock fronts and contact discontinuities, the representation of rarefaction regions, and the capture of interactions between different wave types. Particular attention is also given to numerical diffusion, the preservation of wave structures, and the occurrence of spurious oscillations.

\medskip
\noindent{\textbf{Computational Cost of 2D Riemann Problems:}} To complement the assessment of the computed flow structures, Table~\ref{tab:cpu} reports the CPU times of all three methods for the ten configurations. These measurements allow the computational efficiency of the proposed generalized scheme to be compared with that of UCS2, while IMEX-RK2 provides an additional reference for numerical solutions. Together, the density contours and CPU-time measurements provide a basis for assessing the wave-capturing performance and computational cost of the proposed scheme across a range of two-dimensional Riemann problems.

The differences in CPU time partly reflect the distinct spatial and source-term discretization strategies adopted by the three methods. In particular, the proposed generalized scheme and UCS2 employ genuinely multidimensional formulations and associated quadrature procedures, whereas the IMEX-RK2 implementation uses a dimension-by-dimension treatment, as expected, and it takes less computational cost. These differences contribute to the variation in computational effort reported in Table~\ref{tab:cpu}.
\begin{table}[htb]
\centering
\caption{CPU time (in seconds) for the two-dimensional Riemann problem with $N \times N=400 \times 400$, comparing the UCS2, IMEX-RK2, and proposed schemes with $\alpha=1/10$.}
\label{tab:cpu}
\begin{tabular}{cccc}
\hline
RP problems & UCS2 & n+1/10 & IMEX RK2\\
\hline
RP2 & 8801.7626 & 9271.8498 & 231.6196\\
RP3 & 35882.5056 & 34478.6762 & 936.2753\\
RP4 & 10949.4991 & 10957.7424 & 289.9436\\
RP5 & 10282.6925 & 10188.5300 & 267.9521\\
RP6 & 12905.2317 & 14217.17010 & 348.4659\\
RP11 & 13135.6505 & 13164.2404 & 347.9327\\
RP12 & 11555.0773 & 10631.3495 & 292.0908\\
RP13 & 13170.1841 & 13316.8375 & 348.4941\\
RP14 & 4289.0494 & 4493.4494 & 116.7729\\
RP18 & 8670.6700 & 8820.4881 & 232.9502 \\
\hline
\end{tabular}
\end{table}
\section{Conclusion}\label{sec:conclusion}
A generalized second-order central finite-volume scheme has been developed for one- and two-dimensional hyperbolic balance laws with stiff relaxation source terms. Based on an integral quadrature rule with undetermined coefficients, the formulation incorporates an intermediate-stage parameter \(\alpha\), whose admissible range is determined by the theoretical requirements of the scheme. Under appropriate coefficient relations, the proposed formulation recovers UCS2 as a particular case.

The theoretical analysis has established the uniform consistency of the scheme across relaxation regimes and its asymptotic-preserving (AP) property. The stability analysis has determined the admissible range of \(\alpha\) and the conditions under which the scheme is A-stable and L-stable. The corresponding stability regions have been examined for the stiff source contribution and compared with that of the classical explicit second-order Runge–Kutta (RK2) method for the explicit parameter choice. These results characterize the influence of the scheme parameters on its numerical properties and provide a theoretical basis for selecting \(\alpha\). Numerical experiments have demonstrated second-order convergence for smooth problems in both one and two dimensions. For the one-dimensional tests, the proposed scheme has been compared with UCS2, while exact solutions have been used for error assessment where available and highly resolved IMEX-RK2 solutions have served as reference solutions otherwise. The two-dimensional formulation has been evaluated using ten Riemann problems involving shocks, rarefaction waves, contact discontinuities, and complex wave interactions. Comparisons with UCS2 and IMEX-RK2 have been used to assess the resolution of these flow structures, while CPU-time comparisons with IMEX-RK2 have provided insight into the computational cost of the proposed scheme in the stiff regime. Overall, the proposed formulation provides a unified, parameter-dependent framework for central finite-volume schemes for solving hyperbolic balance laws with stiff relaxation.

% Overall, the proposed formulation provides a unified, parameter-dependent framework for central finite-volume schemes for hyperbolic balance laws with stiff relaxation. The theoretical analysis and numerical experiments demonstrate its consistency and stability properties, as well as its applicability across different relaxation regimes, subject to the derived restrictions on \(\alpha\).
% \medskip
% \noindent{\textbf{Data availability:}} Enquiries about data availability should be directed to the authors.

\medskip
\noindent{\textbf{Declarations:}}\\
\textbf{Conflict of interest:} The authors declare that they have no known competing financial interests or personal
relationships that could have appeared to influence the work reported in this paper. Codes will be available from the corresponding author upon reasonable request.

\medskip
\noindent{\textbf{Funding:}} The authors did not receive any funding support for conducting this study.

\bibliographystyle{ieeetr}
\bibliography{biblio}
%\clearpage
\appendix
\section{}\label{appendix}
The numerical formulations of the UCS2 central scheme and the IMEX-RK2 scheme employed in this study are presented in this appendix. The corresponding one-dimensional and two-dimensional formulations are described separately.
\subsection{UCS2 central scheme}\label{sub:ucs}
The UCS2 central scheme \cite{Liotta} is considered as a reference numerical method in this study. Its one-dimensional and two-dimensional formulations are described below.
\subsubsection{One-dimensional UCS2 scheme}\label{sub:ucs1d}
The one-dimensional UCS2 predictor-corrector formulation consists of two successive stages. The predictor stage is given by
\begin{align}
\bar{q}_j^{\,n+\frac12}
&=
\bar{q}_j^n
+
\frac{\Delta t}{2}
\left(\frac{1}{\varepsilon}
\mathrm{g}(\bar{q}_j^{\,n+\frac12})
-
\frac{\mathrm{f}'_j}{\Delta x}
\right),\\
\bar{q}_j^{\,n+\frac13}
&=
\bar{q}_j^n
+
\frac{\Delta t}{3}
\left(\frac{1}{\varepsilon}
\mathrm{g}(\bar{q}_j^{\,n+\frac13})
-
\frac{\mathrm{f}'_j}{\Delta x}
\right).
\end{align}
which provides an intermediate approximation to the solution. The corrector stage is then used to update this predicted solution and is expressed as
\begin{align}
\bar{q}_{j+\frac12}^{\,n+1}
&=
\frac{1}{2}
(\bar{q}_j^n+\bar{q}_{j+1}^n)
+
\frac{1}{8}
(\bar{q}'_j-\bar{q}'_{j+1})
-
\lambda
\Bigl[
\mathrm{f}(\bar{q}_{j+1}^{\,n+\frac12})
-
\mathrm{f}(\bar{q}_j^{\,n+\frac12})
\Bigr]
\notag\\
&\quad
+
\frac{\Delta t}{\varepsilon}
\left(
\frac{3}{8}\mathrm{g}(\bar{q}_j^{\,n+\frac13})
+
\frac{3}{8}\mathrm{g}(\bar{q}_{j+1}^{\,n+\frac13})
+
\frac{1}{4}\mathrm{g}(\bar{q}_{j+\frac12}^{\,n+1})
\right).
\end{align}
In these expressions, $\lambda$ denotes the ratio of the temporal and spatial step sizes, i.e., $\lambda=\Delta t/\Delta x$. The quantities $\bar{q}'_j$ and $\mathrm{f}'_j$ are evaluated using the MinMod slope limiter. Within the framework of generalized central schemes, the UCS2 formulation is obtained by setting the parameter $\alpha$ to $1/3$.\\
\subsubsection{Two-dimensional UCS2 scheme}\label{sub:ucs2d}
The two-dimensional formulation of the UCS2 scheme is given by\\
\noindent\textbf{Predictor step:}
\label{UCS2:scheme2d} 
Solve for $\bar{q}_{j,k}^{\,n+\frac12}$ and
$\bar{q}_{j,k}^{\,n+\frac13}$ from
\begin{equation}
\bar{q}_{j,k}^{\,n+\frac12} = \bar{q}_{j,k}^{\,n} +\dfrac{\Delta t}{2} \left[\dfrac{1}{\varepsilon}\,\mathrm{g}\!\left(\bar{q}_{j,k}^{\,n+\frac12}\right) -\dfrac{(f_1)^{\dagger}_{j,k}}{\Delta x} -\dfrac{(f_2)^{\ast}_{j,k}}{\Delta y}\right].
\end{equation}
\begin{equation}
\bar{q}_{j,k}^{\,n+\frac13} = \bar{q}_{j,k}^{\,n} +\dfrac{\Delta t}{3}
\left[ \dfrac{1}{\varepsilon}\,\mathrm{g}\!\left(\bar{q}_{j,k}^{\,n+\frac13}\right)- \dfrac{(f_1)^{\dagger}_{j,k}}{\Delta x}-\dfrac{(f_2)^{\ast}_{j,k}}{\Delta y}\right].
\end{equation}
\noindent\textbf{Corrector step:} Then update to $\bar{q}_{j+\frac12,k+\frac12}^{\,n+1}$ by
\begin{equation}
\begin{aligned}
\bar{q}_{j+\frac12,k+\frac12}^{\,n+1}
&=\left\langle \dfrac{1}{4} \left(\bar{q}_{j,\cdot}^{\,n}+\bar{q}_{j+1,\cdot}^{\,n}\right)+\dfrac{1}{8} \left(\bar{q}^{\dagger}_{j,\cdot} -\bar{q}^{\dagger}_{j+1,\cdot} \right)-\lambda\left[ \mathrm{f}_1\!\left(\bar{q}_{j+1,\cdot}^{\,n+\frac12}\right) -\mathrm{f}_1\!\left(\bar{q}_{j,\cdot}^{\,n+\frac12}\right) \right]\right.\\
&\qquad\left. +\frac{\Delta t}{\varepsilon} \left[\dfrac{3}{16} \mathrm{g}\!\left(\bar{q}_{j,\cdot}^{\,n+\frac13}\right) +\dfrac{3}{16} \mathrm{g}\!\left(\bar{q}_{j+1,\cdot}^{\,n+\frac13}\right) \right]\right\rangle_{k+\frac12}\\
&\quad+ \left\langle \dfrac{1}{4} \left(\bar{q}_{\cdot,k}^{\,n}+\bar{q}_{\cdot,k+1}^{\,n}\right) +\dfrac{1}{8}\left(\bar{q}^{\ast}_{\cdot,k} -\bar{q}^{\ast}_{\cdot,k+1}\right)
-\mu \left[\mathrm{f}_2\!\left(\bar{q}_{\cdot,k+1}^{\,n+\frac12}\right)
-\mathrm{f}_2\!\left(\bar{q}_{\cdot,k}^{\,n+\frac12}\right) \right]\right.\\
&\qquad\left. +\frac{\Delta t}{\varepsilon} \left[ \dfrac{3}{16}
\mathrm{g}\!\left(\bar{q}_{\cdot,k}^{\,n+\frac13}\right)+ \dfrac{3}{16}
\mathrm{g}\!\left(\bar{q}_{\cdot,k+1}^{\,n+\frac13}\right) \right] \right\rangle_{j+\frac12} + \dfrac{\Delta t}{4 \varepsilon}\mathrm{g}\!\left(
\bar{q}_{j+\frac12,k+\frac12}^{\,n+1} \right).
\end{aligned}
\end{equation}
\noindent\textbf{Additional definitions.}
We define the mesh ratios
$\lambda = \frac{\Delta t}{\Delta x}$ and $\mu = \frac{\Delta t}{\Delta y}$.
The quantities $q^{\dagger}_{j,k}$ and $q^{\ast}_{j,k}$ denote the discrete slopes in the $x$- and $y$-directions, respectively, approximating
\[
q^{\dagger} \approx \Delta x\, \partial_x q,
\qquad
q^{\ast} \approx \Delta y\, \partial_y q.
\]
Similarly, $(\mathrm{f}_1)^{\dagger}_{j,k}$ and $(\mathrm{f}_2)^{\ast}_{j,k}$ approximate the derivatives of the flux functions.

To prevent spurious oscillations near discontinuities, all slopes are computed using the MinMod limiter \eqref{eq:minmod}.
\[ 
\langle q_{\cdot,j} \rangle_{j+\frac12} := \frac12 (q_{j,k}+q_{j+1,k}), \quad \langle q_{j,\cdot} \rangle_{k+\frac12} := \frac12 (q_{j,k}+q_{j,k+1}), 
\] 
denotes the staggered average in the $x$- and $y$-directions respectively.
\subsection{IMEX Runge-Kutta Scheme.}\label{sub:imex}
The IMEX-RK2 scheme \cite{PareschiRusso} is used as a reference time-integration method in this study. Its one-dimensional and two-dimensional formulations are presented below.
\subsubsection{One-dimensional IMEX-RK2 scheme}\label{sub:imex1d}
% Let the spatial domain be $\Omega=[a,b]\subset\mathbb{R}$, which is discretized into $N$ uniform cells
% \[
% I_j=[x_{j-\frac12},x_{j+\frac12}], \qquad j=1,\dots,N,
% \]
% with cell size $\Delta x=(b-a)/N$ and cell centers $x_j=\tfrac12(x_{j-\frac12}+x_{j+\frac12})$. The time interval $[0,T]$ is discretized into time levels $t^n=n\,\Delta t$, where $\Delta t$ denotes the time step.

% The spatial discretization is based on a conservative finite-volume formulation for the balance law~\eqref{eq:hyperbolic}. In particular, the interface fluxes are computed using a second-order piecewise linear reconstruction with MinMod slope limiting, and the resulting semi-discrete operator is integrated in time by the IMEX scheme.
The IMEX-RK2 scheme is characterized by the following parameters:
\begin{equation*}
\gamma = 1 - \dfrac{1}{\sqrt{2}},
\qquad
c = \dfrac{1}{2\gamma},
\end{equation*}
and by the explicit and implicit Butcher tableaux
\[
\begin{array}{c|cc}
0 & 0 & 0 \\
c & c & 0 \\
\hline
& 1-\gamma & \gamma
\end{array}
\qquad
\begin{array}{c|cc}
\gamma & \gamma & 0 \\
1-\gamma & 1-\gamma & \gamma \\
\hline
& 1-\gamma & \gamma
\end{array}.
\]

Let $\bar{q}_j^n$ denote the numerical approximation of the cell average over $I_j$ at grid point $x_j$ and time $t^n$. The first stage is given by
\begin{equation}
\bar{q}_j^{(1)} = \bar{q}_j^n - \gamma
\frac{\Delta t}{\Delta x} \bigl( \mathrm{f}_{j+\frac12}(\bar{q}^n) - \mathrm{f}_{j-\frac12}(\bar{q}^n) \bigr) + \frac{\gamma\Delta t}{\varepsilon}\,
\mathrm{g}\!\left(\bar{q}_j^{(1)}\right),
\end{equation}
where the flux contribution is evaluated explicitly and the source term implicitly.

Following the Boscarino reformulation, an auxiliary explicit state is introduced,
\begin{equation}
\bar{q}_j^{(e)} = \left(1-\frac{c}{\gamma}\right)q_j^n + \frac{c}{\gamma}\bar{q}_j^{(1)},
\end{equation}
which is used to evaluate the explicit flux in the second stage.

The second stage reads
\begin{equation}
\bar{q}_j^{n+1} = \bar{q}_j^{\ast} - \gamma \frac{\Delta t}{\Delta x}
\bigl( \mathrm{f}_{j+\frac12}(\bar{q}^{(e)}) - \mathrm{f}_{j-\frac12}(\bar{q}^{(e)}) \bigr)
+\frac{\gamma\Delta t}{\varepsilon}\, \mathrm{g}\!\left(\bar{q}_j^{n+1}\right),
\end{equation}
where
\begin{equation}
\bar{q}_j^{\ast} = \left(1-\frac{1-\gamma}{\gamma}\right)\bar{q}_j^n
+ \frac{1-\gamma}{\gamma}\bar{q}_j^{(1)}.
\end{equation}
\subsubsection{Two-dimensional IMEX-RK2 scheme}\label{sub:imex2d}
% Let the two-dimensional spatial domain be
% $\Omega=[a,b]\times[c,d]\subset\mathbb{R}^2$,
% which is discretized into $N_x\times N_y$ uniform Cartesian cells
% \[
% I_{j,k}
% =
% [x_{j-\frac12},x_{j+\frac12}]
% \times
% [y_{k-\frac12},y_{k+\frac12}],
% \qquad
% j=1,\ldots,N_x,\quad
% k=1,\ldots,N_y,
% \]
% with cell sizes $\Delta x=\dfrac{b-a}{N_x}$,
% and $\Delta y=\dfrac{d-c}{N_y}$, and cell centers
% $x_j=\dfrac{1}{2}\left(x_{j-\frac12}+x_{j+\frac12}\right)$,
% and $y_k=\dfrac{1}{2}\left(y_{k-\frac12}+y_{k+\frac12}\right)$.

% The time interval $[0,T]$ is discretized as
% $t^n=n\Delta t$, where $\Delta t$ denotes the time step. Let $\bar{q}_{j,k}^n$ denote the
% cell average of $q$ over $I_{j,k}$ at time $t^n$.
% The two-dimensional finite-volume spatial operator is defined by
% \begin{equation}
% \mathcal{L}_{j,k}(q)
% =
% -\dfrac{
% \mathrm{f}_{1,j+\frac12,k}(q)-\mathrm{f}_{1,j-\frac12,k}(q)
% }{\Delta x}
% -\dfrac{
% \mathrm{f}_{2,j,k+\frac12}(q)- \mathrm{f}_{2,j,k-\frac12}(q)
% }{\Delta y},
% \label{eq:2d_spatial_operator}
% \end{equation}
% where $\mathrm{f}_1$ and $\mathrm{f}_2$ denote the numerical fluxes in the $x$- and $y$-directions, respectively. The same second-order piecewise linear reconstruction with MinMod slope limiting is employed in both spatial directions.\\
% The same IMEX RK2 coefficients
% \begin{equation*}
% \gamma=1-\frac{1}{\sqrt{2}},
% \qquad
% c=\frac{1}{2\gamma},
% \label{eq:2d_imex_coefficients}
% \end{equation*}
% as in the one-dimensional formulation are used. .
For the two-dimensional hyperbolic relaxation system \eqref{2d:EQuation}, we employ the following finite-volume IMEX-RK2 discretization. 
The two-dimensional finite-volume spatial operator is defined by
\begin{equation}
\mathcal{L}_{j,k}(q)
=
-\dfrac{
\mathrm{f}_{1,j+\frac12,k}(q)-\mathrm{f}_{1,j-\frac12,k}(q)
}{\Delta x}
-\dfrac{
\mathrm{f}_{2,j,k+\frac12}(q)- \mathrm{f}_{2,j,k-\frac12}(q)
}{\Delta y},
\label{eq:2d_spatial_operator}
\end{equation}
where $\mathrm{f}_1$ and $\mathrm{f}_2$ denote the numerical fluxes in the $x$- and $y$-directions, respectively. A second-order piecewise linear reconstruction with MinMod slope limiting is employed in both spatial directions. Furthermore, the IMEX-RK2 coefficients are retained from the one-dimensional formulation.

The two-stage IMEX-RK2 scheme is constructed as follows. The first stage is given by
\begin{equation}
\bar{q}_{j,k}^{(1)}
=
\bar{q}_{j,k}^{n}
+
\gamma \Delta t\,\mathcal{L}_{j,k}(\bar{q}^{n})
+
\frac{\gamma\Delta t}{\varepsilon}\,
\mathrm{g} \left(\bar{q}_{j,k}^{(1)}\right).
\label{eq:2d_first_stage}
\end{equation}
Following the Boscarino reformulation, the auxiliary state
\begin{equation}
\bar{q}_{j,k}^{(e)} = \left(1-\frac{c}{\gamma}\right)\bar{q}_{j,k}^{n}
+\frac{c}{\gamma}\bar{q}_{j,k}^{(1)},
\label{eq:2d_explicit_state}
\end{equation}
is introduced for the explicit flux evaluation in the second stage. We
also define
\begin{equation}
\bar{q}_{j,k}^{*} = \left(1-\frac{1-\gamma}{\gamma}\right)\bar{q}_{j,k}^{n} + \frac{1-\gamma}{\gamma}\bar{q}_{j,k}^{(1)}.
\label{eq:2d_star_state}
\end{equation}
The second stage is then written as
\begin{equation}
\bar{q}_{j,k}^{n+1}=\bar{q}_{j,k}^{*}+\gamma\Delta t\,\mathcal{L}_{j,k}(\bar{q}^{(e)})+\frac{\gamma\Delta t}{\varepsilon}\, \mathrm{g} \left(\bar{q}_{j,k}^{n+1}\right).
\label{eq:2d_second_stage}
\end{equation}
The above two-dimensional IMEX-RK2 formulation is used for the numerical approximation of the system \eqref{2d:EQuation}.

\end{document}